\documentstyle[epsf]{article}

\newcommand{\ee}{\end{equation}}
\newcommand{\be}{\begin{equation}}
\def\adots{\cdot^{\displaystyle{~\cdot}^{\displaystyle{~\cdot}} }}
\def\bin#1#2{\pmatrix{ #1 \cr #2 \cr}}  

\begin{document}

\centerline {\bf\huge Stokes Matrices and Monodromy of the  }
\vskip 0.3 cm 
\centerline {\bf\huge Quantum Cohomology of Projective Spaces}

\vskip  0.4 cm
\centerline{\bf\large Davide Guzzetti}

\vskip 0.3 cm
\centerline{ International School for Advanced Studies (SISSA - ISAS)}
\centerline {Via Beirut, 2--4, ~~~  34014
Trieste, 
Italy.}
\centerline{ E-mail: guzzetti@sissa.it, Fax:  +39 040-3787-528}
\vskip 0.4 cm
\centerline{ August 1998}
\begin{abstract}   
In this paper we compute  Stokes matrices and monodromy for the quantum
cohomology of  projective spaces. 
This problem can be
formulated in a ``classical'' framework, as the problem of computation
of Stokes matrices and monodromy of (systems of) differential
equations with regular and irregular singularities.  We prove 
that  the Stokes' matrix of the quantum cohomology coincides  with the Gram 
matrix in the theory of 
derived categories of coherent sheaves.   We also study   the 
monodromy group of the quantum cohomology and we show that it is related 
 to hyperbolic triangular groups. 
\end{abstract}

\section{ Introduction}
\vskip 0.2 cm

 In this paper we  compute Stokes' matrices  and monodromy group for the 
Frobenius manifold given by the quantum cohomology of the projective space 
${\bf CP}^d$. Our main motivation is to study the  links between 
quantum cohomology  and the theory of coherent sheaves. 
 
 Stokes matrices first appeared in the theory of WDVV equations of associativity in the 
paper \cite{Dub4} by B. Dubrovin. WDVV equations were formulated in a geometrical setting: 
the theory of Frobenius manifolds. From then on, the notion of Frobenius manifold has been 
largely studied in many papers, of which we cite \cite{Dub1}, \cite{Dub2}, \cite{Dub3}.  

 The Stokes' matrix is a part of the {\it monodromy data} for a
semisimple Frobenius 
manifold. Monodromy data can serve as natural moduli of  semisimple 
Frobenius manifolds. More precisely, any local chart of the atlas
covering 
the manifold is reconstructed 
from the monodromy data. The glueing of the local charts is described by the action of 
the braid  group on the data, particularly, on the central connection 
matrix and on 
the Stokes' matrix \cite{Dub1}
 \cite{Dub2}. 

One well-known example of Frobenius 
manifolds is the quantum cohomology of smooth 
projective varieties \cite{Manin} \cite{KM} \cite{RT} \cite{MDS} .

It was conjectured \cite{Dub3} that the Stokes matrix for the quantum
cohomology of a  good Fano variety $X$  
is equal to the Gram matrix 
of the bilinear 
form $\chi(E,F):= \sum_k ~(-1)^k~$dim$~Ext^k(E,F)$  
computed on a  full collection of  exceptional objects in the 
derived category $Der^b(Coh(X))$ 
of coherent sheaves on $X$. 
More precisely, let $Der^b(Coh(X))$ be the derived 
category of coherent sheaves on a smooth projective variety $X$ of dimension 
$d$. 
An object $E$ of  
$Der^b(Coh(X))$ is called {\it exceptional} if $Ext^i (E,E)=0$ for 
$0<i<d$, $Ext^0(E,E)= 
{\bf C}$  and $Ext^d(E,E)$ is of the smallest dimension (if $X$ is a 
projective space, then $Ext^d(E,E)=0$).
  A collection  $\{ E_1,...,E_s \}$ of exceptional 
objects is an {\it exceptional collection} if for any $1\leq m<n\leq s$  we
have $Ext^{i}(E_n,E_m)=0$ for any  $i\geq 0$, $Ext^i(E_m,E_n)=0$ for any 
 $i\geq 0$  except possibly 
for one  value of $i$.   
A {\it full exceptional collection} is an exceptional collection 
which generates 
$Der^b(Coh(X))$ as a triangulated category. This theory is developed in 
 \cite{Rud1} \cite{Rud2} \cite{BP}.   We say that a Fano 
variety is {\it good} if it has a full exceptional collection. 

 It is known that $X={\bf CP}^d$ is good, the collection of sheaves on ${\bf CP}^d$ 
$\{ {\cal O}(n) \}_{n \in {\bf Z}}$ is exceptional, and $\{E_1,E_2,...,E_{d+1}\}
:=\{ {\cal O}, {\cal O}(1), ~...,
{\cal O}(d) \}$ is a full exceptional collection \cite{Bl}, \cite{GR}. In this case, 
$s_{ij}=\chi\left({\cal O}(i-1),{\cal O}(j-1)\right)$, 
$i,j=1,2,...,d+1$  has the ``canonical form'':
$$ 
  s_{ij}= \bin{d+j-i}{j-i}~~~~i<j$$
$$ 
  s_{ii}=1,~~~~~s_{ij}=0~~~~i>j$$
The inverse to this matrix has entries  $a_{ij}$ 
$$ 
   a_{ij}=(-1)^{j-i} \bin{d+1}{j-i}~~~~i<j
$$
$$ 
  a_{ii}=1,~~~~~a_{ij}=0~~~~i>j$$ 
This matrix is equivalent to the one above with respect to the action
of the braid group. We will also call it ``canonical''.

 The mentioned conjecture claims that the Stokes matrix of the quantum cohomology of 
${\bf CP}^d$ is equal to the above Gram matrix (modulo the action of the braid group: 
remarkably, this action on the Stokes matrix for the Frobenius manifold coincides 
with the natural action of the braid group on the collections 
 of exceptional objects 
\cite{Zas} \cite{Rud1}). 

\vskip 0.2 cm 

 This conjecture has its origin in the paper by Cecotti
and 
Vafa \cite{CV}, where another Stokes matrix introduced in \cite{Dub5} 
 for the $tt^*$ equations was found in the case of the ${\bf CP}^2$ 
topological $\sigma$ model. It was suggested, on physical arguments, 
that the entries of the Stokes' matrix
 $S=\pmatrix{ 1&x&y\cr
              0&1&z\cr
              0&0&1 \cr
             } $ 
 are integers. They must satisfy a Diophantine equation $x^2+y^2+z^2 -xyz=0$ 
 whose integer solutions 
$(x,y,z)$ are all equivalent to (3,3,3) modulo the action of the braid group. 
The authors of \cite{CV} also suggested that their matrix must coincide with the Stokes 
matrix defined in the theory of WDVV equations of associativity, that is, in the 
geometrical theory of Frobenius manifolds for 2D topological field theories \cite{Dub4}, \cite{Dub1}.

Later, in \cite{Zas}, the links between $N=2$ supersymmetric
field theories and the theory of derived categories were further investigated
and the 
coincidence of 
$\chi(E_i,E_j)$ with the Stokes matrix of $tt^*$ for ${\bf CP}^d$ 
was conjectured. 

 The conjecture may probably be derived from more general conjectures by 
Kontsevich in the framework of categorical mirror symmetry. To my knowledge,  
the subject was discussed in \cite{Konpis} (I thank  B. Dubrovin for this 
reference). 

\vskip 0.2 cm
 The main result of this paper is the proof (Theorem 2, $2^{\prime}$) 
that the conjecture about coincidence of the 
Stokes matrix for quantum cohomology of 
${\bf CP}^d$  and the Gram matrix $\chi(E_i,E_j)$ of a full
 exceptional collection 
in $Der^b(Coh({\bf CP}^d))$ is true. In this way, we generalize to any $d$ the result obtained in \cite{Dub2} for $d=2$. 
 
 We remark that it has not yet been proved that the Stokes' matrix 
for $tt^*$ equations and the Stokes' matrix for the corresponding 
Frobenius manifold coincide. This point deserves further investigation.

\vskip 0.2 cm 

We also study the structure of the monodromy group of the quantum cohomology of 
${\bf CP}^d$. The notion of {\it monodromy group} of a Frobenius manifold was introduced 
 in \cite{Dub1}. We prove (Theorem 3) that for $d=3$ the group is isomorphic to 
the subgroup of 
orientation preserving transformations in the hyperbolic triangular group $[2,4,\infty]$. 
In \cite{Dub2} it was proved that for $d=2$ the monodromy group is isomorphic to the 
direct product of the subgroup of orientation preserving transformations in $[2,3,\infty]$ and the cyclic group of order 2, $C_2=\{\pm \}$. Our numerical calculations also suggest that for any $d$ even the monodromy group may be 
 isomorphic to the orientation preserving transformations in 
$[2,d+1,\infty]$, and for any $d$ odd to the direct product of 
 the orientation 
preserving transformations in $[2,d+1,\infty]$ by $C_2$.

\vskip 0.2 cm 
\noindent
{\bf Acknowledgments}
\vskip 0.2 cm 
 I am indebted to Prof. B Dubrovin who introduced me to the problem
 and constantly gave me suggestions and advice. I also thank M. Bertola and
 M. Mazzocco for useful discussions. 
\vskip 0.3 cm

\section{ The system corresponding to ${\bf CP}^{k-1}$}
\vskip 0.2 cm

 We introduce here the linear system of differential equations 
whose Stokes matrices are the Stokes matrices for the quantum cohomology 
of ${\bf CP}^{k-1}$ (we use the more convenient choice $k=d+1$).  

In the quantum cohomology of ${\bf CP}^{k-1}$ we choose 
flat coordinate $t^1,..,t^k$ for 
the symmetric non degenerate {\it bilinear form} $<~,~>$ :
$$
   \eta_{\alpha \beta}=<\partial_{\alpha},\partial_{\beta}>= 
\delta_{ \alpha+\beta,k+1}
~~~\hbox{ where }\partial_{\alpha}={\partial \over \partial
t^{\alpha}}
$$    
Let $\eta$ be the matrix $(\eta_{\alpha \beta})$. In flat coordinates the 
 {\it Euler vector field} is 
$$ 
   E= \sum_{\alpha \neq 2} ~(1-q_{\alpha})t^{\alpha}{\partial \over 
\partial t^{\alpha} }
+~k {\partial \over \partial t^2}
$$
$$ q_1=0,~q_2=1,~q_3=2,~...,~q_{k}=k-1
$$
the {\it multiplication} is  
$$\partial_{\alpha}\cdot \partial_{\beta} = 
c^{\gamma}_{\alpha \beta}(t) \partial_{\gamma}  
~~~~~\hbox{ where }~~ c_{\alpha \beta \gamma}(t):=\eta_{\delta \gamma} 
c^{\delta}_{\alpha \beta}(t) = \partial_{\alpha}\partial_{\beta}
\partial_{\gamma} F(t)$$
$$
F(t)= {1 \over 6} \sum_{\alpha+\beta+\gamma=k+2}~t^{\alpha} 
t^{\beta} t^{\gamma} + \sum_{l=1}^{\infty} \Phi_{l}(t^3,...,t^k)
e^{lt^2}$$
$$
 \Phi_l(t^3,...,t^k)= \sum_{n=2}^{\infty} 
    \left[  \sum_{\alpha_1+...+\alpha_n= (l+1)(k-1)+l-1+n } I(l;\alpha_1,...,\alpha_n) 
{t^{\alpha_1}...t^{\alpha_n} \over n!}\right]
, ~~~~~~ I(1,k,k)=1
$$
and  the {\it unity vector field } is $e = {\partial \over \partial
t^1}$. Finally, let 
$${\mu}= \hbox{diag}(\mu_1,...,\mu_k)=\hbox{diag}(-{k-1\over
2},-{k-3 \over 2},...,{k-3 \over 2},{k-1 \over 2}),
~~~~~\mu_{\alpha}=q_{\alpha}-{d\over 2}$$

Consider the system of differential equations determining deformed
flat coordinates (see \cite{Dub1} \cite{Dub2}):
\be
 \partial_z \xi= ({\cal U}(t)+{1\over z}{\mu}) \xi
\label{zzzz}
\ee
$$ \partial_{\alpha} \xi = z C_{\alpha}(t) \xi $$
where $z\in{\bf C}$, $\partial_z:={\partial \over \partial z}$ 
and  $\xi$ is a column vector of components 
$$ 
\xi^{\alpha}= \eta^{\alpha\beta}\partial_{\beta}\tilde{t}(t,z)
$$ 
Here $(\eta^{\alpha\beta})=(\eta_{\mu \nu})^{-1}$ and $\tilde{t}(t,z)$ is 
one of the $k$ (deformed) flat coordinates. ${\cal U}(t)$ is the matrix of 
multiplication by the Euler vector field $E(t)$, and $C_{\alpha}(t)$ is
the matrix of entries  $(C_{\alpha})^{\beta}_{~\gamma}:=
c^{\beta}_{\alpha \gamma}$. The monodromy data of the system
(\ref{zzzz}) are, by definition, the {\it monodromy data} of the quantum
cohomology of ${\bf CP}^{k-1}$ in the local chart containing $t$.  
Let us compute $C_2(t)$ and ${\cal U}(t)$ at the 
semisimple point $(0,t^2,0,...,0)$: 
$$
  E \cdot  ~\partial_{\beta}=E^{\gamma} 
c^{\alpha}_{\gamma \beta} ~\partial_{\alpha} 
= E^2  c^{\alpha}_{2 \beta} ~\partial_{\alpha}=
$$
$$
   =k~c^{\alpha}_{2 \beta} \partial_{\alpha} \equiv  {\cal U}^{\alpha}_{~~\beta} 
\partial_{\alpha}
$$
Moreover $c_{2 \alpha \beta} (0,t^2,0,...,0)= \partial_2 
\partial_{\alpha} \partial_{\beta} F (0,t^2,...,0)$. This immediately 
yields 
$$
C_2(0,t^2,0,...,0) = \pmatrix{ 0 &   &       &     &  e^{t^2} \cr
                     1 & 0            \cr
                       & 1 & 0          \cr
                       &   & \ddots& \ddots \cr
                       &   &       &   1 & 0  \cr
              },
~~~~
{\cal U}(0,t^2,0,...,0)
 = \pmatrix{ 0 &   &       &     & k e^{t^2} \cr
                     k & 0            \cr
                       & k & 0          \cr
                       &   & \ddots& \ddots \cr
                       &   &       &   k  & 0  \cr
              }
$$ 
Let $y_{\alpha} 
:= \eta_{\alpha \beta} \xi^{\beta} \equiv \partial_{\alpha} \tilde{t}$. It 
satisfies 
\be
\partial_z y=\left(\hat{\cal U}(t^2)  - {1 \over z} \mu \right) y
\label{10p}
\ee
\be
   \partial_2 y= z \hat{C}_2(t^2) y
\label{zzz1}
\ee
where
$$ \hat{C}_2(t^2):=\eta~C_2(0,t^2,...,0)~\eta^{-1}=
 \pmatrix{ 0 &1        \cr
                      &0&1             \cr
                       &  & 0&1          \cr
                       &   & &\ddots& \ddots \cr
                       &   & &      &   0 &1 \cr                      
                      e^{t^2} &   & &      &    &  0  \cr
              }
$$
$$
\hat{\cal U}(t^2):=\eta~{\cal U}(0,t^2,...,0)~\eta^{-1}
=\pmatrix{ 0 &k       \cr
                      &0&k             \cr
                       &  & 0&k          \cr
                       &   & &\ddots& \ddots \cr
                       &   & &      &   0 &k \cr                      
                      ke^{t^2} &   & &      &    &  0  \cr
              }
$$

\vskip 0.2 cm
\noindent
{\bf Lemma 1}: {\it Let ${y}(t^2,z)=(y_1(t^2,z),...,y_k(t^2,z))^T$ be  a column
vector 
solution of the above 
system  (\ref{10p}) (\ref{zzz1}).  With the following substitution
\be
 y_{\alpha}(t^2,z)= {1 \over k^{\alpha-1}} ~z^{{k-1\over 2}-\alpha+1}~ 
                    (z \partial_z)^{(\alpha-1)}~\varphi(t^2,z)   ~~~
\label{40p}
\ee
\be
  \equiv z^{{k-1\over 2}-\alpha+1}~
  \partial_2^{\alpha-1}\varphi(t^2,z),~~~ \alpha=1,2,...,k
\label{zzz2}
\ee
the above system is equivalent to the equations 
\be
(z \partial_z)^{k} \varphi = (k z)^k e^{t_2} ~\varphi
\label{30p}
\ee
\be
\partial_2^k\varphi =z^k e^{t_2} ~\varphi
\label{zzz3}
\ee
}
\vskip 0.2 cm
\noindent
The proof is a simple calculation we leave to the reader.

 \vskip 0.2 cm
The substitution of the lemma implies $\partial_2 \varphi= {1\over k} 
z \partial_z \varphi$. Then 
$$ 
  e^{t^2\over k} {\partial \over \partial e^{t^2 \over k}}\varphi(t^2,z) =
z  {\partial \over \partial z } \varphi(t^2,z)
$$
which implies (with abuse of notation)
$$\varphi(t^2,z) \equiv \varphi(ze^{t^2\over k})$$
Namely, $\varphi$ (at $(0,t^2,...,0)$) depends on one argument 
$w=ze^{t^2\over k}$ and satisfies the {\it generalized hypergeometric
equation} 
\be 
   \left(w {d \over dw}\right)^k \varphi(w) = (kw)^k~\varphi(w)
\label{30pir}
\ee
The equation is equivalent to the system 
 \be
     {d Y \over dw}= \left[ \hat{\cal U}+{\hat{\mu}\over w} \right] Y
\label{10pir}
\ee
where 
\be
 Y_n(w)= {1 \over k^{n-1}} ~w^{{k-1\over 2}-n+1}~ 
                    (w \partial_w)^{(n-1)}~\varphi(w)   ~~~n=1,2,...,k
\label{40pir}
\ee 
$$
   \hat{\cal U}:=\hat{\cal U}(0)= \left[ \matrix{ 0 & k &   &           &          &\cr
                               & 0 & k &           &          & \cr 
                               &   & 0 & k         &          & \cr 
                               &   &   & \ddots    &   \ddots & \cr
                               &   &   &           &   \ddots & k \cr
                             k &   &   &           &          & 0 \cr 
}   
                          \right]
$$
$$
\hat{\mu}:=-\mu=\hbox{ diag}\left({k-1\over 2}, {k-3\over 2}, {k-5 \over 2}, ..., -{k-3
\over 2}, -{k-1 \over 2} \right)
$$
 The system (\ref{10pir}) may also be interpreted as the system (\ref{10p}) with
 $t^2=0$. We will return
 later (section 4)  on the connection between its monodromy data and
 the monodromy data of the system (\ref{10p}).

\vskip 0.3 cm 
Let us study system (\ref{10pir}). We change notation and  
choose the more familiar letter $z$
instead of $w$.  So, the system (\ref{10pir}) is re-written as 
\be
     {d Y \over dz}= \left[ \hat{\cal U}+{\hat{\mu}\over z} \right] Y
\label{10}
\ee
and (\ref{40pir}) (\ref{30pir}) become
 \be
 Y_n(z)= {1 \over k^{n-1}} ~z^{{k-1\over 2}-n+1}~ 
                    (z \partial_z)^{(n-1)}~\varphi(z)   ~~~n=1,2,...,k
\label{40}
\ee
\be 
   \left(z {d \over dz}\right)^k \varphi(z) = (kz)^k~\varphi(z)
\label{30}
\ee
The point $z=0$ is e fuchsian singularity, and $z=\infty$ is a
singularity of the second kind.  
(\ref{10}) has a fundamental matrix solution $Y_0(z)$
whose behaviour at $z=0$ is 
$$
   Y_0(z)= (I+O(z)) z^{\hat{\mu}}z^{R}~~~~~~
    R=\left[ \matrix{ 0 & k &   &           &          &\cr
                               & 0 & k &           &          & \cr 
                               &   & 0 & k         &          & \cr 
                               &   &   & \ddots    &   \ddots & \cr
                               &   &   &           &   \ddots & k \cr
                              &   &   &           &          & 0 \cr 
}   
                          \right]
$$
and the monodromy for a counter-clockwise loop around the origin is
$e^{2 \pi i (\hat{\mu} +R)}$.

  The characteristic polynomial of the matrix $\hat{\cal U}$ is $0=\det(
\hat{\cal U} - u)= (-u)^k+(-1)^{k+1} k^k$. It  has $k$ distinct
eigenvalues $u_{n}=k~e^{2 \pi i (n-1) \over k}$, $n=1,...,k$. The
equations for the eigenvector ${\bf x}_n$  corresonding to $u_n$, namely
$\hat{\cal U}{\bf x}_n=u_n~{\bf x}_n$, written for the components ${x^1}_n$,
...,${x^k}_n$ of the column vector ${\bf x}_n$ are
$$
           {x^{l+1}}_n=e^{2 \pi i (n-1) \over k}~{x^l}_n ~~~~ l=1,2,...,k-1
$$
$$
          {x^1}_n= e^{2 \pi i (n-1) \over k} {x^k}_n
$$
With the choice ${x^1}_n= e^{i \pi (n-1) \over k}$ we get ${\bf x}_n=
(e^{i \pi (n-1) \over k},e^{i 3 \pi (n-1) \over k},e^{i 5 \pi (n-1) \over k},...,
e^{-{i \pi (n-1) \over k}})^T$ ($T$ stends for transpose). The  matrix 
$$
   X= {1 \over \sqrt{k}}~[ {\bf x}_1 |   {\bf x}_2 |...|{\bf
   x}_k]= {1 \over \sqrt{k}}~({x^j}_n)~~~~~{x^j}_n=e^{(2j-1)i \pi {n-1 \over
   k}}~~j,n=1,2,...,k
$$
puts ${\cal U}$ in diagonal form:
$$
  U= X^{-1} \hat{\cal U} X =
\hbox{diag}(u_1,u_2,...,u_n,...,u_k)~~~~u_n=k~e^{2 \pi i (n-1) 
                                                        \over k}
$$
We stress that $u_i \neq u_j$ for $i \neq j$. 
The system (\ref{10}) is transformed by the gauge $X$ in an equivalent
form 
\be
     {d \tilde{Y} \over dz}= \left[ U+{V\over z} \right] \tilde{Y}
\label{20}
\ee
$$ 
  \tilde{Y}= X^{-1}Y,~~~U=X^{-1}\hat{\cal U}X,~~~V= X^{-1} \hat{\mu} X
$$
Observe that 
$$ 
\eta \hat{\mu} + \hat{\mu} \eta=0 
~~~~~~~      X X^T = \eta^{-1}
$$
This implies that $V$ is skew-symmetric
$$
V+V^T=0
$$
With the gauge $X$, $Y_0(z)$ transforms into 
$$ 
\tilde{Y}_0(z)= (X^{-1} +O(z)) ~z^{\hat{\mu}}z^R, ~~~~z\to 0
$$


\vskip 0.3 cm
\section{Asymptotic Behaviour and Stokes' Phenomenon}
\vskip 0.2 cm

 Our aim is to explicitely compute a Stokes' matrix for the above
system (\ref{20}), or for the system (\ref{10}).  
The system (\ref{20}) has  formal solution
$$ 
   \tilde{Y}_F = \left[ I + {F_1 \over z}+{F_2 \over z^2}+ ... \right]~e^{z~U}
$$
where $F_j$'s are $k \times k$ matrices. It is a well known result
that fundamental matrix  solutions exist which have $\tilde{Y}_F$ as 
asymptotic expansion for $z \to \infty$ in some ``admissible'' sectors of the
complex 
plane of
angular width greater then $\pi$. In order to find such sectors we
need the so called {\it Stokes' rays}, defined by 
$$ 
   \Re  ~((u_r-u_s)z)=0 ,~~r \neq s~~~  \Rightarrow ~~~\arg z= -\arg
   (u_r-u_s)+{\pi \over 2} + m \pi ,~~~m=0,1 
$$
 There exists a {\it unique}  solution of the system 
asymptotic to $\tilde{Y}_F$ in a sector greater then $\pi$
and bounded by the first two Stokes' rays we meet extending over $\pi$
the
angular width of the sector. The general theory of Stokes' phenomenon
is found in the classical paper by W. Balser, W.B. Jurkat, D.A. Lutz
\cite{BJL1}.   Stokes matrices are also a natural monodromy datum in
the theory of isomonodromy deformation 
 \cite{ITS} \cite{JMU} \cite{JM1} \cite{JM2} \cite{Dub1} \cite{DM}. 

\vskip 0.2 cm
A possible choice for the labelling of the rays is the following: we
call $R_{rs}$ the Stokes' ray 
$$ 
   R_{rs}= \{ z=- i \rho (\bar{u_r} - \bar{u_s}), ~~~\rho>0 \}~~~~~r
   \neq s
$$
 \vskip 0.3 cm
\noindent
{\bf Lemma 2}: {\it  For $r<s$ the Stokes' rays of the system
(\ref{20}) are 
$$ 
         R_{rs}= \left\{z= \rho ~\exp \left( i\left[ {2 \pi \over k} - {\pi
                    \over k} (r+s)\right] \right),~~~\rho>0 \right\}
$$
$$
        R_{sr}= -R_{rs}
$$
} 
\vskip 0.2 cm
\noindent
{\bf Proof}: Just compute 
$$ -i (\bar{u_r} -\bar{u_s})= -i ( e^{-i {2 \pi \over k}(r-1) }-e^{-i{2
\pi \over k}(s-1)})=
$$
$$
   = 2 \sin \left({\pi \over k}(s-r)\right) ~e^{\left(i\left[ {2 \pi \over k} - {\pi
                    \over k} (r+s)\right] \right)}
$$ 
 Then we note that $ \sin \left({\pi \over k}(s-r)\right)$ is positive
 because $0<s-r\leq k-1$. 

\rightline{$\Box$}

\vskip 0.3 cm
\noindent
{\bf  Remark 1}: $R_{rs}=R_{pq}$ for $r+s=p+q$. $R_{12}$ is at $\arg z =
- {\pi \over k}$, $R_{13}$ is at $\arg z = - {2 \pi \over k}$, and so
on. For $r+s = k+2$ the corresponding $R_{rs}$'s are at $\arg z =
-\pi$ and the $R_{sr}$'s are at $\arg z =0$.   $R_{k-1,k}$ is at the
angle $-2 \pi +{3 \pi \over k}$ or, equivalently, at ${3 \pi \over
k}$. See figure 1. 

\begin{figure}
\epsfxsize=15cm
\centerline{\epsffile{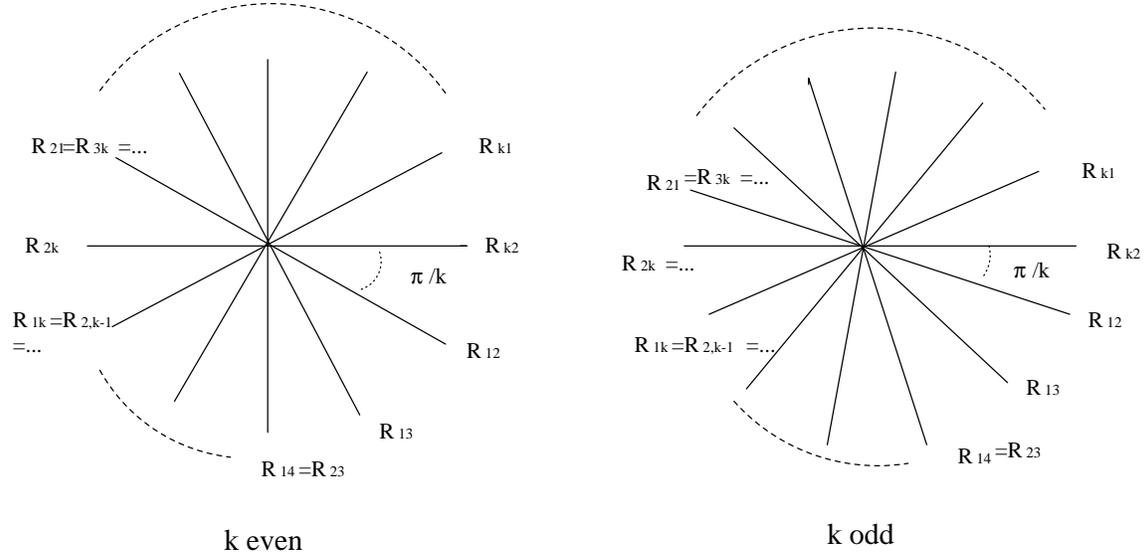}}
\caption{ Stokes' rays}
\end{figure}

\vskip 0.2 cm

  We choose two ``admissible'' overlapping sectors in a canonical
  way. Let $l$ be an ``admissible'' line through the origin, namely a
  line not containing Stokes' rays. For our purposes we take 
$$
  l= \{ z ~|~ z= \rho e^{i \epsilon}, ~~~\rho \in{\bf R},~~~ 0<\epsilon<
  {\pi \over k} \}  
$$
$l$ has a natural orientation inherited from {\bf R}. We call $\Pi_R$
and $\Pi_L$ the half planes to the right/left of $l$ w.r.t its
orientation.
$$ 
  \Pi_R=\{ - \pi +\epsilon < \arg z <
\epsilon\}~~~\Pi_L=\{\epsilon<\arg z < \pi + \epsilon \}$$    
We then define two different ``admissible'' sectors ${\cal S}_L$,
${\cal S}_R$ which
contain $l$ 
$$ 
  {\cal S}_L=\{ z\in {\bf C}~|~ 0<arg z < \pi + {\pi \over k}
  \}\supset
              \Pi_L
$$
$$
    {\cal S}_R=\{ z\in {\bf C}~|~ - \pi <arg z < {\pi \over k} \}
    \supset
              \Pi_R
$$
  We call che
corresponding  solutions $\tilde{Y}_L(z)$ and $\tilde{Y}_R(z)$. 

\vskip 0.2 cm
\noindent
{\bf Definition: } The {\it Stokes' matrix} of the system (\ref{20}) 
with respect to the admissible line $l$ is the connection matrix $S$  
such that 
$$ 
    \tilde{Y}_L(z)=\tilde{Y}_R(z)S~~~~~~0<\arg z<{\pi \over k}
$$
On the opposite overlapping region one can prove (as a consequence of
the skew-symmetry of $V$, see \cite{Dub1}) that 
$$
 \tilde{Y}_L(z)=\tilde{Y}_R(z~e^{-2\pi i})S^T~~~~~~\pi<\arg z<\pi+{\pi \over k}
$$
In \cite{BJL1} $S$ is called a `` Stokes' multiplier''. The
terminology in this field changes from one author to the other...   

\vskip 0.2 cm
\noindent 
{\bf Definition: } We call  {\it central connection
 matrix} the connection matrix $C$ such  that
$$ 
   \tilde{Y}_0(z)=\tilde{Y}_{R}(z) C   ~~~z\in \Pi_R $$


\vskip 0.2 cm

  It is clear that the system (\ref{10}) has solutions $Y_0(z)=X
\tilde{Y}_0(z)$, and  $Y_L(z)= X
\tilde{Y}_L(z)$,  $Y_R(z)= X \tilde{Y}_R(z)$ asymptotic to 
$X \tilde{Y}_F(z)$ as $z \to \infty$ in ${\cal S}_L$ and ${\cal S}_R$
respectively, which are connected by the same $S$ and $C$. 
 
\vskip 0.2 cm

 In order to compute the entries of $S$ explicitely, we use the
reduction of  
(\ref{10}) to the generalized hypergeometric equation (\ref{30}).  
If
$\varphi^{(1)}(z)$, ..., $\varphi^{(k)}(z)$ is a basis of $k$ linearly
independent solutions of (\ref{30}), then the matrix $Y(z)$ of entries
$(n,j)$ defined by   
\be
Y^{(j)}_n(z):= {1 \over k^{n-1}} ~z^{{k-1\over 2}-n+1}~ 
                    (z \partial_z)^{(n-1)}~\varphi^{(j)}(z)
\label{entries} 
\ee
is a fundamental matrix for
(\ref{10}).

\vskip 0.3 cm
\noindent
{\bf Lemma 3:} {\it 
                    The generalized hypergeometric equation (\ref{30})
has two bases of linerly independent solutions 
$\varphi_{L/R}^{(1)}(z)$, ..., $\varphi_{L/R}^{(k)}(z)$ having
                    asymptotic behaviours 
$$ 
   \varphi_{L/R}^{(n)}={1 \over \sqrt{k}} {e^{i{\pi \over k}(n-1)}
   \over 
          z^{k-1 \over 2}}
                         ~ \exp\left[k~e^{i{2 \pi \over k}(n-1)} ~z
   \right]~\left(1 + O\left({1 \over z}\right) \right),~~~~~z \to \infty
$$
in ${\cal S}_L$ and ${\cal S}_R$ respectively. Let $\Phi(z)$ denote
 the row vector 
$ [ \varphi^{(1)}(z),...,\varphi^{(k)}(z)]$.  The fundamental
 matrices  $Y_L(z)$, $Y_R(z)$ of (\ref{10}) 
 are expressed through formula (\ref{entries}) in terms of $\Phi_L(z)$
 and $\Phi_R(z)$ and 
$$
Y_L(z)=Y_R(z)~S~~~~~~~0<\arg z < {\pi \over k}
$$
if and only if 
$$ 
  \Phi_L(z)= \Phi_R(z) S ~~~~~~~0<\arg z < {\pi \over k}
$$
}
\vskip 0.2 cm 
\noindent
{\bf Proof:} Simply observe that for a fundamental solution  in
${\cal S}_L$ or ${\cal S}_R$ (we omit subscripts $L$, $R$)
$$ 
  Y(z)= \left[ \matrix{ 
                         z^{k-1\over 2} \varphi^{(1)}&  ... & 
         z^{k-1\over 2} \varphi^{(k)}\cr
       \vdots & ... & \vdots \cr 
} 
\right] = X \tilde{Y}
$$
which is asymptotic, for $z \to \infty$, to 
$$
\sim \left[  \matrix{  1 & e^{i {\pi \over k}}&e^{i {2\pi \over k}}
  & e^{i {3\pi \over k}}& ...& e^{i {k-1 \over k}\pi } \cr 
   \vdots &\vdots &\vdots & \vdots & ... & \vdots \cr
                       }
             \right]  ~
                        \left[
\matrix{
\exp(kz)&                              &  & \cr
        &\exp(k e^{2 \pi i \over k}z) &     &  \cr
        &                             & \ddots &  \cr
        &                             &        &  
                           \exp( k e^{2 \pi i (k-1)\over k}) \cr
}
\right]
$$
Now , the first row   of $Y(z)$ is $z^{k-1 \over 2} \Phi(z)$

\rightline{$\Box$}

\vskip 0.2 cm
 The Stokes' matrix $S$ has entries  
$$ 
                 s_{ii}=1 
$$
$$ 
               s_{ij}=0~~~\hbox{ if } R_{ij}\subset \Pi_R
$$
This follows from the fact that on the overlapping reagion $0<\arg z
<{\pi \over k}$ there are no Stokes' rays and 
$$ 
   e^{zU}~S~e^{-zU} \sim I, ~~~~z \to \infty,~~~~
\hbox{ then }~~ e^{z(u_i-u_j)}~s_{ij}\to \delta_{ij}
$$
 Moreover, $\Re \left(z(u_i-u_j)\right)>0$ to the
left of the ray $R_{ij}$, while  $\Re \left(z(u_i-u_j)\right)<
0$ to the right (the natural orientation on $R_{ij}$, from
$z=0$ to $\infty$ is understood). This implies
$$ 
   \left|e^{zu_i} \right| > \left|e^{zu_j} \right|~~~~~ \hbox{ and  }
                      e^{z(u_i-u_j)}\to \infty ~~\hbox{ as } z \to \infty 
$$
on the left, while on the right 
$$ 
   \left|e^{zu_i} \right| < \left|e^{zu_j} \right|~~~~~ \hbox{ and  }
                      e^{z(u_i-u_j)}\to 0 ~~\hbox{ as } z \to \infty 
$$
With this observations in mind, we prove the following

\vskip 0.3 cm 
\noindent 
{\bf Lemma 4 :} {\it  $S$ has a column whose entries are all zero but
one.  More precisely:
\par \noindent
 For $k$ even
$$  s_{i,{k \over 2}+1}=0~~~~\forall i \neq {k \over 2}+1 ,
~~~~~~~ s_{{k \over 2}+1,{k \over 2}+1}=1 $$
For $k$ odd 
$$
s_{i,{k+1 \over 2}}=0~~~~\forall i \neq {k+1 \over 2} ,
~~~~~~~ s_{{k +1\over 2},{k+1 \over 2}}=1 $$
}

\vskip 0.2 cm 
\noindent
{\bf Proof: } Let us determine $n$ such that $s_{in}=0$ for any $i \neq
n$ and $s_{nn}=1$. We need to find all rays in $\Pi_R$. We start with
$R_{rs}$ with $r<s$. We know that for $r+s=k+2$ the ray is the
negative real half-line (at angle $-\pi$). Then $R_{rs} \subset \Pi_R$
for $r+s\leq k+1$ ($r<s$). Then, in $\Pi_R$ we have 
$$ 
  \left. \matrix{ R_{12} & R_{13} & ... & R_{1k} \cr 
                 R_{23} & R_{24} & ... & R_{2, k-1} \cr 
                 R_{34} & R_{35} & ... & R_{3,k-2} \cr
                 \vdots & \vdots &     &          \cr
                 R_{ab} &        &     &          \cr
               }
                \right. 
$$
where 
$
         R_{ab} =
                          R_{{k \over 2},{k \over 2} +1}$ for $k$ 
                          even, and 
                         $ R_{{k-1 \over 2},{k+3\over 2}}$  for $k$ 
                          odd. 
In $\Pi_R$ we have also $R_{rs}$ with $r+s\geq k+2$ and $r>s$. For 
fixed $n$ we require $R_{in} \subset \Pi_R$ for any
$i$. Namely, 
$$
     \forall i<n ~~~i+n \leq k+1 ~,~~~~~~\forall i>n ~~~i+n\geq k+2
$$ 
This yields $n= {k \over 2}+1 $ for $k$ even, $n={k+1 \over 2}$ for
$k$ odd. 

\rightline{$\Box$}                      
              
\vskip 0.3 cm 
Let $n(k)$ be $ {k \over 2}+1 $, or  ${k+1 \over 2}$. 
Lemma 4 implies that the $n(k)^{th}$ columns of ${Y}_L$ and
${Y}_R$ coincide. In particular, their asymptotic representation holds
for $-\pi < \arg z < \pi +{\pi \over
k}$. Actually, this domain can be further enlarged, up to 
$$  
    -{\pi \over k} -\pi  < \arg z < \pi +{\pi \over k} ~~~~k \hbox{
    even} 
$$
$$ 
   -\pi < \arg z < \pi +  {2\pi \over k} ~~~~k \hbox{ odd } 
$$ 
To see this recall that $\left|e^{zu_i} \right|< \left|e^{zu_j}
\right|$ on the right of $R_{ij}$, and conversely on the left. Then it
is easy to see that for $k$ even $\left|\exp ({z~u_{{k \over 2}+1}})
\right|$ dominates  all exponentials in the sector $ -{\pi \over
k} - \pi < \arg z < {\pi \over k} - \pi$, while for $k$ odd 
$\left|\exp ({z~u_{{k +1\over 2}}})\right|$ dominates  all exponentials in
the sector $ \pi < \arg z< \pi + {2 \pi \over k}$. 
 
The first entry of the $n(k)^{th}$ column is
$\varphi^{(n(k))}_L(z)\equiv \varphi^{(n(k))}_R(z)$ times $z^{k-1
\over 2}$. Then $\varphi^{(n(k))}$ has the estabilished asymtotic
behaviour on the enlarged domains above.

\vskip 0.3 cm
 We now introduce an integral representation for a solution
 $\varphi(z)$ of the generalized hypergeometric equation which will
 allow us to compute the entries of $S$.

\vskip 0.2 cm 
\noindent
{\bf Lemma 5:} 
              { \it  The function 
$$ 
     g^{(n)}(z)= {1 \over (2 \pi )^{k+1 \over 2} ~e^{i \pi ( {k \over
     2}-n-1)}} ~ \int_{-c -i \infty}^{-c+i
     \infty}ds~\Gamma^k(-s)~e^{-i\pi k s} ~e^{i2(n-1)\pi s}~z^{ks} 
$$
defined for ${\pi \over 2} -2(n-1){\pi \over k}< \arg z < {3 \pi \over
2} - 2 (n-1) {\pi \over k} $, $z \neq 0$ and for any positive number
$c>0$, is a solution of the generalized hypergeometric equation
(\ref{30}) ( the
path of integration is a vertical line through $-c$). It has asymptotic behaviour 
$$ 
  g^{(n)}(z)\sim {1 \over \sqrt{k}}~{e^{i{\pi \over k}(n-1)} \over
  z^{k-1\over 2}}~ \exp \left(k e^{i{2 \pi \over k}(n-1)}z
  \right)~~~~~~z\to \infty
$$
In particular, for $n(k)= {k\over 2}+1$ ($k$ even), or $n(k)={k+1 \over 2}$
($k$ odd), the analytic continuation of $g^{(n(k))}(z)$ has the above
asymptotic  behaviour in the domains 
$$ 
      - \pi - {\pi \over k} < \arg z< \pi + {\pi \over k}~~~~~k
      \hbox{ even }
$$
$$
  - \pi < \arg z < \pi + { 2 \pi \over k}   ~~~~~k \hbox{ odd } 
$$ 
and coincides with the solution $\varphi^{(n(k))}_L\equiv
\varphi^{(n(k))}_R$ appearing in the first rows of the fundamental
matrices $Y_L$ and 
$Y_R$ of the system (\ref{10}).    

 The following identity holds 
\be
   \sum_{m=0}^k ~(-1)^{m-k}~\pmatrix{ k \cr m \cr} ~g^{(n)}(z~e^{i {2 \pi
\over k}m}) =0
\label{60}
\ee     
} 
 
We'll sketch the proof in Appendix 1. 

\vskip 0.2 cm 
\noindent 
{\bf Remark 2:} 
Observe that for   basic solutions of the  hypergeometric equation
  $\Phi_{L/R}= [ \varphi^{(1)}_{L/R},..., \varphi^{(k)}_{L/R}]$, 
$$
 \varphi^{(n)}(z )\sim 
~{1 \over \sqrt{k}}  {e^{i {\pi \over k} (n-1)}\over  z^{k-1\over
  2} }\exp( k  e^{i {2\pi \over k}
  (n-1)} z)
$$
on some sector, and 
$$ 
 \varphi^{(n)}(z e ^{2 \pi i \over k})\sim 
(-1)~{1 \over \sqrt{k}}  {e^{i {\pi \over k} ([n+1]-1)}\over  z^{k-1\over
  2} }\exp( k  e^{i {2\pi \over k}
  ([n+1]-1)} z)
$$
like $- \varphi^{(n+1)}$, on the sector rotated by ${- 2 \pi \over k}$. Note however
that $\varphi^{(k)} (ze^{2 \pi i \over k}) \sim  (\sqrt{k} z^{k-1\over
  2} )^{-1}~e^{kz} $, like $\varphi^{(1)}(z)$. 

 Also, note that $g^{(n)}(z e^{2\pi i \over k})= - g^{(n+1)}(z)$ (we
 mean analytic continuations), with
 asymptotic behaviour on rotated domain.

\section{ Monodromy Data of the Quantum Cohomology of  ${\bf CP}^{k-1}$}

\vskip 0.2 cm
Let us return to the system (\ref{10p}).
\be
\partial_z y=\left(\hat{\cal U}(t^2)  + {1 \over z} \hat{\mu} \right) y
\label{tuono}
\ee
 In this
section we use the original notation $w=ze^{t^2\over k}$.
The system has a fundamental matrix 
$$
 y_0(t^2,z)=(I+A_1(t^2)z+A_2(t^2)z^2+...)~z^{\hat{\mu}}~z^R,~~~~z\to 0
$$
where $R$ is the same of system (\ref{10}). The series appearing in
the solutions converges near $z=0$. 
The matrix $\hat{\cal U}(t^2)$ has eigenvalues and eigenvectors
$$
   u_n(t^2)=e^{i{2\pi\over k}(n-1)}~e^{t^2 \over k} ~\equiv~u_n e^{t^2
   \over k}~~~~~n=1,...,k$$
$$
   {\bf x}_n(t^2) \hbox{ of entries } {x^j}_n(t^2)=              
e^{i(2j-1)(n-1){\pi   \over k}}~e^{{2j-1-k \over 2k}t^2}
~\equiv {x^j}_n e^{{2j-1-k \over 2k}t^2}$$
Let $X(t^2)=({x^j}_n(t^2))$. With the gauge $y=X(t^2)~\tilde{y}(t^2,z)$
we obtain the
equivalent system 
\be
    \partial_z \tilde{y} = \left[ U(t^2)+{V(t^2)\over z}\right]
    \tilde{y}
\label{tuono1}
\ee
$$
   U(t^2)=X^{-1}(t^2)~\hat{\cal
   U}(t^2)~ X(t^2)=\hbox{diag}(u_1(t^2),...,u_k(t^2))
$$
  $$ V(t^2) =X^{-1}(t^2)~\hat{\mu} ~X(t^2)) 
$$
$$
V(t^2)^T+V(t^2)=0
$$
The skew symmetry of $V$ follows from $\eta\hat{\mu} +\hat{\mu}\eta=0$
and from the choice of the normalization
of the eigenvectors, such that $X(t^2)^T~X(t^2)=\eta^{-1}$.

Let us fix an initial point $t_0=(0,t_0^2,0,...,0)$. 
The system (\ref{tuono1}) has fundamental matrices $y_R(t_0^2,z)$,
$y_L(t_0^2,z)$, which are asymptotic to the formal solution 
$$ 
  \tilde{y}_F(t_0^2,z)=\left[ I+{F_1(t_0^2)\over z}+{F_2(t_0^2)\over z}+
... \right]~e^{z~U(t_0^2)}
$$
in the sectors 
$$ 
  {\cal S}_L(t_0)=\{ z\in {\bf C}~|~ 0<\arg \left[z \exp\left({t_0^2\over
k}\right)\right]   
< \pi + {\pi \over k}
  \}
$$
$$
    {\cal S}_R(t_0)=\{ z\in {\bf C}~|~ - \pi <
\arg \left[z \exp\left( {t_0^2\over k}\right)\right]    < {\pi \over k} \}
    $$
and 
$$ \tilde{y}_L(t_0^2,z)=y_R(t_0^2,z)~ S ~~~~~~~~~~~0<\arg \left[z
\exp\left({t_0^2\over k}\right)\right]<{\pi \over k}
$$
with respect to the admissible line
$$ 
   l_{t_0}:= \{z~|~z=\rho \exp \left({i \epsilon} -{\Im m~ t_0^2\over k}
   \right),~~\rho>0\}
$$
  The Stokes' matrix is precisely the matrix $S$ of system
  (\ref{10}) with respect to the admissible line $l_{t_0}$. 
Also the central connection matrix defined by
$$
   y_0(t_0^2,z)=y_R(t_0^2,z)~C ~~~~~~~~~~ - \pi <\arg \left[
z e^{t_0^2 \over k} \right] < {\pi \over k}
$$
is the the same of the system (\ref{10}).

\vskip 0.2 cm
\noindent
{\bf Definition: } $C$ and $S$, together with 
$\hat{\mu}$, $R$, and $e={\partial \over
\partial t^1}$  are the {\it monodromy 
data} of the quantum cohomology of ${\bf CP}^{k-1}$ in a local chart 
containing $t_0$. 
\vskip 0.2 cm

Recall that we fixed a point $t_0=(0,t_0^2,...,0)$. When we consider a
point $t$  away from $t_0$, the system (\ref{tuono}) acquires the
general form 
\be
 \partial_z y = \left[ \hat{\cal U}(t) +{\hat{\mu}\over z} \right]y
\label{zzzz1}
\ee
where  $\hat{\cal U}(t^1,..,t^k)= \eta~ {\cal U}(t)~
\eta^{-1}$ and  $y_{\alpha}^{(j)}(t^1,..,t^k;z)=\partial_{\alpha}
\tilde{t}^{j}(t,z)$. 

The admissible line $l_{t_0}$ must be considered {\it
fixed once and for all}. Instead, the Stokes' rays change. This is
because they are functions of the eigenvalues $u_1(t)$, ..., $u_k(t)$
of the matrix $\hat{\cal U}(t^1,..,t^k)$. For example, if just $t^2$ 
varies, while $t^1=t^3=...=t^k=0$, the system (\ref{tuono}) has Stokes' rays 
$$
  R_{rs}(t^2)=\{z~|~=\rho~ \exp\left(i{2\pi \over k} -i{\pi \over
  k}(r+s) -i {\Im m~t^2 \over k} \right),~~\rho>0 \}
$$
The dependence of the coefficients of the system (\ref{zzzz1}) on $t$
is isomonodromic \cite{Dub1} \cite{Dub2}. Then $\mu$ and $R$ are the same
for any $t$. $S$ and $C$ do not change if we move in a sufficently small 
neighbourhood of
$t_0$.  Problems arise when some Stokes' rays cross $l_{t_0}$. $S$ and $C$
must be modified by an action of the braid group. 
 We will return to this point later.


\vskip 0.2 cm 
\section{ Computation of S  }

\vskip 0.2 cm
 To compute $S$, we factorize it  in ``Stokes' factors''. 
Our fundamental matrix
 $Y_L$ has the required asymptotic form on the sector between
 $R_{k2} ~~(\arg z=0)$ and $ R_{1k}~~(\arg z= \pi +{\pi \over
 k})$. $Y_R$ has the same behaviour between $R_{2k}~~(\arg z= -\pi)$
 and $R_{k1} ~~( \arg z = {\pi \over k})$.

\begin{figure}
\epsfxsize=9cm
\centerline{
\epsffile{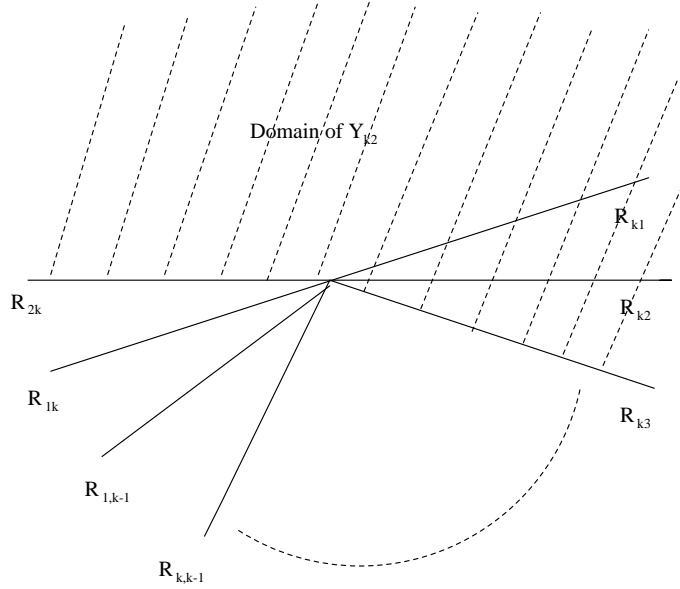}}
\caption{ Sector where $Y_{k2}$ has the asymptotic behaviour
$X~e^{zU}$ for $z \to \infty$}
\end{figure}

 Of course, we can consider fundamental matrices with the same
asymptotic 
behaviour on other 
 sectors of angular width less then $\pi +{\pi \over k}$ and bounded
 by two Stokes' rays. We introduce the following notation: consider a
 fundamental matrix of (\ref{10}) having the required asymptotic
behaviour on such a sector.  If we go all over the sector  clockwise 
 we meet Stokes rays belonging to the sector at each
 displacement of ${\pi \over k}$. Let $R_{ij}$  be  the last ray we
 meet before reaching the boundary (the boundary is still a Stokes ray {\it
 not } belonging to the sector). Then we will call the fundamental
 matrix $Y_{ij}$. For example,
 $Y_L=Y_{k1}$ and $Y_R=Y_{1k}$. See figure 2.

Sometimes, a different labelling is
 used in the literature. The rays must be enumerated as in figure
 3. The numeration refers to the line $l$: the rays  are labelled in
counter-clockwise order starting from the first one in $\Pi_R$ (which
will be $R_0$; then $R_0$, $R_1$, ..., $R_{k-1}$ are in $\Pi_R$, and
$R_k$, ..., $R_{2k-1}$ are in $\Pi_L$).  
    For our particular  choiche of $l$ , $R_0\equiv R_{1k}$ 
(at $\arg z = -\pi + {\pi \over k}$); $R_1\equiv R_{1,k-1}$ follows 
 counter-clockwise... Then we proceed untill we reach $R_{k-1}\equiv R_{k2}$ before
 crossing $l$, and so on. The fundamental matrices are labelled as
 we prescribed above, namely $Y_j$ if its sector contains $R_j$ as the
 last ray met going all over the sector clockwise before the
 boundary. The sector itself is denoted by ${\cal S}_j$. See figure 3.

\begin{figure}
\epsfxsize=9cm
\centerline{
\epsffile{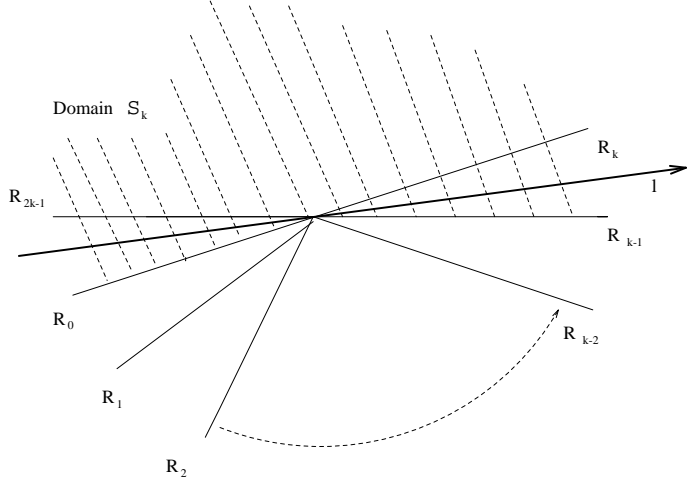}}
\caption{ Sector ${\cal S}_k$ and labelling of Stokes' rays}
\end{figure}

We define {\it Stokes' factors} to be the connection matrices $K_j$
 such that  
$$ 
    Y_{j+1} (z) = Y_j(z) K_j
$$
on the overlapping region of width $\pi$. We warn the reader that
also the Stokes' factors will be labelled with both conventions above,
according to our convenience (for example $K_0 \equiv K_{1k}$). 

 As a consequence of the above definitions, we can factorize $S$ as
 follows 
$$
  Y_L= Y_{k1}=Y_{k2} K_{k2}= Y_{k3}K_{k3}K_{k2}=...$$
$$ 
  ...= Y_{1k} K_{1k}K_{1,k-1}K_{k,k-1} K_{k,k-2}...K_{k3} K_{k2}
\equiv Y_R S$$
Then 
\be 
  S=  K_{1k}K_{1,k-1}K_{k,k-1} K_{k,k-2}...K_{k3} K_{k2}  
\label{facto}
\ee
We  observe that, being the first row of $Y(z)$ equat to  $z^{k-1 \over 2}
\Phi(z)$, the following holds:
$$ 
   \Phi_{j+1}(z) = \Phi_j K_j
$$

\vskip 0.3 cm
{\bf Remark 3:} The Stokes' factors of the system (\ref{10}) and of the
gauge-equivalent system (\ref{20}) are the same. From the
skew-symmetry of $V$ it follows that $K_{ji}=(K_{ij}^{-1})^T$. 
\vskip 0.3 cm

Before computing the Stokes factors explicitely, we show that just two
of them are enough to compute all the others. Let
$$ 
  F(z):= \left(  {1 \over \sqrt{k}} {1 \over
                                             z^{k-1 \over 2} }
                         \exp(k z),~ 
                  {1 \over \sqrt{k}} {e^{{i \pi \over k}} \over
                                             z^{k-1 \over 2} }
                         \exp(k e^{2\pi i \over k} z),~
                 ...,~
                {1 \over \sqrt{k}} {e^{{i \pi \over k}(k-1)} \over
                                             z^{k-1 \over 2} }
                         \exp(k e^{{2\pi i \over k}(k-1)} z)
\right)
$$ 
$$ = (F^{(1)}(z),~F^{(2)}(z),~...,~F^{(k)}(z))$$
be the row vector whose entries are the first terms of the asymptotic
expansions of an actual solution $\Phi(z)$ of the generalized
hypergeometric equation. By a straightforward computation  we see that
$$ 
     F(z e^{2\pi i \over k })= F(z) T_F~,~~~~
T_F=\pmatrix{ 0 &     &   &    &...& 1 \cr
              -1&  0  &   &     &&    \cr
                &  -1 & 0 &     &&   \cr
                &     & -1 &    &&        \cr
                &     &    & \ddots&\ddots& \vdots\cr             
                &     &    &      &    -1 &0\cr   
}
$$
We use now the convention of enumeration of Stokes' rays 
$R_0$, $R_1$, ... starting
from $l$ (see above). Let $\Phi_m(z)$ be an actual solution of the
hypergeometric equation having asymptotic behaviour $F(z)$ on ${\cal
S}_m$. 
$$ 
   \Phi_m(z) \sim F(z) ~~~z \to \infty ~~~ z \in {\cal S}_m 
$$
Then 
$$
\Phi_{m+2}(ze^{2\pi i \over k})\sim F(z e^{2 \pi i \over k}
)= F(z) T_F  ~,~~~~ z \in {\cal S}_m
$$
Namely, 
$$ \Phi_{m+2}(ze^{2\pi i \over k})T_F^{-1}\sim F(z) ~,~~~~~ z \in
{\cal S}_m$$
then, by unicity of actual solutions having asymptotic behaviour $F(z)$ in a
sector wider then $\pi$, we have
$$ 
  \Phi_{m+2}(ze^{2\pi i \over k})= \Phi_m(z) T_F~,~~~~~z \in {\cal S}_m
$$
\vskip 0.2 cm 
\noindent
{\bf Lemma 6:} {\it For any $p \in {\bf Z}$ 
               $$ K_{m+2p}= T_F^{-p}K_m T_F^p$$
}
{\bf Proof: } For $z \in {\cal S}_m \cap {\cal S}_{m+1}$ 
$$ 
\Phi_{m+1}(z) = \Phi_m(z)~ K_m = \Phi_{m+2}(ze^{2\pi i \over
k})~T_F^{-1}K_m =$$ 
$$ = \Phi_{m+3}(ze^{2\pi i \over k})~K_{m+2}^{-1}T_F^{-1}K_m 
 = \Phi_{m+1}(z) ~~T_F~ K_{m+2}^{-1}~T_F^{-1}~K_m $$
Then $K_{m+2}=~T_F^{-1}~K_m~T_F$. By induction we prove the lemma. 

\rightline{$\Box$}

\vskip 0.3 cm
>From the lemma it follows that just two Stokes' factors are enought to
compute all the others. We are ready to give a concise formula for
$S$: 
\vskip 0.3 cm 
\noindent
{\bf Theorem 1: }  {\it Let $l$ be an admissible line (i.e. not
containing Stokes' rays), and let us enumerate the rays in
counter-clockwise order starting from the first one in $\Pi_R$ (which
will be $R_0$, then $R_0$, $R_1$, ..., $R_{k-1}$ are in $\Pi_R$, and
$R_k$, ..., $R_{2k-1}$ are in $\Pi_L$). Then the Stokes' matrix
for (\ref{10}), (\ref{20}), (\ref{30}) and $k>3$  is 
\vskip 0.2 cm
\be
  S = \left\{ \matrix{ 
                       \left(K_0~K_1~T_F^{-1} \right)^{k \over 2}
                       ~T_F^{k\over 2} ~\equiv~ T_F^{k \over 2}~\left(
                                        T_F^{-1}~K_{k-2}~K_{k-1}
                       \right)^{k \over 2} ,  ~~~k \hbox{ even } \cr
\cr
\cr
                      \left(K_0~K_1~T_F^{-1} \right)^{k-1 \over 2}
                       ~K_0~T_F^{k-1\over 2} ~\equiv~ T_F^{k-1 \over
                                                       2}~K_{k-1}\left(
                                        T_F^{-1}~K_{k-2}~K_{k-1}
                       \right)^{k-1 \over 2} ,  ~~~k \hbox{ odd } \cr
         } \right.  
\label{S100}
\ee
}

\vskip 0.3 cm 
\noindent 
{\bf Proof: } $S=K_0~K_1~K_2~...~K_{k-1}$,  $K_{2p}= T_F^{-p}~K_0~
T_F^p$ and $K_{2p+1}= T_F^{-p}~K_1~T_F^p$. Then 
$$ 
S=K_0 K_1 ~( T_F^{-1}K_0T_F)~( T_F^{-1}K_1T_F)~K_4 K_5
~...~K_{k-1}$$
$$ 
  = (K_0 K_1 T_F^{-1})~K_0K_1~T_F~K_4 K_5 ~...~ K_{k-1}$$
$$ 
 = (K_0 K_1 T_F^{-1})~K_0K_1~T_F~( T_F^{-2} K_0 T_F^2)~(T_F^{-2} K_1
 T_F^2)~K_6~...~K_{k-1}
$$ 
$$
 = (K_0 K_1 T_F^{-1})(K_0 K_1 T_F^{-1})~K_0 K_1 T_F^2~K_6~...~ K_{k-1}
$$ 
 Now observe that $K_{k-1}= T_F^{-({k \over 2}-1)}~K_1~T_F^{{k\over 2}
 -1}$ for $k$ even, while  $K_{k-1}= T_F^{-({k-1 \over
 2})}~K_0~T_F^{k-1\over 2}$ 
for $k$ odd. Then, for $k$ even 
$$ 
S = (K_0 K_1 T_F^{-1})^{{k \over 2} -2}~K_0K_1T_F^{{k \over
2}-2}~K_{k-2}K_{k-1}$$
$$
     = (K_0 K_1 T_F^{-1})^{{k \over 2} -2}~K_0K_1T_F^{{k \over
2}-2} ~ T_F^{-({k \over 2}-1)}K_0T_F^{{k\over 2}-1}~ 
T_F^{-({k \over 2}-1)}K_1T_F^{{k\over 2}-1}=  
\left(K_0~K_1~T_F^{-1} \right)^{k \over 2}
                       ~T_F^{k\over 2}
$$
For $k$ odd: 
$$ 
 S= (K_0 K_1 T_F^{-1})^{{k -3\over 2} }~ K_0 K_1 T_F^{k-3 \over 2}
 ~K_{k-1}
$$ 
$$ 
=  (K_0 K_1 T_F^{-1})^{{k -3\over 2} }~ K_0 K_1 T_F^{k-3 \over 2}~
T_F^{-({k-1 \over
 2})}K_0T_F^{k-1\over 2}= \left(K_0~K_1~T_F^{-1} \right)^{k-1 \over 2}
                       ~K_0~T_F^{k-1\over 2}
$$ 
If instead we write the Stokes' factors in term of $K_{k-2}$ and
$K_{k-1}$ we obtain the other two formulas in the same way. 

\rightline{$\Box$}
\vskip 0.3 cm 
\noindent
{\bf Remark 4: } For our particular choiche of $l$, $K_0\equiv K_{1k}$,
$K_1\equiv K_{1,k-1}$, $K_{k-2}\equiv K_{k3}$ and $K_{k-1}\equiv
K_{k2}$.  For $k=3$ 
$$ 
  S=T_FK_{32}~(T_F^{-1}K_{12}K_{32})=K_{13}K_{12}K_{32}$$

\vskip 0.3 cm
It is now worth deriving some properties of the monodromy of $Y(z)$
(for (\ref{10})) and $\Phi(z)$ (for (\ref{30})), 
which will be usefull later. Consider $\Phi_m(z)$ with asymptotic
behaviour $F(z)$ on ${\cal S}_m$. Then 
$$ 
   \Phi_m(z)=  \Phi_{m-2}(z)~K_{m-2}~K_{m-1}  \equiv \Phi_m(ze^{2 \pi i
   \over k})~T_F^{-1} ~K_{m-2}~K_{m-1}
$$ 
On the other hand 
$$  \Phi_m(z)=  \Phi_{m+2}(z)~K_{m+1}^{-1}~K_{m}^{-1}  \equiv 
\Phi_m(ze^{-{2 \pi i
   \over k}})~T_F ~K_{m+1}^{-1}~K_{m}^{-1}
$$ 
This proves the following 
\vskip 0.3 cm 
\noindent 
{\bf Lemma 7: } {\it  The basic solution $\Phi_m(z)$ of the
generalized hypergeometric equation (\ref{30}) with asymptotic
behaviour $F(z)$ on ${\cal S}_m$, satisfies the identity
$$ 
   \Phi_m(z e^{2 \pi i \over k})= \Phi_m(z) ~T_m
$$ 
where 
$$  
      T_m:= ~K_{m-1}^{-1}~K_{m-2}^{-1}~T_F~=~T_F~K_{m+1}^{-1}~K_m^{-1}
$$ 
}
\vskip 0.3 cm
\noindent
{\bf Corollary 1: } { \it The monodromy (at $z=0$) of $\Phi_m(z)$ is    
$$ 
          \Phi_m(z e^{2\pi i})=\Phi_m(z)~(T_m)^k
$$
}
\vskip 0.2 cm
Now, for our particular choiche of the line $l$ and for $m= k$, 
$\Phi_m(z)=\Phi_L(z)$. For the solution $Y_L(z)$ of (\ref{10}), the
relations $Y_R(z)=Y_L(z)~S^{-1}$ ($0<\arg z<{\pi \over k}$), 
$Y_L(z)=Y_R(z e^{-2\pi i})~S^{T}$ ($\pi <\arg z <\pi+ {\pi \over k}$)  
immediately imply 
$$ 
        Y_L(ze^{2 \pi i })= Y_L(z)~S^{-1} ~S^T
$$ 
Recall that the $(n,j)$-th entry of $Y(z)$ is  $ Y_{n,j}(z)\equiv 
                      Y^{(j)}_n(z)= {1 \over k^{n-1}} ~z^{{k-1\over 2}-n+1}~ 
                    (z \partial_z)^{(n-1)}~\varphi^{(j)}(z)$, 
and observe that $(z e^{2\pi i})^{k-1 \over 2}= (-1)^{k-1}
                      ~z^{k-1\over2}$ 
Then, from Corollary 1 we get the following: 

\vskip 0.2 cm 
\noindent
{\bf Corollary 2: } { \it Let $T$ be the $k$-monodromy matrix of
$\Phi_L$ (namely, $T\equiv T_k$ for our choiche of $l$). Then 
$$ 
   T^k = (-1)^{k-1}~S^{-1}~S^T
$$ 
}

 Our formula (\ref{S100}) allows us to easily compute $S$. 
The recipe is simply to take $K_{k3}$, $K_{k2}$ (which we are going to
compute explicitely) and substitute them into 
\be
S= T_F^{k \over
 2}~(T_F^{-1}~K_{k3}~K_{k2})^{k \over 2}=T_F^{k \over 2} ~T^{-{k \over
 2}}
\label{S150}
\ee
 for $k$ even, or into 
\be
S= 
  T_F^{k -1\over
 2}~K_{k2}~(T_F^{-1}~K_{k3}~K_{k2})^{k-1 \over 2}=  T_F^{k -1\over
 2}~
K_{k2}~T^{-{k-1 \over 2}}
\label{S200}
\ee
for $k$ odd.

\vskip 0.3 cm
{\bf Computation of Stokes' factors: } We need to distinguish between
$k$ odd and even. 
In the following  $g(z)$ will mean  $g^{(n(k))}(z)$ ($n(k)={k\over
2}+1$ or ${k+1 \over 2}$ for $k$ even or odd respectively). 

\vskip 0.2 cm 
\noindent
{\bf $k$ odd:} 
         $$ 
              g(z)= \varphi^{({k+1 \over 2})}_{L}(z)\equiv 
                       \varphi^{({k+1 \over 2})}_{R}(z)$$
with asymptotic behaviour on 
$$ 
           - \pi < \arg z < \pi +{2 \pi \over k}
$$
If we  iterate the map  $ z \mapsto z e^{ {2\pi i \over k}} $  for
$m=1, 2, ..., {k+1 \over 2}$ times, the domain of $g(z e^{{2\pi i\over
k} m})$ for each $m$ covers ${\cal S}_R$. When we reach  
$m= {k+1 \over 2}$  a new iteration
(i.e. a new rotation of the domain of $-{2 \pi \over k}$) will live
the sector $ - {\pi \over k} < \arg z < {\pi \over k} $ of ${\cal S}_R$
uncovered. The same, if we do   $ z \mapsto z e^{ -{2\pi i  \over k}}
$ the sector $-\pi< \arg z< -\pi +{2 \pi \over k}$ of ${\cal S}_R$
remains uncovered. Then, by remark 2: 
$$ 
  \Phi_{1k}(z)\equiv \Phi_R(z)= 
\Big((-1)^{k-1 \over 2} g(z e^{i{2 \pi \over k}\left( {k+1
  \over 2}\right)}), ~...\hbox{ ${k-3 \over 2}$ unknown terms }...~,
$$
$$
  g(z),~-g(z e^{2 \pi i \over k}), ~g(z e^{4 \pi i\over
k}),~...~,~(-1)^{k-1 
\over 2}
 g(z e^{i{2 \pi \over k}\left( {k-1
  \over 2}\right)})  \Big)
$$
In the same way we see that 
$$ 
\Phi_{k1}(z)\equiv \Phi_L(z)=\Big((-1)^{k-1 \over 2} g(z e^{-i{2 \pi
\over k}
\left( {k-1
  \over 2}\right)}),~-(-1)^{k-1 \over 2} g(z e^{-i{2 \pi \over k}\left( {k-3
  \over 2}\right)}),~...~,-g(ze^{-{2 \pi i \over k }}),
$$ 
 $$
      g(z),~...\hbox{ ${k-1 \over 2}$ unknown terms }... \Big)
$$
and similar expressions for $\Phi_{1,k-1}$, $\Phi_{k,k-1}$,
$\Phi_{k,k-2}$, ..., $\Phi_{k3}$,$\Phi_{k2}$. The unknown terms are 
computed using  the identity
(\ref{60}) and simple considerations on the dominant exponentials 
$|e^{zu_i}|$ on the sectors which remain uncovered in the iterations of 
$ z \mapsto z e^{\pm i{2 \pi \over k}}$. 
\vskip 0.2 cm 
\noindent
{\bf Example:}  
 A simple example will clarify this procedure. Let  $k=7$; then
$g=\varphi^{(4)}$,
$$
 \Phi_{17}=\Phi_R= ( -g(ze^{8 \pi i \over 7}),~?,~?,~g(z),~ -g(ze^{2\pi i
\over 7}),~ g(ze^{4 \pi i \over 7}),~ -g(ze^{6 \pi i \over 7}))
$$
We look for $\varphi^{(3)}_R$. If we take $(-1)g(z e^{-2 \pi i \over
k})$ we miss to cover $-\pi <\arg z< - \pi + {2 \pi \over 7}$ in
${\cal S}_R$. 
On   $-\pi <\arg z< - \pi + { \pi \over 7}$ (between two nearby
Stokes' rays) we have the relations
$|e^{zu_4}|>|e^{zu_5}|>|e^{z u_3}|$. On $-\pi+{ \pi \over 7} <\arg z<
 - \pi + {2 \pi \over 7}$, $|e^{zu_4}|>|e^{zu_3}|$ (later on we will 
simply write $4>3$). Then 
$$ 
  \varphi^{(3)}_R(z)= (-1)g(z e^{-2 \pi i \over
7}) +c_4 \varphi^{(4)}_R(z) +c_5\varphi^{(5)}_R(z)
$$
To find $c_4$, $c_5$ we need another representation for
$\varphi^{(3)}_R$. We consider $ (-1)g(z e^{12 \pi i \over
7})$, which has the correct asymptotic behaviour, but on a domain
which leaves uncovered $ -{3 \pi \over 7}< \arg z < {\pi \over 7}$. 
The relations are: on  $ 0< \arg z < {\pi \over 7}$,
$1>7>2>6>3$; on   $ -{ \pi \over 7}< \arg z < 0$, $1>2>7>3$; on $ -{2
\pi \over 7}< \arg z <- {\pi \over 7}$, $2>1>3$; on  $ -{3 \pi \over
7}<
 \arg z < -{2\pi \over 7}$, $2>3$. Then 
$$
     \varphi^{(3)}_R(z)= (-1)g(z e^{12 \pi i \over
7})+ d_1 \varphi^{(1)}_R(z)+ d_2 \varphi^{(2)}_R(z)+
 d_6 \varphi^{(6)}_R(z)+ d_7 \varphi^{(7)}_R(z)
$$
In the same way one finds 
$$ 
\varphi^{(2)}_R(z)= \left\{ \matrix{
g(z e^{-4 \pi i \over
7}) + a_3 \varphi^{(3)}_R(z)+ a_4 \varphi^{(4)}_R(z)
+ a_5 \varphi^{(5)}_R(z)+ a_6 \varphi^{(6)}_R(z)
                                                  \cr 
g(z e^{10 \pi i \over
7})+ b_1\varphi^{(1)}_R(z)+b_2\varphi^{(4)}_R(z)
                                                 \cr
} 
\right.
$$
$ \varphi^{(i)}$'s are known for $i=1,4,5,6,7$. Using the identity 
(\ref{60}) we compute $a$, $b$, $c$, $d$. We get 
$$ 
  \varphi^{(2)}_R(z)=  g(z e^{-4 \pi i \over
7})
-\bin{7}{1}g(z e^{-2 \pi i \over
7})
+\bin{7}{2}g(z)-\bin{7}{3}g(z e^{2 \pi i \over
7})
+\bin{7}{4} 
g(z e^{4 \pi i \over
7})
$$
$$
\varphi^{(3)}_R(z)=-g(z e^{-2 \pi i \over7})
+
\bin{7}{1} g(z )- \bin{7}{2} g(z e^{2 \pi i \over7})
$$
A similar computation gives $\Phi_{71}=\Phi_L$. 
{\small
$$ \Phi_{71}^T= \left[ \matrix{ 
                               -g(ze^{-{6\pi i \over 7}}) \cr
                               g(ze^{-{4\pi i \over 7}}) \cr
-g(ze^{{-2\pi i \over 7}}) \cr
 g(z)    \cr
-g(ze^{2\pi i \over 7})+\bin{7}{6} g(z) \cr
g(ze^{4\pi i \over 7}) -\bin{7}{6} g(ze^{2\pi i \over
7})+\bin{7}{5}g(z) -\bin{7}{4}g(ze^{-{2 \pi i \over 7}}) \cr
- g(ze^{{6 \pi i \over 7}}) +\bin{7}{6}g(ze^{{4 \pi i \over 7}})
-\bin{7}{5}
  g(ze^{{2 \pi i \over 7}})+\bin{7}{4}g(z) -\bin{7}{3} g(ze^{-{2 \pi i
\over 7}})+
\bin{7}{2} g(ze^{-{4 \pi i \over 7}})\cr    
      }\right]
$$
}
and
{\small 
$$\Phi_{72}^T= \left[ \matrix{
                                 -g(ze^{-{6\pi i \over 7}}) \cr
                               g(ze^{-{4\pi i \over 7}}) \cr
-g(ze^{{-2\pi i \over 7}}) \cr
 g(z)    \cr
-g(ze^{2\pi i \over 7})\cr
g(ze^{4\pi i \over 7}) -\bin{7}{6} g(ze^{2\pi i \over
7})+\bin{7}{5}g(z)  \cr
- g(ze^{{6 \pi i \over 7}}) +\bin{7}{6}g(ze^{{4 \pi i \over 7}})
-\bin{7}{5}
  g(ze^{{2 \pi i \over 7}})+\bin{7}{4}g(z) -\bin{7}{3} g(ze^{-{2 \pi i
\over 7}})\cr 
                     } \right]
$$}
Notice that in each of the last three entries of $\Phi_{72}$ there is
a term missing w.r.t the corresponding entries of $\Phi_{71}$. This
immediately implies
{\tiny
$$ 
    \hbox{\small $K_{72}$} = \pmatrix{ 
                         1 & 0 & 0 & 0 & 0 & 0 & 0 \cr
                         0 & 1 & 0 & 0 & 0 & 0 & \bin{7}{2} \cr
                         0 & 0 & 1 & 0 & 0 & \bin{7}{4}& 0 \cr
                         0 & 0 & 0 & 1 & \bin{7}{6} & 0 & 0 \cr
                         0 & 0 & 0 & 0 & 1 & 0 & 0 \cr
                         0 & 0 & 0 & 0 & 0 & 1 & 0 \cr
                         0 & 0 & 0 & 0 & 0 & 0 & 1 \cr
 } 
$$}
The next step is the computation of $\Phi_{73}$ and $K_{73}$,
through $\Phi_{72}=\Phi_{73}K_{73}$. It is done in the
same way...

\vskip 0.2 cm
 The above procedure is extended to the
 general case. In Appendix 2 we give, for example, the general
 expressions of $\Phi_R$ and $\Phi_L$. The factors of interest are:

\vskip 0.3 cm
{\bf $k$ odd:} 

\vskip 0.2 cm 
\noindent
$(K_{k2})_{j,~k-j+2}=\bin{k}{2(j-1)}$ for $j=2,...,{k+1\over 2}$. 
$(K_{k2})_{j,j}=1$ for $j=1,...,k$. All the other entries are zero. 
\vskip 0.2 cm 
\noindent
$(K_{k3})_{2,1}=-\bin{k}{1}$; $(K_{k3})_{j,~k-j+3}=\bin{k}{2j-3}$ for 
$j=3,...,{k+1 \over 2}$. $(K_{k3})_{j,j}=1$ for  $j=1,...,k$. 
 All the other entries are zero. Namely:
\vskip 0.3 cm 

{\tiny 
$$ \hbox{\small $K_{k2}$}=\pmatrix{
1& &&&&&&&& &0 \cr
& 1&&&&&&&& &\bin{k}{2}\cr
&&1&&&&&&&\bin{k}{4}&\cr
&&&1&&&&&\bin{k}{6}&&\cr
&&&&\ddots&&&\adots&&&\cr
&&&&&1&\bin{k}{k-1}&&&& \cr
\cr
&&&&&&1&&&\cr
\cr
 &    &&&&&&\ddots&&\cr
  &     &&&&&&&1&\cr
   &     &&&&&&&&1\cr
    &    &&&&&&&&&1\cr}
$$
$$
\hbox{\small $K_{k3}$}=\pmatrix{ 
1&&&&&&&&&&&&0 \cr 
\cr
-\bin{k}{1}&1&&&&&&&&&& &0\cr
&&1&&&&&&& &&&\bin{k}{3}\cr
&&&1&&&&&&&&\bin{k}{5}&\cr
&&&&1&&&&&&\bin{k}{7}&&\cr
&&&&&\ddots&&&&\adots&&&\cr
&&&&&&1&0&\bin{k}{k-2}&&&& \cr
\cr
&&&&&&&1&0 \cr
\cr
&&&&&&&&1\cr
\cr
   & && &&&&&&\ddots&&\cr
    &  && &&&&&&&1&\cr
     &   &&&&&&&&&&1\cr
      &  &&&&&&&&&&&1\cr}
$$
}

\vskip 0.3 cm
{\bf $k$ even: }

\vskip 0.2 cm
\noindent
$(K_{k2})_{j,~k-j+2}=\bin{k}{2(j-1)}$ for $j=2,...,{k\over 2}$. 
$(K_{k2})_{j,j}=1$ for $j=1,...,k$. All the other entries are zero. 
\vskip 0.2 cm 
\noindent
$(K_{k3})_{2,1}=-\bin{k}{1}$; $(K_{k3})_{j,~k-j+3}=\bin{k}{2j-3}$, for 
$j=3,...,{k \over 2}+1$.  $(K_{k3})_{j,j}=1$ for  $j=1,...,k$. 
All the other entries are zero. Namely
\vskip 0.3 cm

{\tiny
$$
\hbox{\small $K_{k2}$}=\pmatrix{ 
1&&&&&&&&&&&0 \cr 
\cr
&1&&&&&&& &&&\bin{k}{2}\cr
&&1&&&&&&&&\bin{k}{4}&\cr
&&&1&&&&&&\bin{k}{6}&&\cr
&&&&\ddots&&&&\adots&&&\cr
&&&&&1&0&\bin{k}{k-2}&&&& \cr
\cr
&&&&&&1&0 \cr
\cr
&&&&&&&1\cr
\cr
    && &&&&&&\ddots&&\cr
      && &&&&&&&1&\cr
        &&&&&&&&&&1\cr
        &&&&&&&&&&&1\cr}
$$
$$ \hbox{\small $K_{k3}$}=\pmatrix{
1& &&&&&&&& &&0 \cr
\cr
-\bin{k}{1}&1& &&&&&&&& &0\cr
&& 1&&&&&&&& &\bin{k}{3}\cr
&&&1&&&&&&&\bin{k}{5}&\cr
&&&&1&&&&&\bin{k}{7}&&\cr
&&&&&\ddots&&&\adots&&&\cr
&&&&&&1&\bin{k}{k-1}&&&& \cr
\cr
&&&&&&&1&&&\cr
\cr
 &    &&&&&&&\ddots&&\cr
  &    & &&&&&&&1&\cr
   &    & &&&&&&&&1\cr
    &    &&&&&&&&&&1\cr}
$$
}


\vskip 0.3 cm 
\section{ Reduction of $S$ to ``Canonical'' Form }

Some examples of computations of $S$ are in Appendix 2. 
The reader may observe that $S$ is not in a nice upper triangular
form (see also Lemma 4) and quite strange numbers
(complicated combinations of sum and products of binomial
coefficients) appear. 

Some natural operations are allowed on the Stokes matrices of a
Frobenius manifold:

\vskip 0.2 cm
a) Permutations. Let us consider the system 
$$ 
 { d \tilde{Y} \over dz}= \left[ U + {1 \over z} V  \right]~\tilde{Y} 
$$ 
where $U = $ diag$( u_1,~u_2,~...,~u_k)$. Let $\sigma:~(1,2,..,k)$ 
$\mapsto$ $(\sigma(1),\sigma(2),...,\sigma(k))$ be a permutation. It 
is represented by an invertible matrix $P$ such that 
$$ P~U~P^{-1}~=~\hbox{
diag}(u_{\sigma(1)},u_{\sigma(2)},...,u_{\sigma(k)})$$ 
The system 
$$ {d {\cal  Y} \over dz}= \left[PUP^{-1}+{1 \over z} PVP^{-1} \right]
~{\cal Y}
$$ 
has solutions 
$${\cal Y}_{L/R}(z) := P~\tilde{Y}_{L/R}(z)~P^{-1} \sim 
( I+O({1\over z}))~e^{z PUP^{-1}}~~~~{\cal Y}_{L}(z)={\cal Y}_{R}(z)
~PSP^{-1}$$
$$
{\cal Y}_0(z):= P \tilde{Y}_0(z)=(PX^{-1}+O(z))z^{\mu}z^R,~~~~{\cal
Y}_0(z)= {\cal Y}_R(z) ~PC
$$
$S$ and $C$ are then transformed in $PSP^{-1}$ and $PC$. 
For a suitable $P$, $PSP^{-1}$ is  upper triangular. 
As a general result \cite{BJL1}, the good permutation is the
one which put $u_{\sigma(1)}$, ..., $u_{\sigma(k)}$ in lexicographical
ordering w.r.t. the oriented line $l$. The effect
of the permutation $P$ corresponds to a change of coordinates in the 
given chart, consisting in the permutation $\sigma$ of the
coordinates. 

\vskip 0.2 cm
b) Sign changes  of the entries. The construction at point a) is
repeated, but now a diagonal matrix ${\cal I}$ with $1$'s or $-1$'s on the diagonal
takes the place of $P$. In ${\cal I} S {\cal I}^{-1}$ some entries
change  sign. Note that ${\cal I}U{\cal I}^{-1}\equiv U$. 
Moreover, the formulae \cite{Dub1} \cite{Dub2}  
which define a local chart of the manifold in
terms of monodromy data are not affected by $S\mapsto {\cal I}S{\cal
I}$, $C \mapsto {\cal I} C$.  

 \vskip 0.2 cm 
c) Action of the braid group. 

\noindent
We first recall that the  braid group is generated by $k-1$ elementary 
braids $\beta_{12}$, $\beta_{23}$, ..., $\beta_{k-1,k}$,  with
relations: 
$$ 
    \beta_{i,i+1}\beta_{j,j+1}=\beta_{j,j+1}\beta_{i,i+1}~\hbox{ for
} i+1 \neq j,~j+1\neq i$$
$$ 
 \beta_{i,i+1}\beta_{i+1,i+2}\beta_{i,i+1}=
\beta_{i+1,i+2}\beta_{i,i+1} 
     \beta_{i+1,i+2} $$ 
This abstract group is realized as the fundamental group of (${\bf C}^k
\backslash$diagonals)$/{\cal S}_k$ $:= \{(u_1,...,u_k)~|~u_i\neq u_j 
\hbox{ for } i\neq j \}/{\cal S}_k$, where ${\cal S}_k$ is the
symmetric group of order $k$. 

 In section 4 we proved that the Stokes' rays of the systems
 (\ref{tuono}) (\ref{tuono1}), and more generally the rays of the
 system (\ref{zzzz1}), depend on the point $t$ of the manifold. This is
 equivalent to the fact that the eigenvalues $u_1(t)$, ..., $u_k(t)$
 of  $\hat{\cal U}(t)$ depend on $t$.
When $t$ changes, $u_1(t)$, ..., $u_k(t)$ change their position in the
 complex plane. Consequently, Stokes' rays move. Let us start from the
 point $t_0=(0,t_0^2,...,0)$. If we move in a sufficently small
 neighbourhood of $t_0$, the rays slightly change their positions. But
 if we move sufficently far away from $t_0$, some Stokes' rays cross the
 {\it fixed} admissible line $l_{t_0}$. Then, we must change
 ``Left''  and ``Right'' solutions of
 (\ref{zzzz1}). Then also $S$ and $C$ change.   This is the reason why
 the monodromy data $S$ and $C$ 
change when we go from one local chart to another. 

The motions of the points $u_1(t)$, ..., $u_k(t)$  as $t$ changes 
represent transformations of
the braid group. Actually, a braid $\beta_{i,i+1}$ can be represented  as an ``elementary'' 
 deformation 
consisting  of  a permutation of $u_i$, $u_{i+1}$ moving 
counter-clockwise (clockwise or counter-clockwise is a matter of
convention).

 Suppose $u_1$, ..., $u_k$ are already  in
lexicographical order w.r.t. $l$, so that $S$ is upper
triangular (recall that this configuration can be reached by a
suitable permutation $P$). The effect on $S$ 
of the deformation of $u_i$,
$u_{i+1}$ representing $\beta_{i,i+1}$ is the following : 
$$ 
   S \mapsto~S^{\beta_{i,i+1}}:=A^{\beta_{i,i+1}}(S)~S~A^{\beta_{i,i+1}}(S)
$$
where 
$$ 
   \left(A^{\beta_{i,i+1}}(S) \right)_{nn}=
   1~~~~~~n=1,...,k~~~n\neq~i,~i+1
$$
$$
 \left(A^{\beta_{i,i+1}}(S) \right)_{i+1,i+1}=-s_{i,i+1}
$$
$$
 \left(A^{\beta_{i,i+1}}(S) \right)_{i,i+1}=
 \left(A^{\beta_{i,i+1}}(S) \right)_{i+1,i}=1
$$
and all the other entries are zero. 

  For the inverse braid $\beta_{i,i+1}^{-1}$ ($u_i$ and $u_{i+1}$ move
  clockwise) the representation is
 $$ 
   \left(A^{\beta_{i,i+1}^{-1}}(S) \right)_{nn}=
   1~~~~~~n=1,...,k~~~n\neq~i,~i+1
$$
$$
 \left(A^{\beta_{i,i+1}^{-1}}(S) \right)_{i,i}=-s_{i,i+1}
$$
$$
 \left(A^{\beta_{i,i+1}^{-1}}(S) \right)_{i,i+1}=
 \left(A^{\beta_{i,i+1}^{-1}}(S) \right)_{i+1,i}=1
$$
and all the other entries are zero.

 For a generic braid $\beta$, which is a product of $N$ elementary braids
(for some $N$) $\beta=\beta_{j_1,j_1+1}...\beta_{j_N,j_N+1}$, the
 action is 
$$ 
  S\mapsto S^{\beta}=
  A^{\beta}(S)~S~\left[A^{\beta}(S)\right]^T
$$
where 
$$ 
 A^{\beta}(S)=
 A^{\beta_{i_N,i_N+1}}(S^{
 \beta_{j_1,j_1+1}...\beta_{j_{N-1},j_{N-1}+1} })
~...~  A^{\beta_{i_2,i_2+1}}(S^{\beta_{i_1,i_1+1}}) A^{\beta_{i_1,i_1+1}}(S)
$$
We remark that $S^{\beta}$ is still upper triangular.

\vskip 0.3 cm 
   In figure 4 we have drawn some lines  $L_j=\{ \lambda= u_j + \rho
   e^{i({\pi \over 2} -\epsilon)},~~\rho>0\}$ ($0<\epsilon<{\pi
   \over k}$ is the
   angle of $l$), which help us to visualize the
   topological effect of the braids action ( they are the branch cuts for the  fuchsian system which will be
   introduced in section 8). We are going to prove that the  braid whose effect is to set the deformed
   points in cyclic order and the cuts in the configuration of figure 4
   (namely, the last two
   cuts remain unchanged, the others are alternatively
   ``inverted''), brings $S$ in a {\it canonical form}: 
$s_{i,i+1}=\bin{k}{1}$ $\forall
   i=1,...,k-2$; $ s_{i,i+2}=\bin{k}{2}$ $\forall
   i=1,...,k-3$; ...; $ s_{i,i+n}=\bin{k}{n}$ $\forall 
i=1,...,k-n+1$;  $ s_{i,k}=-\bin{k}{k-i}$. Namely 
 {\small 
\be
S^{\beta}= \pmatrix{ 1&\bin{k}{1}&\bin{k}{2}&\bin{k}{3}&\bin{k}{4}&
...& -\bin{k}{k-1}\cr
 & 1 & \bin{k}{1} & \bin{k}{2} & \bin{k}{3}& ... & -\bin{k}{k-2}\cr
 &   & 1          &\bin{k}{1} & \bin{k}{2} & ... & - *           \cr
 &   &            & 1          & \bin{k}{1}&... & - *         \cr
&&&&&& \cr
 &   &            &            &    1      &    &  - *       \cr
 &&&&&\ddots & \cr
&&&&&&1\cr
}
\label{SSCAN}
\ee
}
 Note that the last column is negative. Its sign is inverted by $S
\mapsto {\cal I} S {\cal I}$, where ${\cal I}:=$ diag$(1,1,...,-1)$. 
\begin{figure}
\epsfxsize=14cm
\centerline{
\epsffile{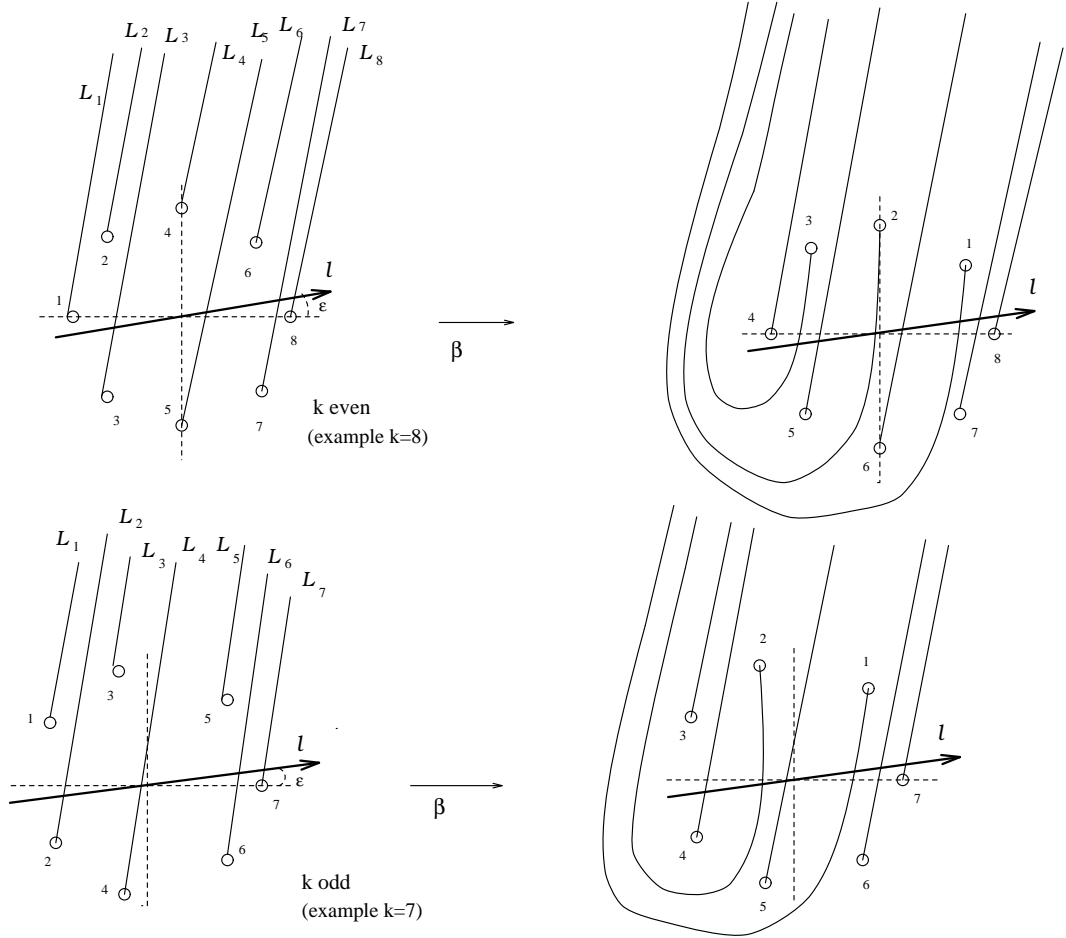}}
\caption{ Effect of the braid which brings $S$ to the canonical form on
the lines $L_j$}
\end{figure}

\begin{figure}
\epsfxsize=15cm
\epsffile{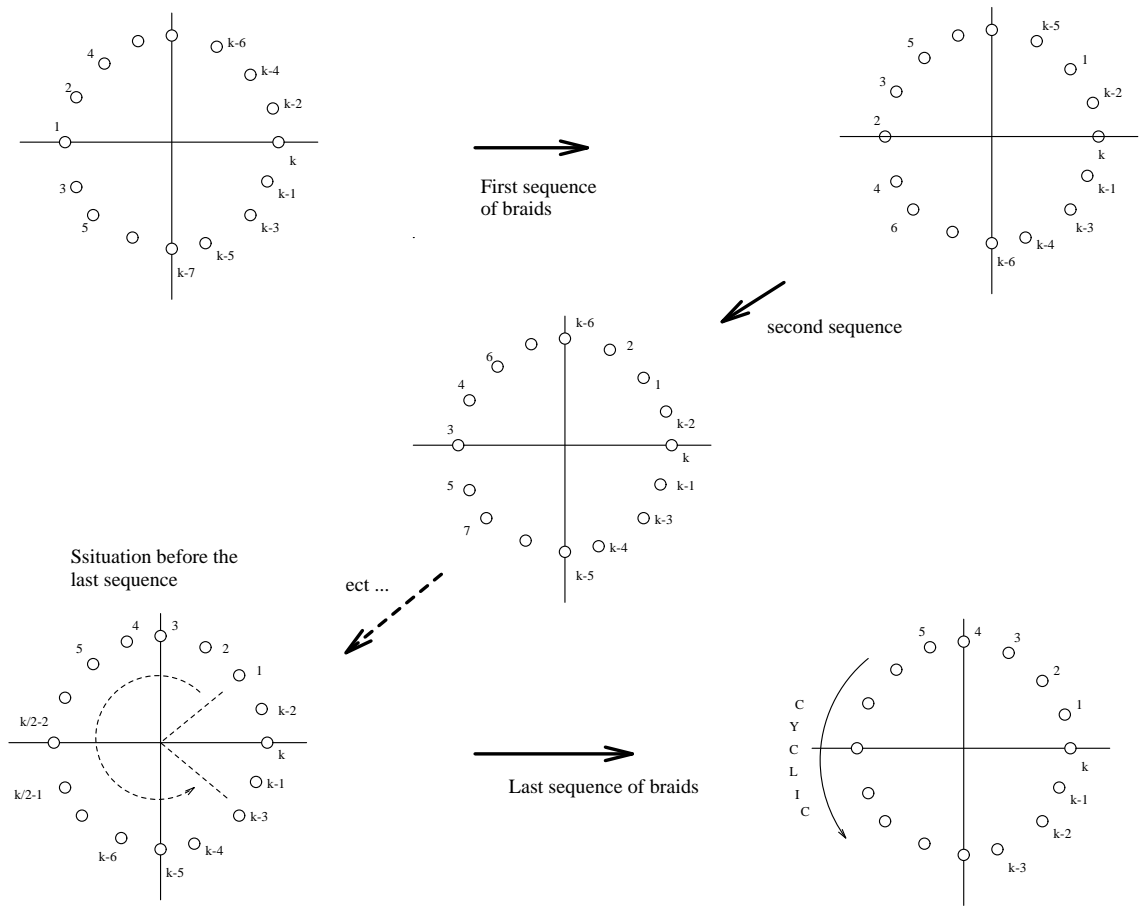}
\caption{ Effect of the sequences of braids which bring $S$ to the
canonical form (the figure refers to $k$ even)}
\end{figure}

\vskip 0.2 cm
\noindent
{\bf Lemma 8:} {\it Let the points $u_j$ ($j=1,..,k$) be in
lexicographical order w.r.t the admissible line $l$. Then the braid 
$$ 
\beta:=
(\beta_{k-5,k-4}~\beta_{k-6,k-5}~...~\beta_{12})~(\beta_{k-6,k-5}~\beta_{k-7,k-6}
~... ~\beta_{23})~(\beta_{k-7,k-6}~...~ \beta_{34})~...$$ $$...~\beta_{{k \over
2}-2,{k \over 2}-1}~
(\beta_{k-3,k-2}~\beta_{k-4,k-3}~...~\beta_{12})$$
for $k$ even, or 
$$ 
\beta:=
(\beta_{k-5,k-4}~\beta_{k-6,k-5}~...~\beta_{12})~(\beta_{k-6,k-5}~\beta_{k-7,k-6}
~...~ \beta_{23})~(\beta_{k-7,k-6}~...~ \beta_{34})~...$$ $$...~(\beta_{{k-3 \over
2},{k-1 \over 2}}~\beta_{{k-5 \over 2},{k-3 \over 2}})~
(\beta_{k-3,k-2}~\beta_{k-4,k-3}~...~\beta_{12})$$
for $k$ odd, brings the points in cyclic counter-clockwise order, 
 $u_1$ being the first point  in $\Pi_L$ 
(figure 4, right side, or figure 5). 
}
\vskip 0.2 cm
Note that we have collected the braids in ${k\over 2}-1$ ($k$ even), or
${k-3 \over 2}$ ($k$ odd) sequences $(...)$. 

\vskip 0.2 cm
\noindent
 {\bf Proof: } Let $k$ be even. The first braid 
$\beta_{k-5,k-4}$ iterchanges $u_{k-4}$ and $u_{k-5}$. The second
 braid interchanges $u_{k-5}$ and $u_{k-6}$. One easily sees that the
 effect of the  first sequence of braids 
$(\beta_{k-5,k-4}~\beta_{k-6,k-5}~...~\beta_{12})$ is to bring $u_1$
 in  the (old) position of $u_{k-4}$, $u_{k-4}$ in the position of
 $u_{k-5}$, $u_{k-5}$ in the position of $u_{k-6}$, ..., $u_4$ in the
 position of $u_3$ and $u_2$ in the position of $u_1$ (figure
 5). $u_{k}$, $u_{k-1}$, $u_{k-2}$, $u_{k-3}$ are not moved. 

The second sequence of braids $(\beta_{k-6,k-5}~\beta_{k-7,k-6}
~... ~\beta_{23})$ acts in a similar way, bringing $u_2$ in $u_{k-5}$,
$u_{k-5}$ in $u_{k-6}$, ..., $u_3$ in $u_2$.  $u_{k}$, 
$u_{k-1}$, $u_{k-2}$, $u_{k-3}$, $u_{k-4}$ are not moved.

 We go on in this way. After the
 action of 
$$(\beta_{k-5,k-4}~\beta_{k-6,k-5}~...~\beta_{12})~(\beta_{k-6,k-5}~\beta_{k-7,k-6}
~... ~\beta_{23})~(\beta_{k-7,k-6}~...~ \beta_{34})~...~\beta_{{k \over
2}-2,{k \over 2}-1}~$$ the points are as in figure 5: $u_k$ is on the
positive real axis, $u_{k-2}$ is the first point met in
counter-clockwise order, $u_1$ is the second, $u_2$ is the third; the
points are in cyclic order up to $u_{k-3}$; finally, $u_{k-1}$ is the
last point before reaching again the positive real axis from below. 

 Then, $ 
(\beta_{k-3,k-2}~\beta_{k-4,k-3}~...~\beta_{12})$ brings $u_1$ in
$u_{k-2}$, $u_{k-2}$ in $u_{k-3}$, $u_{k-3}$ in $u_{k-4}$, and so
on. The cyclic order is reached. 

For $k$ odd the proof is similar. 

\rightline{$\Box$}

\vskip 0.2 cm
 
A careful consideration of the effect of the braid $\beta$ on the
lines $L_j$ (which we leave as an exercise for the reader) shows that
they are alternatively inverted as in figure 4. To reconstruct
uniquely this configuration we just need to know the oriented line
$l$, namely,  its angle $\epsilon$ w.r.t the positive real axis. 
The points $u_{k-1}$, $u_k$ and 
the lines  $L_{k-1}$ and $L_k$ are unchanged 
(angle ${\pi \over 2}-\epsilon$). The line at 
$u_1$ starts in the opposite  direction,  it goes
around $u_2,...,u_{k-2}$ without intesecting other cuts, and then goes
to $\infty$ with the original asymptotic direction ${\pi \over
2}-\epsilon$. 
Moving in the
direction opposite to that of $l$ we meet $u_{k-2}$. Its line 
has the original direction  ${\pi \over 2}
-\epsilon$.  
Then we meet 
$u_2$, and the corresponding line starts with opposite direction, goes
arouns $u_3,...,u_{k-3}$ and then goes to $\infty$ with asymptotic
direction ${\pi \over 2}-\epsilon$. And so on.

\vskip 0.2 cm 
Now we find the matrix representation for the braid $\beta$. 

\vskip 0.2 cm
\noindent 
{\bf Proposition 1:} { \it The braid $\beta $ of Lemma 8 has the following
matrix representation: 
{\tiny
$$ 
  \hbox{ \small $A^{\beta}(S)$}
               = \pmatrix{0&0&0&0&0&&&&0&0&0&0& 1&0 &0 \cr 
                          &&&&&&&&0&0& 1&0& \bin{k}{1}& 0&0\cr
                          &&&&&&&&1&0&\bin{k}{1}&0&\bin{k}{2}&0& 0\cr        
                          &&&&&&&&\bin{k}{1}&0&\bin{k}{2} & 0&\bin{k}{3} &0&
                          0\cr
                         &&&&&&&&*&.&*&.&*&.&. \cr
                         &&&&&&&... &*&.&*&.&*&.&. \cr
                         &&&&&&&&*&.&*&.&*&.&. \cr
                          &&&&&&&&*&.&*&.&*&.&. \cr
0 &0&0&1&0&&&&*&.&*&.&*&.&. \cr
0 &1&0&\bin{k}{1}&0&&&...&*&.&*&.&*&.&. \cr 
1&\bin{k}{1}&0&\bin{k}{2}&0&&&...&*&.&*&.&*&.&. \cr
0 &0&1&\bin{k}{3}&0&&&&*&.&*&.&*&.&. \cr
&&&&&&&&&&&& \cr
0 &0&0&0&1&&&&*&.&*&.&*&.&. \cr
 &&&&&&&&*&.&*&.&*&.&. \cr
                         .&.&.&.&.&&&... &*&.&*&.&*&.&. \cr
                         .&.&.&.&.&&&&*&.&*&.&*&.&. \cr
                         .&.&.&.&.&&&&*&.&*&.&*&.&. \cr
&&&&&&&&\bin{k}{k-7}&0&\bin{k}{k-6}&0&\bin{k}{k-5}&0&0 \cr
&&&&&&&&0&1&\bin{k}{k-5}&0&\bin{k}{k-4}&0&0 \cr 
&&&&&&&&0&0&0&1&\bin{k}{k-3}&0&0 \cr
&&&&&&&&&&&& \cr
&&&&&&&&0&0&0&0&0&1&0\cr
0&0&0&0&0&&&&0&0&0&0&0&0&1\cr
}
$$
}
for $k$ even. 
  
{\tiny 
$$ 
  \hbox{\small $A^{\beta}(S)$}=
                \pmatrix{0&0&0&0&0&&&&0&0&0&0& 1&0 &0 \cr 
                          &&&&&&&&0&0& 1&0& \bin{k}{1}& 0&0\cr
                          &&&&&&&&1&0&\bin{k}{1}&0&\bin{k}{2}&0& 0\cr        
                          &&&&&&&&\bin{k}{1}&0&\bin{k}{2} & 0&\bin{k}{3} &0&
                          0\cr
                         &&&&&&&&*&.&*&.&*&.&. \cr
                         &&&&&&&... &*&.&*&.&*&.&. \cr
                         &&&&&&&&*&.&*&.&*&.&. \cr
                          &&&&&&&&*&.&*&.&*&.&. \cr
0 &0&0&0&1&0&&&*&.&*&.&*&.&. \cr
0&0 &1&0&\bin{k}{1}&0&&...&*&.&*&.&*&.&. \cr 
1&0&\bin{k}{1}&0&\bin{k}{2}&0&&...&*&.&*&.&*&.&. \cr
0 &1&\bin{k}{2}&0&\bin{k}{3}&0&&&*&.&*&.&*&.&. \cr
0 &0&0&1&\bin{k}{4}&0&&&*&.&*&.&*&.&. \cr
&&&&&&&&&&&&&\cr
0 &0&0&0&0&1&&&*&.&*&.&*&.&. \cr
                         .&.&.&.&.&&&... &*&.&*&.&*&.&. \cr
                         .&.&.&.&.&&&&*&.&*&.&*&.&. \cr
                         .&.&.&.&.&&&&*&.&*&.&*&.&. \cr
&&&&&&&&\bin{k}{k-7}&0&\bin{k}{k-6}&0&\bin{k}{k-5}&0&0 \cr
&&&&&&&&0&1&\bin{k}{k-5}&0&\bin{k}{k-4}&0&0 \cr 
&&&&&&&&0&0&0&1&\bin{k}{k-3}&0&0 \cr
&&&&&&&&&&&&&\cr
&&&&&&&&0&0&0&0&0&1&0\cr
0&0&0&0&0&&&&0&0&0&0&0&0&1\cr
}
$$
}
 for $k$ odd. }
  
The $``*''$ means $\bin{k}{j}$, and $j$ increases by
 one when we move downwards row by row.

\vskip 0.3 cm 
\noindent
{\bf Proof:} We need some steps. 
\vskip 0.2 cm
  1)  For an upper triangular Stokes' matrix $S$ (with entries
  $s_{ij}$) the entries of the matrix $A^{\beta}$ which are different
  from zero are the following:

\vskip 0.2 cm
\noindent  
{\bf $k$ odd:}

\vskip 0.2 cm
$$
  A^{\beta}_{k-1,k-1}=A^{\beta}_{k,k}=1
$$
 For  $j$ even, $2+2i\leq j \leq k-2$, and $i=1,2,...,{k\over 2}-1$  
$$
A^{\beta}_{{k\over 2}-i,2i}=1$$
$$
A^{\beta}_{{k\over 2}-i,j}= -s_{2i,j}
+
 \sum_{n=1}^{{j\over 2}-i-1}~(-1)^{n+1} F_{n,i,j}
$$
$$
F_{n,i,j}= 
\sum_{\alpha_1=1}^{{j\over 2}-i-n}  \left(
\sum_{\alpha_2=1}^{{j\over 2}-i-\alpha_1-n+1} \left(
\sum_{\alpha_3=1}^{{j\over 2}-i-\alpha_1-\alpha_2-n+2} \left(
...
\sum_{\alpha_p=1}^{{j\over 2}-i-\sum_{r=1}^{p-1}\alpha_r -n+p-1}
\left(...
\sum_{\alpha_n=1}^{{j\over 2}-i-\sum_{r=1}^{n-1}\alpha_r -1}
f(i,j,\alpha)  \right) \right) \right) \right)
$$
$$
f(i,j,\alpha)=s_{2i,2i+2\alpha_1}~s_{2i+2\alpha_1,2i+2\alpha_1+2\alpha_2}~s_{2i+2\alpha_1+
 2\alpha_2,2i+2\alpha_1+2\alpha_2+2\alpha_3}~...~
s_{2i+2\alpha_1+...+2\alpha_n,j}
 $$
For $j$ even, $2i+2\leq j \leq
k-2$ and $i=0,1,...,{k\over 2}-2$ 
$$
A^{\beta}_{{k\over 2}+i,2i+1}=1$$
$$
A_{{k\over 2}+i,j}=-s_{2i+1,j}+\sum_{n=1}^{{j\over 2}-i-1}(-1)^{n+1}
G_{n,i,j} $$
 $$
G_{n,i,j}= 
\sum_{\alpha_1=1}^{{j\over 2}-i-n}  \left(
\sum_{\alpha_2=1}^{{j\over 2}-i-\alpha_1-n+1} \left(
\sum_{\alpha_3=1}^{{j\over 2}-i-\alpha_1-\alpha_2-n+2} \left(
...
\sum_{\alpha_p=1}^{{j\over 2}-i-\sum_{r=1}^{p-1}\alpha_r -n+p-1}
\left(...
\sum_{\alpha_n=1}^{{j\over 2}-i-\sum_{r=1}^{n-1}\alpha_r -1}
g(i,j,\alpha)  \right) \right) \right) \right)
$$
$$
g(i,j,\alpha)=s_{2i+1,2i+2\alpha_1}~s_{2i+2\alpha_1,2i+2\alpha_1+2\alpha_2}~s_{2i+2\alpha_1+
 2\alpha_2,2i+2\alpha_1+2\alpha_2+2\alpha_3}~...~s_{2i+2\alpha_1+...+2\alpha_n,j}
 $$
More explicitely, 
{\small
\be
  A^{\beta}(S)= \pmatrix{ .&.&.&.&.&.&... \cr
                           .&.&.&.&.&.&... \cr                      
                          .&.&.&.&0&0&... \cr
                         .&.&0&0&0&1&...\cr        
                          .&0&0&1&0&-s_{46}& ...\cr
0&1&0&-s_{24}&0&-s_{26}+s_{24}s_{46}&... \cr 
1&-s_{12}&0&-s_{14}+s_{12}s_{24}&0&-s_{16}+(s_{12}s_{26}+s_{14}s_{46})
-s_{12}s_{24}s_{46}&... \cr
0&0&1&-s_{34}&0&-s_{36}+s_{34}s_{46}&... \cr
. &.&0&0&1&-s_{56}&... \cr
 .&.&.&.&0&0&... \cr
                         .&.&.&.&.&.&... \cr
                           .&.&.&.&.&.&... \cr
}
\label{As1}
\ee
%
%
}

\vskip 0.2 cm
\noindent
 {\bf $k$ odd:}

\vskip 0.2 cm 

$$  A_{k-1,k-1}=A_{k,k}=1$$
For $j$ odd, $2i+1 \leq j \leq k-2$, $i=1,...,{k-3 \over 2}$:
$$ 
  A_{{k-1\over 2},2i}=1
$$
$$
A_{{k-1\over 2}+i,j}=-s_{2i,j}+\sum_{n=1}^{{j-1\over 2}-i} (-1)^{n+1}
H(n,i,j)
$$
$$
H(n,i,j)=\sum_{\alpha_1=0}^{{j-1\over 2}-i-n}\left(\sum_{\alpha_2=1}^{
{j-1\over 2}-i-\alpha_1-n+1} \left( ...\sum_{\alpha_p=1}^{{j-1\over
2}-i -\sum_{\alpha_r=1}^{p-1}\alpha_r-n+p-1} 
\left( ... \sum_{\alpha_n=1}^{{j-1\over
2}-i -\sum_{\alpha_r=1}^{n-1}\alpha_r-1} h(i,j,\alpha) \right)\right)
\right)
$$
$$
h(i,j,\alpha)= s_{2i,2i+1+2\alpha_1}~s_{2i+1+2\alpha_1,
2i+1+2\alpha_1+2\alpha_2 }~....~s_{2i+1+2\alpha_1+...+2\alpha_n,j}
$$
For $j$ odd, $2i+3\leq j \leq k-2$, $i=0,1,...,{k-3\over 2}$:
$$ 
  A_{{k-1\over 2} -i , 2i+1}=1$$
$$
A_{{k-1\over 2}-i,j}=-s_{2i+1,j}+\sum_{n=0}^{{j-1\over 2}-i}(-1)^{n+1}
V(n,i,j)
$$
$$
V(n,i,j)=\sum_{\alpha_1=1}^{{j-1\over 2}-i-n} \left(
\sum_{\alpha_2=1}^{{j-1\over
2}-i-\alpha_1-n+1}\left(...\sum_{\alpha_p=1}^{{j-1\over
2}-i-\sum_{r=1}^{p-1}\alpha_r -n+p-1} \left(...
 \sum_{\alpha_n=1}^{{j-1\over
2}-i-\sum_{r=1}^{n-1}\alpha_r -1} v(i,j,\alpha) \right)\right) \right)
$$
$$
v(i,j,\alpha)=
s_{2i+1,2i+1+2\alpha_1}~s_{2i+1+2\alpha_1,2i+1+2\alpha_1+2\alpha_2}
~...~ 
s_{2i+1+2\alpha_1+...+2\alpha_n,j}
$$
  More explicitely:
{\small
\be
  A^{\beta}(S)= \pmatrix{ .&.&.&.&.&.&. &... \cr
                          .&.&.&.&.&.&0& ... \cr
                         .&.&.&.&0&0&1 &...\cr        
                          .&.&0&0&1&0&-s_{57}& ...\cr
0&0&1&0&-s_{35}&0&-s_{37}+s_{35}s_{57}&... \cr 
1&0&-s_{13}&0&-s_{15}+s_{13}s_{35}&0&-s_{17}+(s_{13}s_{37}+s_{15}s_{57})
-s_{13}s_{35}s_{57}&... \cr
0&1&-s_{23}&0&-s_{25}+s_{23}s_{35}&0&-s_{27}+(s_{23}s_{37}+s_{25}s_{57})
-s_{23}s_{35}s_{57}&... \cr
.&0&0&1&-s_{45}&0&-s_{47}+s_{45}s_{57}&... \cr
.&.&.&0&0&1&-s_{67}&...\cr
 .&.&.&.&.&0&0&... \cr
                         .&.&.&.&.&.&.&... \cr
                           .&.&.&.&.&.&.&... \cr
}
\label{As2}
\ee
}

In order to prove the above expressions, we have to find each matrix
$A^{\beta_{i,i+1}}$ corresponding to the elementary braid $\beta_{i,i+1}$
appearing in $\beta$. This means computing its  entry $(i+1,i+1)$.  
This is a rather complicated problem, since the
entry $(i+1,i+1)$ of a given 
$A^{\beta_{i,i+1}}$ is minus the entry $(i,i+1)$ of the
Stokes' matrix resulting from the action of the elementary braids
acting before $\beta_{i,i+1}$, which in general is a sum of products
of the elements $s_{ij}$ of the initial Stokes matrix $S$. 

  First we recall that  $S\mapsto
  A^{\beta_{i,i+1}}~S~A^{\beta_{i,i+1}}$ has the following effect on
  the entries of $S$:
$$ 
   s_{n,i}\mapsto s_{n,i+1} ,~~~~s_{n,i+1}\mapsto
   s_{n,i}-s_{i,i+1}s_{n,i+1},~~~~ n=1,2,...,i-1 $$
$$s_{i,i+1}\mapsto -s_{i,i+1}$$
$$ s_{i,n}\mapsto s_{i+1,n},~~~~s_{i+1,n}\mapsto
s_{i,n}-s_{i,i+1}s_{i+1,n},
~~~~ n=i+2,...,k $$
 while all the other entries of $S$ remain unchanged.     

    We start from $A^{\beta_{k-5,k-4}}$, whose non trivial entry is simply 
$-s_{k-5,k-4}$. Its action on $S$ brings $s_{k-6,k-5}$ to
$s_{k-6,k-4}$ (the reader may compute all the elements of
$S^{\beta_{k-5,k-4}}$). 
 
 Then, the entry of $A^{\beta_{k-6,k-5}}$ is $-s_{k-6,k-4}$.
 Proceeding in this way, the reader may check that for the first
 sequence of braids $(\beta_{k-5,k-4}\beta_{k-6,k-5}~...~\beta_{12})$
 the entries $(i,i+1)$ of the matrices are: 
$$ -s_{i,k-4}~\hbox{ for } A^{\beta_{i,i+1}}
$$
 Now, observe that 
{\small
$$
\pmatrix{ 1&  \cr
           &\ddots & \cr
           &       &1 \cr
           &       & & 0&1\cr
           &       & & 1 & x_1 \cr 
           &       & & & & & 1&   \cr
           &       & & & & & & 1 \cr
           &       & & & & &  && \ddots  \cr
           &       & & & & &  && & 1 \cr
}
\pmatrix{ 1&  \cr
           &\ddots & \cr
           &       &1 \cr
           &       & &1 \cr
           &&       & & 0&1\cr
           &&       & & 1 & x_2 \cr 
           &&       & & & & & 1&   \cr
           &&       & & & & &  & \ddots  \cr
           &&       & & & & &  & & 1 \cr
}
=\pmatrix{ 1&  \cr
           &\ddots & \cr
           &       &1 \cr
           &       & &0 & 0 & 1 \cr
           &       & &1 & 0 & x_1\cr
           &       & &0 & 1 & x_2 \cr 
           &       & &  &   &    & 1&   \cr
           &       & &  &   &    &  &   \ddots  \cr
           &       & &  &   &    &  &          & 1 \cr
}
 $$
}
and recall that   $A^{\beta_1
 \beta_2}= A^{\beta_2} A^{\beta_1}$. This implies for the first
sequence of braids: 
{\small
$$ m_1:= A^{\beta_{k-5,k-4}~\beta_{k-6,k-5}~...~\beta_{12}}
    \pmatrix{
                     0&       & && & 1 &   &   \cr 
                     1&0      & && & * &   &     \cr
                      & 1     &0&& & * &  &     \cr
                      &       & \ddots&\ddots &&\vdots &      &       \cr 
                      &        &     &  \ddots    &0& *
           & & \cr 
                      &       &      &             &  1   &*& &    \cr
                      &&&&& &1& \cr 
                      &&&&& & &1 \cr
                      &&&&& & & &1 \cr
                      &&&&& & & & &1 \cr
}
$$}
the entries $*$ are exactly those of the $A_{i,i+1}$'s, namely $-s_{1,k-4}$, $-s_{2,k-4}$, ...,
$-s_{k-5,k-4}$ from the top to the bottom of the $(k-4)^{th}$ column.

For the second sequence of braids $(\beta_{k-7,k-6}~...~\beta_{34})$, the entries are 
$$  -s_{i-1,k-6}~\hbox{ for } A^{\beta_{i,i+1}}
$$
and, as above:
{\small   
$$    m_2:=
A^{\beta_{k-6,k-5}~\beta_{k-7,k-6}
~... ~\beta_{23}}= \pmatrix{1&           &       & && & 0 &   &   \cr
           &          0&       & && & 1 &   &   \cr 
           &          1&0      & && & * &   &     \cr
           &           & 1     &0&& & * &  &     \cr
           &           &       & \ddots&\ddots &&\vdots &      &       \cr 
           &           &        &     &  \ddots    &0& *
                                                         & & \cr 
           &          &       &      &             &  1   &*& &    \cr
           &          &&&&& &1& \cr 
           &          &&&&& &&1 \cr
           &          &&&&& && &1 \cr 
           &          &&&&& && &  &1 \cr
           &          &&&&& && &  & &1 \cr
}$$}
 and the entries $*$ are  those of the $A_{i,i+1}$'s (namely,
$-s_{1,k-6}$, ..., $-s_{k-7,k-6}$ from the top to the bottom). 
 
For the third sequence  $(\beta_{k-6,k-5}~...~\beta_{23})$, they are 
$$  -s_{i-2,k-8}~\hbox{ for } A^{\beta_{i,i+1}}
$$
and: {\small 
$$         m_3:=  A^{\beta_{k-7,k-6}~...~ \beta_{34}}=
   \pmatrix{1&  &           &       & && & 0   \cr 
                 & 1&           &       & && & 0    \cr
           &&          0&       & && & 1 &   &   \cr 
           &&          1&0      & && & * &   &     \cr
           &&           & 1     &0&& & * &  &     \cr
           &&           &       & \ddots&\ddots &&\vdots &      &       \cr 
           &&           &        &     &  \ddots    &0& *
                                                         & & \cr 
           &&          &       &      &             &  1   &*& &    \cr
           &&          &&&&& &1& \cr 
           &&          &&&&& &&1 \cr
           &&          &&&&& && &1 \cr 
           &&          &&&&& && &  &1 \cr
           &&          &&&&& && &  & &1 \cr
           &&          &&&&& && &  & & &1 \cr
}    $$
}

 And so on. We reach the last but one ``sequence'', namely
 $\beta_{{k\over 2}-2,{k\over 2}-1}$ for $k$ even, or $(\beta_{{k-3\over
 2},{k-1\over 2}}\beta_{{k-5 \over 2},{k-3 \over 2}})$  for $k$
 odd. The entries are $-s_{12}$, or $-s_{23}$,$-s_{13}$ respectively. 
 Then 
{\small 
$$
m_{{k \over 2} -2}:= A^{ \beta_{   {k \over 2}-2,{k\over 2}-1   } }=
\pmatrix{ 1&  \cr
           &\ddots & \cr
           &       &1 \cr
           &       & & 0&1\cr
           &       & & 1 & -s_{12}    \cr 
           &       & & & & & 1&   \cr
           &       & & & & &  & \ddots  \cr
           &       & & & & &  & & 1 \cr
}
$$}
or
{\small
$$
 m_{k-5 \over 2}:=
A^{\beta_{{k-3\over
 2},{k-1\over 2}}\beta_{{k-5 \over 2},{k-3 \over 2}}}=
\pmatrix{ 1&  \cr
           &\ddots & \cr
           &       &1 \cr
           &       & & 0&0&1\cr
           &       & & 1 &0 & -s_{13}    \cr 
           &       & &   &1 & -s_{23}     \cr
           &       & &   &  &  &1  \cr
           &       & &   &  & &  & \ddots  \cr
           &       & &   &  & &  &    & 1 \cr
}
$$
}
 The entries for the last braid are more complicated, because the
 entries on the first upper sub-diagonal of the Stokes' matrix have
 been shuffled 	by the   preceeding braids. We  give the result
 ($A_{i,i+1}$  stends for $A^{\beta_{i,i+1}}$ 
{\small
$$
   A_{k-3,k-2} :~~~ -s_{k-3,k-2}$$
$$
   A_{k-4,k-3} :~~~ -s_{k-5,k-2}+s_{k-5,k-4}s_{k-4,k-2}$$
$$
   A_{k-5,k-4} :~~~ -s_{k-7,k-2}+(s_{k-7,k-4}s_{k-4,k-2}+
   s_{k-7,k-6}s_{k-6,k-2})-s_{k-7,k-6}s_{k-6,k-4}s_{k-4,k-2}$$
$$ \vdots$$
$$
   A_{{k \over 2}-1,{k \over 2}} :~~~
   -s_{1,k-2}+(s_{12}s_{2,k-2}+s_{14}s_{4,k-2} +...+
  s_{1,k-4} s_{k-4,k-2})~+~...~+~(-1)^{{k\over2}-1}
   s_{12}s_{24}s_{48}...s_{k-4,k-2}$$
$$
 A_{{k \over 2}-2,{k \over 2}-1} :~~~
 -s_{2,k-2}+(s_{24}s_{4,k-2}+s_{26}s_{4,k-2} +...+
   s_{2,k-4}s_{k-4,k-2})~+~...~+~(-1)^{{k\over2}-2}
   s_{24}s_{48}...s_{k-4,k-2}$$
$$ 
     \vdots $$
$$ 
     A_{34} :~~~
-s_{k-8,k-2}+(s_{k-8,k-4}s_{k-4,k-2}+s_{k-8,k-6} s_{k-6,k-2}) -
s_{k-8,k-6}s_{k-6,k-4} s_{k-4,k-2}
$$
$$ A_{23} :~~~
-s_{k-6,k-2}+s_{k-6,k-4}s_{k-4,k-2}
   $$
$$
  A_{12} :~~~
-s_{k-4,k-2}
$$
}
For $k$ odd
{\small
$$
   A_{k-3,k-2} :~~~ -s_{k-3,k-2}$$
$$
   A_{k-4,k-3} :~~~ -s_{k-5,k-2}+s_{k-5,k-4}s_{k-4,k-2}$$
$$
   A_{k-5,k-4} :~~~ -s_{k-7,k-2}+(s_{k-7,k-4}s_{k-4,k-2}+
   s_{k-7,k-6}s_{k-6,k-2})-s_{k-7,k-6}s_{k-6,k-4}s_{k-4,k-2}$$
$$ \vdots$$
$$
   A_{{k -1\over 2},{k+1 \over 2}} :~~~
   -s_{2,k-2}+(s_{23}s_{3,k-2}+s_{25}s_{5,k-2} +...+
  s_{2,k-4} s_{k-4,k-2})~+~...~+~(-1)^{{k-3\over2}}
   s_{23}s_{35}s_{57}...s_{k-4,k-2}$$
$$
 A_{{k-3 \over 2},{k-1 \over 2}} :~~~
 -s_{1,k-2}+(s_{13}s_{3,k-2}+s_{15}s_{5,k-2} +...+
  s_{1,k-4} s_{k-4,k-2})~+~...~+~(-1)^{{k-3\over2}}
   s_{13}s_{35}s_{57}...s_{k-4,k-2}$$
$$ 
     \vdots $$
$$ 
     A_{34} :~~~
-s_{k-8,k-2}+(s_{k-8,k-4}s_{k-4,k-2}+s_{k-8,k-6} s_{k-6,k-2}) -
s_{k-8,k-6}s_{k-6,k-4} s_{k-4,k-2}
$$
$$ A_{23} :~~~
-s_{k-6,k-2}+s_{k-6,k-4}s_{k-4,k-2}
   $$
$$
  A_{12} :~~~
-s_{k-4,k-2}
$$
}
and: 
$$ 
  A^{\beta_{k-3,k-2}~\beta_{k-4,k-3}~...~\beta_{12}}   =    \pmatrix{
                     0&       & && & 1 &   &   \cr 
                     1&0      & && & * &   &     \cr
                      & 1     &0&& & * &  &     \cr
                      &       & \ddots&\ddots       &      &\vdots   \cr 
                      &       &       &  \ddots     &  0   &*     &      \cr 
                      &       &       &             &  1   &*     &     &    \cr
                      &&&&& &1&0 \cr 
                      &&&&& &0&1 \cr
}$$
which we call $m_{{k\over 2}-1}$ for $k$ even, $m_{k-3\over 2}$ for
                     $k$ odd.  

Then, for  $k$ even 
$$A^{\beta}=m_{{k\over 2}-1}m_{{k\over 2}-2}~...~m_3m_2m_1$$
and for $k$ odd 
$$A^{\beta}=m_{{k-3\over 2}}m_{{k-5\over 2}-2}~...~m_3m_2m_1$$

Doing a careful computation, we obtain (\ref{As1}) and (\ref{As2}). 

\vskip 0.3 cm
2) The second step consists of expressing $A^{\beta}$ in terms of the
entries of the Stokes' factors of $S$, which are simply binomial
coefficients. First we prove the following

\vskip 0.2 cm 
\noindent
{\bf Lemma 9:} {\it  Given an upper triangular $k \times k$ matrix $S$, with entries
$s_{ii}=1$, we can uniquely determine numbers $a_{ij}$ such that, for
$k$, $i$, $j$ all even or all odd: {\small
$$
   s_{ij}=
a_{ij}~+~(a_{i,i+2}a_{i+2,j}+a_{i,i+4}a_{i+4,j}+~...~+a_{i,j-2}a_{j-2,j})~+$$
$$
+~(a_{i,i+2}a_{i+2,,i+4}a_{i+4,j}+a_{i,i+2}a_{i+2,i+6}a_{i+6,j}+~...~+a_{i,j-4}a_{j-4,j-2}
a_{j-2,j})~+.... ~+~ a_{i,i+2}a_{i+2,i+4}a_{i+4,i+6}~...~a_{j-2,j}
$$
}

If $k$ is even, but $i$ is odd, just replace in the formula $i+2$ with $i+1$, $i+4$
with $i+3$, ect. If $k$ is even, but $j$ is odd, just replace $j-2$ with $j-1$, $j-4$
with $j-3$, ect. 
 
If $k$ is odd, but $i$ is even, or $j$ is even, just do the same
replacements as above. More explicitely:
{\tiny
 $$   \hbox{\small $S$}:=\pmatrix{ \cr
                 1& a_{12} &  a_{13}+a_{12}a_{23} & 
                   a_{14}+a_{12}a_{24} &
                   \matrix{a_{15}+(a_{12}a_{25}+a_{14}a_{45})\cr
                  \cr +
                   a_{12}a_{24}a_{45}\cr} &
                   \matrix{a_{16}+(a_{12}a_{26}+a_{14}a_{46})\cr\cr 
                         +
                   a_{12}a_{24}a_{46}\cr} & \matrix{  a_{17}+
                              (a_{12}a_{27}+a_{14}a_{47}
                   +a_{16}a_{67}) \cr\cr
                   +(a_{12}a_{24}a_{47}+a_{12}a_{26}a_{67}+a_{14}a_{46}a_{67}) \cr\cr
                   + a_{12}a_{24}a_{46}a_{67} \cr} & ....\cr
                                                            \cr \cr
                   \cr 
                   & 1 & a_{23} & a_{24} & a_{25} + a_{24}a_{45} &
                   a_{26}+ a_{24}a_{46} & a_{27}+(a_{24}a_{47}+
                       a_{26}a_{67})+a_{24}a_{46}a_{67}& .... \cr
                    \cr \cr 
                   &&1&a_{34}&     a_{35} + a_{34}a_{45} &
                   a_{36}+ a_{34}a_{46} & a_{37}+(a_{34}a_{47}+
                       a_{36}a_{67})+a_{34}a_{46}a_{67}& .... \cr
                    \cr \cr 
&&&1& a_{45} & a_{46} & a_{47} a_{46}a_{67} & .... \cr \cr \cr
&&&&1&a_{56} & a_{57}+a_{56}a_{57}& .... \cr \cr\cr
&&&&&1& a_{67}&....\cr\cr\cr
&&&&&&1&.... \cr \cr\cr
&&&&&&&.... \cr
}$$
}
for $k$ even. 

{\tiny
 $$   \hbox{\small $S$}:=\pmatrix{ \cr
                   1& a_{12} &  a_{13} & 
                   a_{14}+a_{13}a_{34} &
                   a_{15}+a_{13}a_{35} &
                   a_{16}+(a_{13}a_{36}+a_{15}a_{56}) 
                         +
                   a_{13}a_{35}a_{56} & a_{17}+
                              (a_{13}a_{37}+a_{15}a_{57})
                      +a_{13}a_{35}a_{57}& ....\cr
                                                             \cr
                   \cr 
                   & 1 & a_{23} & a_{24}+a_{23}a_{34} & a_{25} + a_{23}a_{35} &
                   a_{26}+
                   (a_{23}a_{36}+a_{25}a_{56})+a_{23}a_{35}a_{56} &
                   a_{27} + (a_{23}a_{37}+a_{25}a_{57})+a_{23}a_{35}a_{57}& .... \cr
                    \cr \cr 
                   &&1&a_{34}&     a_{35}  &
                   a_{36}+ a_{35}a_{56} & a_{37}+a_{35}a_{57}
                             & .... \cr
                    \cr \cr 
&&&1& a_{45} & a_{46}+a_{45}a_{56} & a_{47}+ a_{45}a_{57} & .... \cr \cr \cr
&&&&1&a_{56} & a_{57}& .... \cr \cr\cr
&&&&&1& a_{67}&....\cr\cr\cr
&&&&&&1&.... \cr \cr\cr
&&&&&&&.... \cr
}$$
}
for $k$ odd. 
}

\vskip 0.2 cm 
\noindent
{\bf Proof: } 
We have to solve a non linear system $F_{ij}(a)=s_{ij}$. 

The sum of the differences between the indices of the
factors $a$ in $s_{ij}$ is equal to the difference of the indices of
$s_{ij}$, namely $j-i$. 

>From this it follows that the trems non-linear in the $a$'s occur on
sub-diagonals which lie above all the sub-diagonals containing the
factors of the non linear terms. 

Then, the system $F_{ij}(a)=s_{ij}$ is uniquely solvable, starting
from the first sub-diagonal and successively determining all the
$a_{rs}$ going up diagonal by diagonal. 

\rightline{$\Box$}

\vskip 0.2 cm
\noindent
{\bf Corollary 3:} {\it With the above factorization, the matrix
$A^{\beta}$ becomes: 
 {\small
\be  A^{\beta}(S)= \pmatrix{ .&.&.&.&.&.&. &.&... \cr
                           .&.&.&.&.&.&. &0&... \cr                      
                          .&.&.&.&.&0&0& 1&... \cr
                         .&.&.&0&0&1&0&-a_{68}&...\cr        
                          .&0&0&1&0&-a_{46}&0&-a_{48} & ...\cr
0&1&0&-a_{24}&0&-a_{26}&0&-a_{28}
&... \cr 
1&-a_{12}&0&-a_{14}&0&-a_{16}&0&-a_{18}&... \cr
0&0&1&-a_{34}&0&-a_{36}&0&-a_{38}&... \cr
. &.&0&0&1&-a_{56}&0&-a_{58}&... \cr
 .&.&.&.&0&0&1&-a_{78}&... \cr
                         .&.&.&.&.&.&0&0&... \cr
                           .&.&.&.&.&.&.&.&... \cr
}
\label{As3}
\ee}
for $k$ even; and 
{\small 
\be
  A^{\beta}(S)= \pmatrix{ .&.&.&.&.&.&. &... \cr
                          .&.&.&.&.&.&0& ... \cr
                         .&.&.&.&0&0&1 &...\cr        
                          .&.&0&0&1&0&-a_{57}& ...\cr
0&0&1&0&-a_{35}&0&-a_{37}&... \cr 
1&0&-a_{13}&0&-a_{15}&0&-a_{17}&... \cr
0&1&-a_{23}&0&-a_{25}&0&-a_{27}&... \cr
.&0&0&1&-a_{45}&0&-a_{47}&... \cr
.&.&.&0&0&1&-a_{67}&...\cr
 .&.&.&.&.&0&0&... \cr
                         .&.&.&.&.&.&.&... \cr
                           .&.&.&.&.&.&.&... \cr
}
\label{As4}
\ee}
for $k$ odd. In other words, all the non linear terms in formulae
(\ref{As1}) and (\ref{As2}) drop.
}
\vskip 0.2 cm 
\noindent 
{\bf Proof:} Just substitute the factorization of $S$ in (\ref{As1}),
(\ref{As2}).

\rightline{$\Box$}

\vskip 0.2 cm
  In section 5 we computed the Stokes' factors for $S$. If we sum all the
  factors appearing in formula (\ref{facto}) we get a matrix of the form: 
$$M:= 
 K_{1k}+K_{1,k-1}+K_{k,k-1}+...+K_{k3}+K_{k2}= \pmatrix{ k &   &              &&&& 0\cr
            * & k &              &&&0& * \cr
            * & * & k           &&&*&* \cr
            \vdots&\vdots&\vdots&\ddots& \adots&\vdots&\vdots \cr
            * & * & * &\adots &\ddots&\vdots&\vdots \cr
            * & * & 0 &       &      & k   & *\cr
            * & 0 &   &       &      &     & k  \cr
             }$$
 The $*$ are the binomial coefficient appearing in the factors. If we
 know $M$, we can determine all the entries of the single
             Stokes' factors, because if the entry $(i,j)$ is not
             zero for one factor, then it is zero for all the other
              factors. Now, 
we rename  the entries of the factors according to the following rule:
{\small 
    $$
M:=        \pmatrix{ k &   &              &&&&&&&&& 0\cr
            a_{k-2,k}& k &              &&&&&&&&0& a_{k-2,k-1} \cr
            a_{k-4,k} & a_{k-4,k-2} & k &&&&&&&&a_{k-4,k-3}&a_{k-4,k-1} \cr
            \vdots&\vdots&\vdots&  \ddots&&&&&& \adots &\vdots&\vdots& \cr
a_{4k} &a_{4,k-2} & ...& a_{46} &k&
            0&0&0&a_{45}&a_{47}&...&a_{4,k-1} \cr
a_{2,k}&a_{2,k-2}&...& a_{26}&a_{24}&k&0&a_{23}&a_{25}&a_{27}&...&a_{2,k-1}
            \cr
a_{1k} & a_{1,k-2}&
            ...&a_{16}&a_{14}&a_{12}&k&a_{13}&a_{15}&a_{17}&...&a_{1.k-1} \cr
a_{3k}&a_{3,k-2}&...&a_{36}&a_{34}&0&0&k&a_{35}&a_{37}&...&a_{3,k-1}
            \cr
a_{5k}&a_{5,k-2}&...&a_{56}&0&0&0&0&k&a_{57}&...&a_{5,k-1}
            \cr
\vdots&\vdots&\adots& &&&&&&\ddots&...&...\cr
a_{k-3,k}&a_{k-3,k-2}&&&&&&&&&k&a_{k-3,k-1}\cr
a_{k-1,k}& &&&&&&&&&& 1\cr
}$$
}
for $k$ even.
{\small 
$$ 
M= \pmatrix{
k&&&&& & & & &      &
& &0 \cr
a_{k-2,k}&k&&&& & & & &      &
& &a_{k-2,k-1} \cr
 a_{k-4,k}&a_{k-4,k-2}&k&&& & & & &      &
& a_{k-4,k-3}&a_{k-4,k-1} \cr
 \vdots&\vdots&\vdots&\ddots&& & & & &      &
\adots& \vdots&\vdots \cr
a_{5k}&a_{5,k-2}&a_{5,k-4}&...&k&0&0&0&0&a_{56}&
...& a_{5,k-3}&a_{5,k-1} \cr
a_{3k}&a_{3,k-2}&a_{3,k-4}&...&a_{35}&k&0&0&a_{34}&a_{36}&
...& a_{3,k-3}&a_{3,k-1} \cr
a_{1k}&a_{1,k-2}&a_{1,k-4}&...&a_{15}&a_{13}&k&a_{12}&a_{14}&a_{16}&
...& a_{1,k-3}&a_{1,k-1} \cr
a_{2k}&a_{2,k-2}&a_{2,k-4}&...&a_{25}&a_{23}&0&k&a_{24}&a_{26}&...   &
 a_{2,k-3}&a_{2,k-1} \cr
a_{4k}&a_{4,k-2}&a_{4,k-4}&...&a_{45}&0&0&0&k&a_{46}&...   &
 a_{4,k-3}&a_{4,k-1} \cr
\vdots&\vdots&&\adots&& & & & &      &  \ddots    &
 \vdots   &\vdots    \cr
a_{k-3,k}&a_{k-3,k-2}&&&& & & & &      &      &
      k   &a_{ k-3,k-1} \cr
a_{k-1,k}&           &&&& & & & &      &      &
          &    k        \cr
}
$$
}
for k odd. The strange labelling is simply the one such that 
$$P M P^{-1} = \pmatrix{k&a_{12}&a_{13}&a_{14}& ...\cr
                        & k& a_{23}&a_{24}& ... \cr
                        & & k&a_{34}&...\cr
                        & &        & \ddots&  \vdots\cr
                        & &        & &  k \cr}
$$
where the matrix of permutation is 
$$P=\pmatrix{ 
             &&&& & &&1&0&&&\cr
             &&&& &&1&0&0&&&&\cr
             &&&& &0&0&0&1&0&&&\cr
             &&&&&1&0&0&0&0& &&&\cr
             &&&& 0&0&0&0&0&1&0& && \cr
             \vdots&\vdots&
             &\vdots&\vdots&\vdots&\vdots&\vdots&\vdots&\vdots&\vdots& \vdots&&
             \vdots &\vdots \cr
              & 1 &&&&&&&&&&&&0&\cr
             0&0&&&&&&&&&&&&0&1\cr
             1&0& &&&&&&&&&&&0&0\cr
}
$$
for $k$ even (the 1 on the first row is on the ${k\over 2}+1$-th
column)  ; and 
$$P=\pmatrix{ 
             &&&& & &0&1&&&&\cr
             &&&& &&0&0&1&&&&\cr
             &&&& &0&1&0&0&0&&&\cr
             &&&&&0&0&0&0&1& &&&\cr
             &&&&0&1&0&0&0&0&0 && \cr
             \vdots&\vdots&
             &\vdots&\vdots&\vdots&\vdots&\vdots&\vdots&\vdots&\vdots& \vdots&&
             \vdots &\vdots \cr
             0 & 1 &&&&&&&&&&&&0&\cr
             0&0&&&&&&&&&&&&0&1\cr
             1&0& &&&&&&&&&&&0&0\cr
}
$$
for $k$ odd (the 1 on the first row is on the ${k+1\over 2}$-th column).

With this choice of the labelling, the product $S_{upper}:=P \left(
K_{1k}K_{1,k-1}K_{k,k-1}~...~ K_{k3}K_{k2}~\right)P^{-1}$ is precisely
factorized as in lemma 9.

 Then we can write  the entries of $A^{\beta}$
 (formulae (\ref{As3}) (\ref{As4})) from
 the entries of the Stokes factors (which are binomial coefficient). 

The final result is precisely the claim of the proposition.

\rightline{$\Box$} 
\vskip 0.3 cm
We are ready to prove the main result: 
\vskip 0.2 cm 
\noindent
{\bf Theorem 2:} {\it Consider the Stokes matrix $S=T_F^{k \over 2}~T^{-{k
\over 2}}$ ($k$ even) or $S=T_F^{k-1 \over 2} K_{k2} T^{-{k-1 \over
2}}$ ($k$ odd) for the quantum cohomology of ${\bf CP}^{k-1}$ and set
it in the upper triangular form $S_{upper}=P~S~P^{-1}$ by the
permutation $P$. Then, there exists a braid $\beta$ (Lemma 8), represented by a
matrix $A^{\beta}$ (Proposition), which sets $S_{upper}$ in the form
(\ref{SSCAN}).  
The last column is negative, but conjugation by } ${\cal I} =$ 
diag$(1,1,...,1,-1)$ {\it  makes it positive. We reach the canonical 
form:
$$
  s_{ij}=\bin{k}{j-i},~~~~i<j$$
Another conjugation by } diag$(-1,1,-1,1,-1,...)$ {\it 
brings the matrix in the equivalent
canonical form 
$$ s_{ij}=(-1)^{j-i}\bin{k}{j-i}, ~~~~ i<j$$
Finally, by the action of the braid group, the last matrix can be put
in the form 
$$s_{ij}=\bin{k-1+j-i}{j-i},~~~~i<j$$
In all the above matrices 
$$ s_{ii}=1,~~~~s_{ij}=0 ~~~~i>j$$
}
\vskip 0.2 cm
\noindent
{\bf Proof:} 
First, we want to explain which is the braid which brings the upper
triangular matrix
with entries $ s_{ij}=(-1)^{j-i}\bin{k}{j-i}$ in the matrix 
$s_{ij}=\bin{k-1+j-i}{j-i}$. We make use of the following 
known result \cite{Zas}: 

\vskip 0.2 cm 
\noindent
{\it Consider the upper triangular Stokes' matrix S, the  braid 
$\beta=
 \beta_{12}(\beta_{23}\beta_{12})(\beta_{34}\beta_{23}\beta_{12})$~...~
$(\beta_{n-1,n}\beta_{n-2,n-1}...\beta_{23}\beta_{12})$ and the 
permutation $$P=\pmatrix{&&& 1\cr
                        &&1&\cr
                       &\adots&& \cr
                      1&&& \cr}$$ 
Then, the relation 
\be
    \left[S^{-1}\right]^{\beta}=PS^TP
\label{Smeno1}
\ee
holds. 
}
\vskip 0.2 cm

Observe that for the matrix $S$, whose upper triangular part has
entries  
$s_{ij}=\bin{k-1+j-i}{j-i}$, we have   $PS^TP\equiv S$. Moreover,
$S^{-1}$ is upper triangular with entries $
s_{ij}=(-1)^{j-i}\bin{k}{j-i}$. This proves that $S$ and $S^{-1}$ are 
equivalent w.r.t the action of the braid group. 
\vskip 0.2 cm

Let us now prove the theorem staring from $k$ even. We have to prove that
$A^{\beta}~P~T_F^{k \over 2}~T^{-{k\over 2}}~P^{-1} [A^{\beta}]^T$ is in
``canonical form'' (\ref{SSCAN}). The proof ``reduces'' to the computation of
products of matrices explicitely given. We do the products in an
shrewd way. First we rewrite  
$$ S^{\beta} = A^{\beta}~\left( P~T_F~P^{-1}\right)^{k\over 2}~\left(P~ T^{-1}~ P^{-1}
\right)^{k\over 2}~[A^{\beta}]^T$$
and we compute: {\tiny
$$\hbox{\small $P~T_F~P^{-1}$}= \pmatrix{ 
             &&&& & &&1&0&&&\cr
             &&&& &&1&0&0&&&&\cr
             &&&& &0&0&0&1&0&&&\cr
             &&&&&1&0&0&0&0& &&&\cr
             &&&& 0&0&0&0&0&1&0& && \cr
             \vdots&\vdots&
             &\vdots&\vdots&\vdots&\vdots&\vdots&\vdots&\vdots&\vdots& \vdots&&
             \vdots &\vdots \cr
              & 1 &&&&&&&&&&&&0&\cr
             0&0&&&&&&&&&&&&0&1\cr
             1&0& &&&&&&&&&&&0&0\cr
}
\pmatrix{ 0 &     &   &    &...& 1 \cr
              -1&  0  &   &     &&    \cr
                &  -1 & 0 &     &&   \cr
                &     & -1 &    &&        \cr
                &     &    & \ddots&\ddots& \vdots\cr             
                &     &    &      &    -1 &0\cr   
}
\pmatrix{
&&&&&&0&1\cr
&&&&&1&0&0\cr
&&&&\vdots&\cr 
&&0&1 \cr
&1&0&0 \cr
1&0&0&0\cr
0 &0 &1&0\cr
&&0&0\cr
&&&&\vdots\cr
&&&&&0&0&0\cr
&&&&&&1&0 \cr
 }
$$
}
(in $P^{-1}$ the 1 on the first column is on the ${k\over 2}+1$-th
row)
{\small
$$
  = \pmatrix{ 0&-1\cr
              0&0&0& -1 \cr
              -1& 0&0&0& & \cr
               &0&0&0&0&-1 \cr
              &  &-1&0&0&0\cr
              &&&&&&\vdots&\vdots \cr
              &&&&&&&&0 &0&0&-1\cr 
              &&&&&&&&-1&0&0&0\cr
              &&&&&&&&&&1&0 \cr
}
$$}
Thus
$$ \left(P~T_F~P^{-1}\right)^{k\over 2}= \pmatrix{
&&&&&&&-1 \cr
&&&&&&1\cr
&&&&&-1\cr
&&&&\adots\cr
&&&-1 \cr
&&1\cr
&-1\cr
                       1\cr
}
$$
Then, using the expression of $A^{\beta}$ from the proposition: 
{\tiny
$$ 
 \hbox{\small $F:=  A^{\beta}~\left(P~T_F~P^{-1}\right)^{k\over 2}$}= 
\pmatrix{
           0&0&1 \cr
           &0&\bin{k}{1}&0&1 \cr
           &&\bin{k}{2}&0&\bin{k}{1} \cr
           &&\bin{k}{3}&0&\bin{k}{2}&0&1\cr
           &&\vdots    & &\vdots& & \vdots& \cr
           &&&&&&&&...&1\cr
           &&*&&*&&*&&...&\bin{k}{1}&0&1 \cr
           &&*&&*&&*&&...&\bin{k}{2}&0&\bin{k}{1}&-1\cr
           &&*&&*&&*&&...&\bin{k}{3}&-1&0\cr
           &&\vdots&&\vdots&&\vdots \cr
           &&\bin{k}{k-4}&0&\bin{k}{k-5}\cr
           &&\bin{k}{k-3}&-1&0\cr
          \cr           
      0&-1&0\cr
\cr      
     -1&0\cr
}
$$
}
(the -1 on the last column is on the ${k\over 2}$-th row). 

Using the explicit expressions for $K_{k3}$ and  $K_{k2}$ of section 5
we compute
$$ 
  T^{-1}= T_F^{-1}~K_{k3}~K_{k2}=\pmatrix{0&-1\cr
                                          &0&-1\cr
                                          & & \ddots&\ddots \cr
                                          &&&& 0&-1 \cr
                                          &&&&&0&-1 \cr
                                          1&&&&&&0 \cr
}~K_{k3}~K_{k2}=$$
{\tiny
$$
  = \pmatrix{ \bin{k}{1}& -1   && &&&&&&&&& -\bin{k}{2} \cr
  & 0& -1&   &&&&&&&& -\bin{k}{4}&-\bin{k}{3} \cr
  &  &  0&-1      &&&&&&&-\bin{k}{6}&-\bin{k}{5}&0 \cr \cr
  &  &   & \ddots &\ddots &&&&&\adots&\adots \cr
\cr
  &  &   &        &    0  &-1&0& -\bin{k}{k-2}&-\bin{k}{k-3}\cr
  &  &   &        &       & 0& -1&-\bin{k}{k-1} \cr\cr
  &&&&&& 0& -1 \cr\cr
  &&&&&&&   0&-1 \cr
\cr
  &   &   &       &    &&   &   & \ddots&\ddots &\cr\cr\cr\cr\cr
  &   &   &       &      && &   &   &        &&0&-1\cr \cr
  &   &   &      &        &&&    &  &        &&  &0\cr
}
$$
}
Then {\tiny
$$
\hbox{\small $P~T^{-1}~P^{-1}$}= \pmatrix{ 0&0&-1\cr
                           -1&0&-\bin{k}{k-1} \cr
                           0&0&0&0&-1\cr
                           0&-1&-\bin{k}{k-2}&0&-\bin{k}{k-3}&0 \cr
                            &0&0&0&0&0&-1\cr
&&&-1&-\bin{k}{k-4}&0&-\bin{k}{k-5} \cr
&&&  &             & &        & \ddots \cr
&&&  &             & &        &        & 0&0&0&-1& 0\cr
&&&  &             & &        &        & -1&-\bin{k}{4}&0&-\bin{k}{3}& 0\cr
&&&  &             & &        &        &  &0&0&0& 1\cr
&&&  &             & &        &        &  &0&-1&-\bin{k}{2}&
                           \bin{k}{1}\cr
}
$$
}
After this we  computed $F~\left[P~T^{-1}~P^{-1}\right]$,
$F~\left[P~T^{-1}~P^{-1}\right]~ \left[P~T^{-1}~P^{-1}\right]$,
$F~\left[P~T^{-1}~P^{-1}\right]^3 $, ..., 
$ F~\left[P~T^{-1}~P^{-1}\right]^{k\over 2}$. We omit the intermediate
computations and we give the final result: 
{\tiny
$$
\hbox{\small $F_1:= F~\left(P~T^{-1}~P^{-1}\right)^{k\over 2}$}= 
\pmatrix{
&&&&&&&&    &&1&\bin{k}{2}&-\bin{k}{1} \cr
&&&&&&&&    &\bin{k}{4}&0&\bin{k}{3}&-\bin{k}{2} \cr
&&&&&&&    &&\bin{k}{5}&0&\bin{k}{4}&-\bin{k}{3} \cr
&&&&&&&... &&\vdots&\vdots &\vdots      & \vdots        \cr
&&&&&1&\bin{k}{k-6}& ...&...&*&0&*&*\cr
&&&1&\bin{k}{k-4}&0&\bin{k}{k-5}&...&...&*&0 &*&* \cr 
&1&\bin{k}{k-2}&0&\bin{k}{k-3}&0&\bin{k}{k-4}&...&...&*&0 & *&*\cr
1&0&\bin{k}{k-1}&0&\bin{k}{k-2}&0&\bin{k}{k-3}&...&...&*&0&*&* \cr
&&1&0&\bin{k}{k-1}&0&\bin{k}{k-2}&...&...&*&0&*&*\cr
&&&&1&0&\bin{k}{k-1}&...&...&*&0&*&*\cr
&&&&&&1&...&...&*&0&*&*\cr
&&&&&&&...&...&\vdots&\vdots&\vdots&\vdots \cr
&&&&&&&   &&\bin{k}{k-1}&0&\bin{k}{k-2}&-\bin{k}{k-4} \cr
&&&&&&&   &&&1&\bin{k}{k-1}&-\bin{k}{k-2} \cr
&&&&&&&   &&& &&-\bin{k}{k-1} \cr
&&&&&&&   && &&       & 1 \cr
}
$$
}
 Now, multiplying $F_1~\left[A^{\beta}\right]^T$, we obtain precisely
 the ``canonical form'' (\ref{SSCAN}).
\vskip 0.2 cm 
For $k$ odd, we did a similar computation. We omit the detail, but we
indicate the order of multiplications which yielded the most simple
expressions to multiply step by step. Our aim is to compute  $A^{\beta}~P~T_F^{k-1\over
2}~K_{k2}~(T_F^{-1}~K_{k3}~K_{k2})^{k-1\over 2}~P^{-1}~[A^{\beta}]^T$.
First, we computed $PT_FP^{-1}$, then  $(PT_FP^{-1})^{k-1\over 2}$,
then $F:=A^{\beta}~(PT_FP^{-1})^{k-1\over 2}$. After this, we computed
$PK_{k2}P^{-1}$, and $F_1:=F~PK_{k2}P^{-1}$. Finally, we calculated
$m:=P~T_F^{-1}~K_{k3}~K_{k2}~ P^{-1}$ and $F_1~m$, $F_1~m~m$, ...,
$F_2:=F_1~m^{k-1\over 2}$. The matrix $F_2~[A^{\beta}]^T$ proved
to be  in
``canonical form''. 

\rightline{$\Box$}


\section{ Canonical form of $S^{-1}$}

  The matrix $S^{-1}$ such that $Y_R(z)=Y_{L}(z)~S^{-1}$ can be put in
  the same canonical form of $S$, as a consequence of the relation 
(\ref{Smeno1}). The only remarks we want to add concern the braid which
  brings $S^{-1}$ to the canonical form, because it arranges the lines
  $L_j$ in a ``beautiful'' shape. 

\vskip 0.2 cm 
\noindent
{\bf Lemma $8^{\prime}$:} {\it Let the points $u_j$ be in lexicographical order
w.r.t the admissible line $l$. Let us denote
$\sigma_{i,i+1}:=\beta_{i,i+1}^{-1}$. Then, the following  braid
arranges the poins in cyclic clockwise order, $u_1$ being the first
point  in $\Pi_L$ for $k$ even, or the last in $\Pi_R$ for $k$ odd 
(w.r.t the clockwise order) (see figure 6):
$$ 
  \beta^{\prime}:= \left[~( \sigma_{34}\sigma_{56}\sigma_{78}
  ~...~\sigma_{k-3,k-2} \sigma_{k-1,k})~( \sigma_{45}\sigma_{67}\sigma_{89}
  ~...~\sigma_{k-4,k-3} \sigma_{k-2,k-1})~...~(\sigma_{{k\over 2},{k\over
  2}+1}\sigma_{{k\over 2} +2,{k\over 2}+3})~\sigma_{{k\over 2}+1,{k
  \over 2}+2} \right]$$
$$\left[~( \sigma_{{k\over 2}+2,{k\over
  2}+3}\sigma_{{k\over 2}+3,{k\over
  2}+4}~...~\sigma_{k-2,k-1}\sigma_{k-1,k})~( \sigma_{{k\over 2}+2,{k\over
  2}+3}\sigma_{{k\over 2}+3,{k\over
  2}+4}~...~\sigma_{k-2,k-1})~....~(\sigma_{{k\over 2}+ 2,{k\over
  2}+3}\sigma_{{k\over 2}+3,{k\over
  2}+4})~\sigma_{{k\over 2}+ 2,{k\over
  2}+3} \right] 
$$
for $k$ even, and 
$$ 
\beta^{\prime}:=\beta_{12}~ \left[
~(\sigma_{45}\sigma_{67}\sigma_{89}~...~\sigma_{k-3,k-2}
\sigma_{k-1,k})~(\sigma_{56}
\sigma_{78}~...~\sigma_{k-4,k-3}\sigma_{k-2,k-1})~....~(\sigma_{{k+1\over
2},{k+3 \over 2}}\sigma_{{k+5 \over 2},{k+7\over 2}})~\sigma_{{k+3\over
2},{k+5\over 2}}\right]$$
$$
  \left[~(\sigma_{{k+5\over 2},{k+7 \over 2}} \sigma_{{k+7 \over 2},{k+9
\over 2}}~...~\sigma_{k-2,k-1}\sigma_{k-1,k})~(\sigma_{{k+5\over 2},{k+7 \over 2}} \sigma_{{k+7 \over 2},{k+9
\over 2}}~...~\sigma_{k-2,k-1})~....~(\sigma_{{k+5\over 2},{k+7 \over 2}} \sigma_{{k+7 \over 2},{k+9
\over 2}})~\sigma_{{k+5\over 2},{k+7 \over 2}}
\right]
$$
for $k$ odd.
}

\begin{figure}
\epsfxsize=15cm
\epsffile{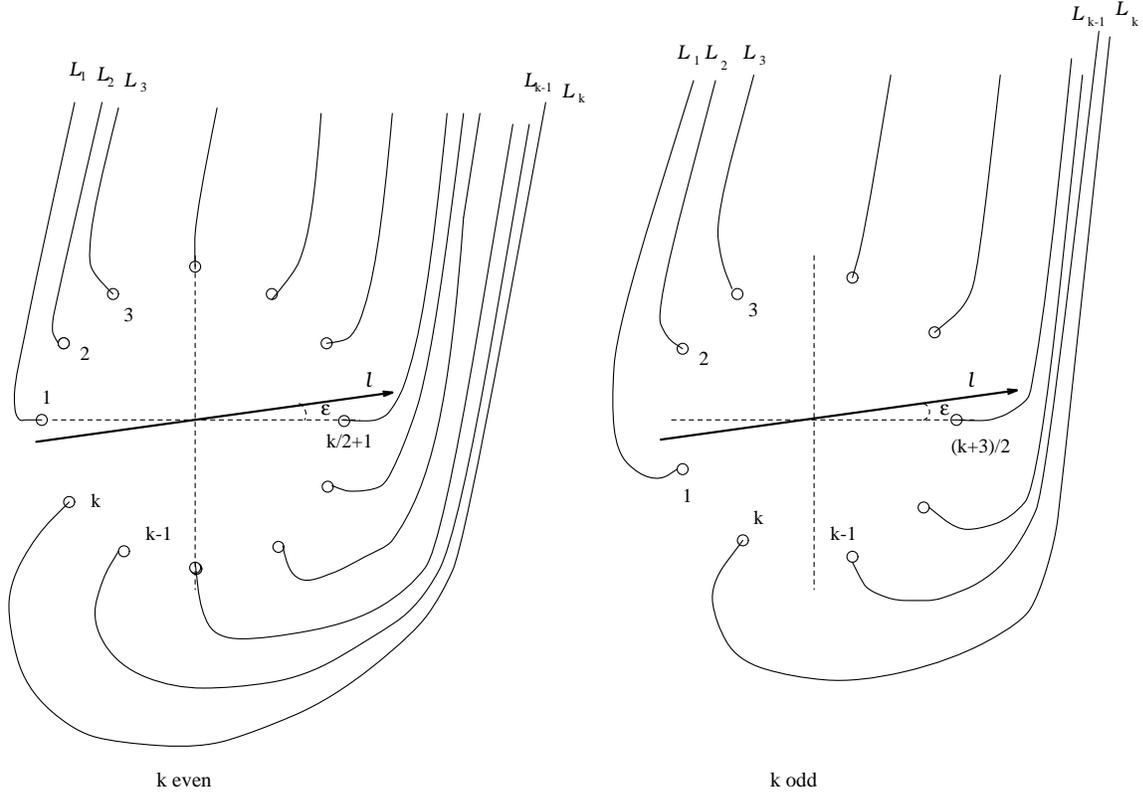}
\caption{ Lines $L_j$ (Branch cuts) after the braid which brings $S^{-1}$ to
canonoical form}
\end{figure}

\vskip 0.2 cm 
 A carefull consideration of the topological effect of  the braid on the lines
 $L_j$ shows that they are arranged as in figure 6.  To reconstruct
 the configuration it is enough to know the admissible line $l$ (at
 angle $\epsilon$ w.r.t. the positive real axis). In fact,  $u_1$ 
is the first point in $\Pi_L$ (in clockwise order) for $k$ even, or
 the last in $\Pi_R$ for $k$ odd. The lines come out of the points in
 centrifugal directions. They go to infinity, without intersections
 (so preserving their lexicographical order w.r.t $l$)  
with the original asymptotic
 direction ${\pi \over 2}-\epsilon$.  
\vskip 0.2 cm
\noindent
{\bf Proposition $1^{\prime}$} {\it The matrix representing
$\beta^{\prime}$ is 
{\tiny
$$ 
A^{\beta^{\prime}}=\pmatrix{ 1\cr
                              &1\cr
 &&-\bin{k}{k-3}&1\cr
 &&-\bin{k}{k-4}&0&-\bin{k}{k-5}&1\cr
&& -\bin{k}{k-5}&0&-\bin{k}{k-6}&  \cr
 &&  \vdots         & & \vdots           &  &...  \cr
 &&  -*          & & -*           & & ...       &1 \cr
 &&  -*          & & -*             & &  ...      &&-\bin{k}{3}&1\cr
 &&  -\bin{k}{{k\over2}-1}   & &  -\bin{k}{{k\over2}-2} & & ...       &&
                              -\bin{k}{2}&0&-\bin{k}{1}&1 \cr
 &&  ~~\bin{k}{{k\over 2}-2}& & ~~\bin{k}{{k\over
                              2}-3}&&...&&~~\bin{k}{1}&0&~~1 \cr
 &&   ~~*& & ~~* & & ...&&~~1  \cr
 &&  ~~\vdots&&~~\vdots&& ... \cr
 &&~~\bin{k}{2} & &~~\bin{k}{1} \cr
 &&~~\bin{k}{1} &0&~~1 \cr
&& ~~1 &0 \cr
}
$$}
for $k$ even, and
 {\tiny
$$ 
A^{\beta^{\prime}}=\pmatrix{ 0&1\cr
                             1&\bin{k}{1}\cr
&&1\cr 
&&&-\bin{k}{k-4}&1\cr
 &&&-\bin{k}{k-5}&0&-\bin{k}{k-6}&1\cr
&&& -\bin{k}{k-6}&0&-\bin{k}{k-7}&  \cr
 &&&  \vdots         & & \vdots           &  &...  \cr
 &&&  -*          & & -*           & & ...       &1 \cr
 &&&  -*          & & -*             & &  ...      &&-\bin{k}{3}&1\cr
 &&&  -\bin{k}{{k-3\over2}}   & &  -\bin{k}{{k-5\over2}} & & ...       &&
                              -\bin{k}{2}&0&-\bin{k}{1}&1 \cr
 &&&  ~~\bin{k}{{k-5\over 2}}& & ~~\bin{k}{{k-7\over
                              2}}&&...&&~~\bin{k}{1}&0&~~1 \cr
 &&&   ~~*& & ~~* & & ...&&~~1  \cr
 &&&  ~~\vdots&&~~\vdots&& ... \cr
 &&&~~\bin{k}{2} & &~~\bin{k}{1} \cr
 &&&~~\bin{k}{1} &0&~~1 \cr
&&& ~~1 &0 \cr
}
$$}
for $k$ odd. 
}

The proposition is proved as in the previous section. Finally, 
the analogous of Theorem 2 holds: 
\vskip 0.2 cm
\noindent
{\bf Theorem $2^{\prime}$} {\it Consider the Stokes matrix $S^{-1}=T^{k \over 2}~T_F^{-{k
\over 2}}$ ($k$ even) or $S=T^{k-1 \over 2} K_{k2}^{-1} T_F^{k-1 \over
2}$ ($k$ odd) for the quantum cohomology of ${\bf CP}^{k-1}$ and set
it in the upper triangular form $S_{upper}^{-1}=P~S^{-1}~P^{-1}$ by the
permutation $P$. Then, there exists a braid $\beta^{\prime}$ (Lemma
$8^{\prime}$),  represented by a
matrix $A^{\beta}$ (Proposition $1^{\prime}$), which sets $S_{upper}^{-1}$ in the
``canonical form''. The entries $s_{ij}$ have minus sign
for $j>{k\over 2}+1,~~i\leq {k\over 2}+1$ and $k$ even, or
$j>{k+3\over 2},~~i\leq {k+3 \over 2}$ and $k$ odd. A suitable
conjugation by ${\cal I}=$ diag$( 1,1,...1,-1,...,-1)$ sets all signs
positive and the final matrix has entries 
  $$ s_{ij}=\bin{k}{j-i}, ~~~s_{ji}=0,~~~\hbox{ for }i<j,~~~~s_{ii}=1
$$ 
}

\vskip 0.2 cm
\noindent
{\bf Proof: } The proof is similar to the one of theorem 2.
\rightline{$\Box$}

\vskip 0.2 cm 
Examples are found in Appendix 2. 


\vskip 0.3 cm 
\section{ Relation between Irregular and Fuchsian systems} 
\vskip 0.2 cm

  Let  us consider the fuchsian system 
\be
 \left( U-\lambda \right) ~{d \phi \over d \lambda} = \left( {1 \over
 2} +V \right) ~\phi 
\label{Fuc20}
\ee
which can also be written 
$$ 
  {d \phi \over d \lambda}= \sum_{j=1}^k~{A_j \over \lambda - u_j}
  ~\phi 
$$
$$ 
   A_j = - E_j ~ \left( {1 \over
 2} +V \right),~~~~~~~(E_j)_{jj}=1,~~\hbox{ otherwise } (E_j)_{nm}=0
$$
Around the point $u_j$ a fundamental matrix has the form 
$$ 
   \left[B_0 + O\left(\lambda - u_j \right)  \right]~(\lambda - u_j)^M
$$ 
where $M=$diag$(-{1\over 2},0,...,0)$ and the columns of $B_0$ are 
the eigenvectors of $A_j$; in particular,  the first column is
$(0,..,0,1,0,...,0)^T$, and $1$ occurs at the $j^{th}$ position. Then,
the system has  $k$ independent vector solutions, of which $k-1$ are
regular near $u_j$ and the last  is 
$$ 
   \phi^{(j)}(\lambda)={ 1 \over \sqrt{ \lambda - u_j}}~\pmatrix{0\cr 
                                \vdots \cr
                                1 \cr
                                \vdots \cr
                                0 \cr
                              } + O \left( \sqrt{ \lambda - u_j}
                                \right)~~~\lambda\to u_j
$$
where $1$ occurs at the $j^{th}$ row. For any $u_j$ we can construct
such a basis of solutions. The branch of $ \sqrt{ \lambda -
u_j}$ is chosen as follows: let us consider an angle $\eta$ with a range of
$2 \pi$, for example $-{\pi \over 2} \leq \eta < {3 \pi
\over 2}$, such that $ \eta \neq \arg(u_i-u_j)$, $\forall i
\neq j$. Then consider the cuts $L_j=\{ \lambda= u_j + \rho e^{i
\eta},~~\rho>0\}$. Actually, the cuts have two sides,  $L_j^{+}=
\{ \lambda= u_j + \rho e^{i
\eta},~~\rho>0\}$ and $L_j^{-}=\{ \lambda= u_j + \rho e^{i(
\eta- 2 \pi)},~~\rho>0\}$. The branch is determined by the
choice $ \log (\lambda -u_j)= \log |\lambda - u_j| + i \eta $ on
$L_j^{+}$ and  $ \log (\lambda -u_j)= \log |\lambda - u_j| + i (\eta- 2\pi) $ on
$L_j^{-}$. On ${\bf C} \backslash \bigcup_j L_j$, $\sqrt{ \lambda
-u_1}$, ..., $\sqrt{ \lambda -u_k}$ are single valued.

For any two (column) vector solutions $\phi(\lambda)$, $\psi(\lambda)$
we define the symmetric bilinear form: 
$$ 
   \left( \phi,\psi \right) := \phi(\lambda)^T (\lambda -U)
   \psi(\lambda) 
$$ 
which is independent of $\lambda$ and $u_1$, ...,$u_k$. 
Let $G$ be the matrix whose entries are $G_{ij}=\left(
\phi^{(i)},\phi^{(j)} \right)$. In particular, $G_{ii}=1$. 
 Then, it can
be proved (see
\cite{BJL2} and also \cite{Dub2} ) that near $u_j$ 
$$ 
   \phi^{(i)}(\lambda) = G_{ij}\phi^{(j)}(\lambda)+ r_{ij}(\lambda)
$$ 
where $  r_{ij}(\lambda)$ is regular near $u_j$. For a
counter-clockwise  loop around $u_j$ the monodromy of $\phi^{(i)}$ is 
 $$ 
  \phi^{(i)} \mapsto R_j \phi^{(i)} := \phi^{(i)} - 2{ \left(
\phi^{(i)},\phi^{(j)} \right) \over \left(
\phi^{(j)},\phi^{(j)} \right)} ~\phi^{(j)} \equiv \phi^{(i)} - 2
G_{ij} ~\phi^{(j)}
$$

 Then, the monodromy group of (\ref{Fuc20}) acts on $\phi^{(1)}$, ...,
 $\phi^{(k)}$ as a reflection group whose Gram matrix is $2G$. In
 particular, $\phi^{(1)}$, ..., $\phi^{(k)}$ are linearly independent
 (and then a basis) if and only if $\det G \neq 0$.

\vskip 0.2 cm
 
 Now consider an oriented line $l$ of argument  $\theta = {\pi \over 2}
 - \eta$, and for any $j$ define the following vector
\be 
   \tilde{Y}^{(j)} = -{\sqrt{z} \over 2 \sqrt{\pi}} \int_{\gamma_j} d
   \lambda ~\phi^{(j)}(\lambda)~e^{\lambda z}
\label{Laplace}
\ee
which is a Laplace transform of  $\phi^{(j)}$ . The path $\gamma_j$
comes from infinity near $L_j^{+}$, encircles $u_j$ and returns to
infinity along $L_j^{-}$. We can define $\Pi_L=\{\theta< \arg z <\theta
+\pi\}$ 
and $\Pi_R=\{ \theta - \pi< \arg z < \theta \}$. 
$\lambda=\infty$ is a regular singularity for
(\ref{Fuc20}), then the integrals exist for $z \in \Pi_L$, and the 
non-singular matrix $\tilde{Y}(z):= \left[\tilde{Y}^{(1)}|
...|\tilde{Y}^{(k)} \right]$ has the asymptotic behaviour 
$$ \tilde{Y}(z)\sim \left(I + O\left( {1 \over
z}\right)\right)~e^{zU}~~~z \to \infty,~~z \in \Pi_L$$
and satisfies the system (\ref{20}). Then it is a fundamental matrix
$\tilde{Y}_L$. Note that $l$ is admissible, since it does not contain
Stokes` rays.

\vskip 0.2 cm 
 It is a fundamental result \cite{Dub2} that the Stokes' matrix of
(\ref{20}) satisfies 
$$ 
     S+S^T = 2G
$$ 
\vskip 0.3 cm


\section{ Monodromy Group of the Quantum Cohomology of ${\bf CP}^{k-1}$}
\vskip 0.2 cm

A system like (\ref{Fuc20})  comes about in the theory of Frobenius
manifolds (replace $U\mapsto U(t)$, $V\mapsto V(t)$). It  determines
flat coordinates $x^1(t,\lambda)$, ..., $x^k(t,\lambda)$ 
for a linear pencil
of metrics $(~,~)-\lambda<~,~>^*$ ($(~,~)$ is the {\it intersection
form}  \cite{Dub1} \cite{Dub2}).  We write
a gauge equivalent form (gauge $X(t^2)$) at the semisimple point
$(0,t^2,0,...,0)$ 
\be
(\hat{\cal U}(t^2)-\lambda)~{d\psi \over d \lambda}= \left({1\over 2} +
\hat{\mu}\right)\psi
\label{sticaz}
\ee
A fundamental matrix $\psi(t,\lambda)$ has entries 
$\psi_{\alpha}^{(j)}(t,\lambda)=\partial_{\alpha}x^{(j)}(t,\lambda)$. 
Moreover, by (\ref{Laplace})
\be
\partial_{\alpha} \tilde{t}^j(t,z)= 
-{\sqrt{z}\over 2 \sqrt{\pi}}
\int d \lambda~\partial_{\alpha}x^j(t,\lambda)~e^{\lambda z}
\label{Basta!!}
\ee

The {\it Monodromy group} of the Frobenius manifold $M$ is the group of
the transformations which
$(x^{1}(t,\lambda),...,x^{k}(t,\lambda))$ 
undergo
when $t$ moves in $M\backslash \Sigma_{\lambda}$, where $\Sigma_{\lambda}= \{ 
t \in M~|$ $\det\left[  (~,~)-\lambda <~,~>\right]=0$ $\}$ is the {\it
discriminant} of the linear pencil. 

Due to formula (\ref{Basta!!}), for ${\bf CP}^{k-1}$ this group 
is generated by the monodromy of the
solutions of (\ref{sticaz}) when $\lambda$ describes loops around 
$u_1(t),...,u_k(t)$ (see \cite{Dub1} \cite{Dub2}). 
To these loops, we must add the effect of the
displacement $t^2\mapsto t^2+2\pi i$. In fact, in this case 
$$[\varphi^{(1)}(ze^{t^2\over k}),..., \varphi^{(k)}(ze^{t^2\over k})]
\mapsto 
[\varphi^{(1)}(ze^{t^2\over k}),..., \varphi^{(k)}(ze^{t^2\over k})]
~T$$
and the same holds for $\tilde{t}(z,(0,t^2,...,0))$.

Then, the  monodromy group of the quantum cohomology of 
 ${\bf CP}^{k-1}$ is generated by the
transformations $R_1$, $R_2$, ..., $R_k$, $T$ introduced in the
preceeding sections.

 We are going to study the structure of the monodromy group of ${\bf
CP}^{k-1}$ for any $k\geq 3$.   
Recall that the matrix $S$ for (\ref{20}) is not upper triangular, because in $U$
  the order of $u_1$, ..., $u_k$ is not lexicographical w.r.t. the
  line $l$. Then, Coxeter identity is $-S^{-1}S^T=$ product of the
  $R_j$'s in the order refered to $l$. For example, for $k=3$,
  $S^{-1}S^T= -R_2 R_3 R_1$, since the lexicographical ordering would be
$u_2$, $u_3$, $u_1$. From the identity $S^{-1}S^T=
  (-1)^{k-1}T^k$ it follows a first general relation in the group
$$ 
   T^{k} = (-1)^k \hbox{ product of }R_j \hbox{'s in suitable order } 
$$

Two cases must now be distinguished.
\vskip 0.2 cm
\noindent 
{\bf $k$ odd:}  As a general result \cite{BJL2}, $\det G =0$ if and
only if $V + {1\over 2} $ has an integer eigenvalue. The eigenvalues
of $V$ are ${k-1 \over 2}$, ${k-3 \over 2}$, ..., $-{k-1 \over
2}$. Then, for $k$ odd, $\det G \neq 0$, and  $\phi^{(1)}$, ...,
$\phi^{(k)}$ are a basis. The matrices $R_j$ are 
$$ 
    R_j = \pmatrix{     1     &           &          &   &    &         \cr
                              &\ddots    &          &   &    &         \cr 
                              &           & \ddots   &   &    &         \cr 
                     -2G_{j1} & -2 G_{j2} & ...      & -1& ...& -2G_{jk} \cr
                              &           &      &   & \ddots &\cr
                              &&&&&1\cr
}
$$
where $S^T+S=2G$. 

In concrete examples, we have  ``empirically''  found other relations like 
$$ 
   R_2= p_1(T,R_1)$$
$$
R_2= p_2(T,R_1)$$
$$ \vdots$$
$$
R_k= p_k(T,R_1)$$
where $p_j(T,R_1)$ means a product of the elements $T$ and $R_1$.  We
have also found the relation 
$$ 
         \left(T R_1 \right)^k = -I
$$
We investigated the following cases: 
\vskip 0.2 cm
\noindent
{\small ${\bf CP}^2$ ($k=3$)
 $$
    R_2=TR_1T^{-1},~~~R_3=T~(R_1 R_2 R_1)~T^{-1}
$$ 
    $$ (TR_1)^3 =-I 
                     $$
$$ 
   T^3=-R_2 R_3 R_1 
$$ 
  
\noindent 
${\bf CP}^4$ ($k=5$)
 $$\left\{ \matrix{
             R_2=TR_1T^{-1},~~~R_3=TR_2T^{-1} \cr
             R_4=T~(R_2R_3R_2)~T^{-1},~~~ R_5=T^{-1}~(R_2R_1R_2)~T \cr
} \right.
$$
$$ (TR_1)^5=-I
$$ 
$$ 
    T^5= - R_3R_4R_2R_5R_1 
$$ 

\noindent 
${\bf CP}^6$ ($k=7$)
$$ \left\{ \matrix{ 
                    R _2=TR_1T^{-1},~~~R_3=TR_2T^{-1},~~~R_4=TR_3T^{-1} \cr
                    R_5=T~(R_3R_4R_3)T^{-1},~~~R_6=T(TR_1)^3~R_2~[T(TR_1)^3]^{-1}\cr
                     R_7=T^{-2}(R_3R_2R_3)T^2 \cr
} \right. 
$$
$$ (TR_1)^7=-I$$
$$T^7 = -R_4R_5R_3R_6R_2R_7R_1$$
%
}                       
 Note that one ralation, for example that for $R_k$, can be derived
 from the others, and that just $R_1$, $T$, $-I$ are enought to
 generate the monodromy group in each of the examples. They satisfy
 ({\it in the examples}) the relations: 
$$ 
   R_1^2=(-I~T~R_1)^k =(-I)^2 =I 
$$ 
$$ R_1(-I)~\left((-I)R_1 \right)^{-1}=I,
~~~ T(-I)~\left((-I)T \right)^{-1}=I
$$
The last two relations mean  simply the commutativity of $-I$ with $R_1$ and
$T$. The relations are not only satisfied, but also ``fulfilled''
(namely, $(-I~T~R_1)^n\neq I$ for $n<k$). Now call 
$$ 
   X:=R_1,~~~Y:=-ITR_1,~~~Z=-I$$
These elements generate the monodromy group of ${\bf CP}^{k-1}$ 
with {\it at least} the relations 
  $$ 
     X^2=Y^k=Z^2=1$$
$$
     (ZX)(XZ)^{-1}=1,~~~(ZY)(YZ)^{-1}=1$$
Note that $Z$ generates the cyclic group $C_2$ of order 2. 

If there were no other relations (which we did not find
``empirically''), 
we would  conclude that the monodromy
group of the quantum cohomology of 
${\bf CP}^k$ ({\it in the examples}) is isomorphic to the
direct product
$$ 
  < X,Y~|~X^2=Y^k=1>\times C_2
$$ 
where $ < X,Y~|~X^2=Y^k=1>$ means the group generated by $X$, $Y$ with
relations $X^2=Y^k=1$. 

\vskip 0.2 cm 
\noindent 
{\bf $k$ even:}  Now $\det G=0$, since $V +{1\over 2}$ has  integer
eigenvalues. $G$ has rank $k-1$ and the eigenspace of its eigenvalue 0
has dimension 1. Let $(z^1,...,z^k)^T$ be an eigenvector of
eigenvalue 1. The vector $v:= \sum_{j=1}^k z^j \phi^{(j)}$ is zero,
because 
$$
\left(v, \phi^{(i)} \right)= \sum_{j=1}^k z^j G_{ji}=0~~~\forall i $$
then 
$$ 
   z^1 \phi^{(1)}+z^{(2)}\phi^{(2)}+...+z^k\phi^{(k)} =0
$$
and $k-1$ of the $\phi^{(j)}$'s are linearly independent. The fuchsian
system (\ref{Fuc20}) has a regular (vector) solution $\phi_0(\lambda)=
\sum_{n=0}^d ~\phi_n \lambda^n$, where $\phi_n$ are constant (column)
vectors, and $\phi_d$ is the eigenvector of $V+{1
\over 2}$ relative to the largest integer eigenvalue less or equal to
zero; this eigenvalue  is precisely $-d$ (see \cite{BJL2}). In our case, $d=0$ 
and $\phi_o(\lambda)=\phi_0$, a constant vector. $\phi^{(1)}$, $\phi^{(2)}$, ..., $\phi^{(k-1)}$,  $\phi_0$ is then a
possible choice for a basis of solutions.  

Observe that in the gauge equivalent form $\psi=X\phi$, $\psi_0$ is 
the eigenvector of ${1\over 2}+\hat{\mu}$ with eigenvalue zero. 
Then 
$$ \psi_0=\pmatrix{0\cr
                   \vdots \cr
                 1 \cr
               \vdots \cr
               0 \cr
}~\equiv ~\pmatrix{\partial_1 x\cr
\partial_2 x \cr
\vdots \cr
\vdots \cr
\partial_k x \cr}
$$
where all the entries are zero but the one at position ${k \over
2}+1$. $x$ is the flat coordinate for $(~,~)-\lambda <~,~>$
corresponding to $\psi_0$. Then, we can chose the following flat
coordinates:
$$ x^{1}(\lambda,t),~x^{2}(\lambda,t),~...~,~x^{k-1}(\lambda,t), ~
t^{{k \over2}+1}$$
The monodromy group then acts on a $k-1$ dimensional space.

Let us determine the reduction of $R_1$, $R_2$, ..., $R_k$, $T$ to the
  $k-1$ dimensional space. The entries of $T$ on the vectors
  $\phi^{(j)}$ are:  $T\phi^{(i)}= \sum_{j=1}^k ~T_{ji} \phi^{(j)}$,
  $i=1,...,k$. 
On
  the new basis $\phi^{(1)}$, ..., $\phi^{(k-1)}$, $\phi_0$ 
the matrices are rewritten  
$$ 
    R_j \phi^{(i)}= \phi^{(i)} - 2 G_{ij} \phi^{(j)}
    ~~~i=1,..,k-1~~ ~j\neq k
$$
  $$    R_j \phi_0 = \phi_0 
                ~~~~~~~~~~~~~~~~~~~~~~~~~~~~~~~j\neq k$$
$$ 
R_k \phi^{(i)} = \phi^{(i)} - 2 G_{ik} \left( -{1 \over
z^k}~\sum_{j=1}^{k-1} ~z^j \phi^{(j)} \right)
~~~ i \neq k$$    
$$ 
T \phi^{(i)} =  \sum_{j=1}^{k-1} ~T_{ji} \phi^{(j)} +T_{ki} 
 \left( -{1 \over
z^k}~\sum_{j=1}^{k-1} ~z^j \phi^{(j)} \right)
$$
Then the matrices assume a reduced form $\pmatrix{ * & * & 0 \cr
                                                    * & * & 0 \cr
                                                    0 & 0 & 1 \cr }$

We studied two examples; besides Coxeter identity $T^k=$product of
$R_j$'s, we found relations similar to the case $k$ odd: 
$$ 
  R_2= p_1(T,R_1)$$
$$
R_2= p_2(T,R_1)$$
$$ \vdots$$
$$
R_k= p_k(T,R_1)$$
and
$$ 
         \left(T R_1 \right)^k = I
$$
Namely: 
\vskip 0.2 cm
{\small
\noindent
${\bf CP}^3$ ($k=4$) 
$$ 
   \left\{ \matrix{ 
                   R_2=TR_1T^{-1},~~~ R_3=T R_2T^{-1} \cr
                  R_4 = T^{-1}~(R_2R_1R_2)~T   \cr
} \right.
$$
$$ 
   (TR_1)^4 =I 
$$ 
$$ 
   T^4 =R_3R_2R_4R_1$$

\noindent
${\bf CP}^5$ ($k=6$) 
$$ 
   \left\{ \matrix{ 
                   R_2=TR_1T^{-1},~~~ R_3=T R_2T^{-1} \cr
                  R_4 = TR_3T^{-1},~~~R_5=T~(R_2R_3R_4R_3R_2)~T^{-1}   \cr
                          R_6= T^{-1}~(R_2 R_1 R_2)~T
} \right.
$$
$$ 
   (TR_1)^6=I
$$
$$ 
   T^6= R_4R_3R_5R_2R_6R_1
$$ 
}

The same remarks of $k$ odd hold here. Call 
$$
     X:=R_1,~~~~Y:=R_1T$$
then, if there were no other hiden relations, the monodromy group of
the quantum cohomology of 
${\bf CP}^k$ ({\it in the examples}) would be isomorphic to 
$$ 
 < X,Y,~|~ X^2=Y^k=1>
$$ 

\vskip 0.2 cm 
Note that $ < X,Y,~|~ X^2=Y^k=1>$ is (isomorphic to) the subgroup of
orientation preserving transformations 
of the hyperbolic triangular group $[2,k,\infty]$. 

\vskip 0.2 cm 
\noindent 
{\bf Lemma 10:} {\it The    subgroup of the orientation preserving
transformations of the hyperbolic triangular group $[2,k,\infty]$ is
isomorphic to the subgroup of $PSL(2,{\bf R})$ generated by
$$ 
      \tau \mapsto -{1 \over \tau} 
$$ 	
$$ 
   \tau \mapsto {1 \over 2 \cos {\pi \over k} - \tau} $$ 
$ \tau \in H:= \{ z \in {\bf C} ~|~\Im z >0\}$ 
} \vskip 0.2 cm 
\noindent
{\bf Proof:  }Consider three integers $m_1$, $m_2$, $m_3$ such that 
          $$ {1 \over m_1} +{1 \over m_2} +{1 \over m_3}<1$$  
               In the Bolyai-Lobatchewsky plane $H$, the 
               triangular group $[m_1,m_3,m_3]$   of hyperpolic 
reflections in the sides of hyperbolic triangles of angles ${\pi \over
               m_1}$, ${\pi \over m_2}$, ${\pi \over m_3}$ is generated
               by three reflections $r_1$, $r_2$, $r_3$ satisfying the
relations 
$$ 
   r_1^2=r_2^2=r_3^2=(r_2r_3)^{m_1}=(r_3r_1)^{m_2}=(r_2 r_1)^{m_3}=1
$$
and the subgroup of orientation preserving transformation 
is generated by $X=r_2r_3$, $Y=r_3
r_1$, Then 
       $$    X^{m_1}=Y^{m_2}=(XY)^{m_3}=1$$
For $m_1=2$, $m_2=k$, $m_3 = \infty$, a fundamental triangular region
is $\{ 0< \Re z < \cos {\pi \over k} \} \cap \{ |z|>1 \}$. Then 
$$ 
  r_1(\tau)=- \bar{\tau},~~~r_2(\tau)= {1 \over
  \bar{\tau}},~~~r_3(\tau)=2 \cos {\pi \over k} - \bar {\tau}
$$ 
The bar means complex conjugation. Then 
$$X(\tau)=-{1 \over \tau},~~~ Y(\tau)= { 1 \over 2 \cos {\pi \over k}
-\tau}$$

\rightline{$\Box$}
\vskip 0.2 cm 
\noindent
{\bf Remark :} The orientation preserving transformations
 of $[2,3,\infty]$ are the {\it
modular group} $PSL(2, {\bf Z})$. 

\vskip 0.3 cm 
\noindent
{\bf Theorem 3:} {\it The monodromy group of the quantum cohomology of
${\bf CP}^2$ is isomorphic to 
\be 
     < X,Y,~|~ X^2=Y^3=1>\times C_2 \cong PSL(2,{\bf Z}) \times C_2
\label{Re2}
\ee
The monodromy group of the quantum cohomology of 
${\bf CP}^3$ is isomorphic to 
\be
< X,Y,~|~ X^2=Y^4=1>\cong \hbox{ orient. preserv. transf. of  } [2,4,\infty
]
\label{Re3}
\ee
}

\vskip 0.2 cm
 The theorem for the  case of ${\bf CP}^2$ is already proved in
 \cite{Dub2}. 
\vskip 0.2 cm
\noindent 
{\bf Proof :} a) ${\bf CP}^2$: 
$$ 
  R_1=\pmatrix{ -1&3&3 \cr
                  0&1&0 \cr
               0&0&1  \cr
             }~~~~
   T=\pmatrix{ 0&0&1 \cr 
             -1 & 3 & 3 \cr
              0 & -1 & 0 \cr
           }
$$
and $X=R_1$, $Y= -I T R_1$ and $Z= -I$ satisfy the relations of
(\ref{Re2}). They act on the column vector ${\bf x}
=\pmatrix{x \cr y \cr z \cr}$. The quadratic form $q(x,y,z)={\bf x}^T
G {\bf x}$ is $R_1$ and $T$ - invariant. Then $T$, $R_1$ act on two
dimensional invariant subspaces $q(x,y,z)=$ constant. On each of these
subspaces we introduce new coordinates $\chi\in {\bf R}$ and
$\varphi\in [0,2\pi)$. Let $\tau = e^{\chi} ~e^{i\varphi}$ and 
$$ 
    x= {a \over 2} ( \tau \bar{\tau} -{3\over 2} ( \tau +\bar{\tau})+1
    ) ~{i \over \tau -\bar{\tau}} $$
$$ 
    y= {a \over 2} ( \tau \bar{\tau} -{1\over 2} ( \tau +\bar{\tau})-1
    ) ~{i \over \tau -\bar{\tau}} $$
$$ 
    z= {a \over 2} ( -\tau \bar{\tau} -{1\over 2} ( \tau +\bar{\tau})+1
    ) ~{i \over \tau -\bar{\tau}} $$
$a \in {\bf R}$, $a \neq 0$. Note that $q(x,y,z)=a^2>0$.
 Then, it is easily verified that 
$$ 
 {\bf x}(-{1 \over \tau})= -X~{\bf x}(\tau)
$$
$$ 
 {\bf x}({1 \over 1-\tau})= Y~{\bf x}(\tau)
$$
$$ 
 {\bf x}(\tau, -a)= Z~ {\bf x}(\tau,a)
$$. 
This implies the 1 to 1 correspondence between the generators of the
modular group and $X$ and $Y$. 

\vskip 0.2 cm 
b) Case of ${\bf CP}^3$. 
$$ R_1=\pmatrix{ -1 & 4 & -10 & 0 \cr 
                0 & 1 & 0 & 0 \cr
                 0 & 0 & 1 & 0 \cr
                0 & 0 & 0 & 1 \cr},~~~
T=\pmatrix{0 & 0 & 1 & 0 \cr 
           -1 & 0 & 3 & 0 \cr
            0 & -1 & 3 & 0 \cr
           0 & 0 & 0 & 1 \cr
            }
$$
The matrices are already written on $\phi^{(1)}$, $\phi^{(2)}$, $\phi^{(3)}$,
$\phi_0$. Recall that the monodromy acts only on $x^1,x^2,x^3$,
because the last flat coordinate is $t^3$. This action is given by
the following  
three dimensional matrices,  acting on a three dimensional space of  
vectors  ${\bf x}=
\pmatrix{ x \cr y \cr
z \cr}$  
$$
  r_1= \pmatrix{   -1 & 4 & 10 \cr
0 & 1 & 0 \cr
 0 & 0 & 1 \cr},~~~t:= \pmatrix{ 0 & 0 & 1 \cr
-1 & 0 & 3 \cr
 0 & -1 & 3\cr
}
$$
We redefine $X= r_1$ and $Y=t r_1$, which satisfy the relations
(\ref{Re3}). 
We proceed as above, defining
$$ 
    x= {a } ( \tau \bar{\tau} -{1\over \sqrt{2}} ( \tau
    +\bar{\tau})+{1\over 3}
    ) ~{i \over \tau -\bar{\tau}} $$
$$ 
    y= {a } (-{2 \over 3} \tau \bar{\tau} -{2 \sqrt{2}\over 3} ( \tau
    +\bar{\tau})+{2 \over 3}
    ) ~{i \over \tau -\bar{\tau}} $$
$$ 
    z= {a } ( -{1 \over 3}\tau \bar{\tau} -{\sqrt{2}\over 6} ( \tau
    +\bar{\tau})+{1\over 3}
    ) ~{i \over \tau -\bar{\tau}} $$ 
$a\neq 0$. Note that ${\bf x}^T g {\bf x}= (8/9) a^2$, where $g$ is the $3 \times
   3 $ reduction of $G$. 
It is easily verified that 
$$
  {\bf x}(-{1\over \tau})= -X~ {\bf x}(\tau)
$$ 
 $$
{\bf x}({1\over \sqrt{2}-\tau})= -Y~ {\bf x}(\tau)
$$
which proves the theorem. 

\rightline{$\Box$}



\vskip 1 cm 
\noindent
{\bf APPENDIX 1: Proof of Lemma 5}
\vskip 0.3 cm
 Let us consider the function 
$$ 
  g(z)= C~ \int_{-c-i \infty}^{-c+i\infty} ds ~\Gamma^k (-s) ~e^{i
  \pi f s } ~z^{ks} 
$$
where $C$ and $f$ are constants to be determined later, $c>0$. The
path of integration is a vertical line through $-c$.

{\bf a)} Domain of definition.  We use Stirling formula 
$$ 
  \Gamma(-s)= e^s e^{-(s+{1 \over 2}) \log(-s)} \sqrt{ 2 \pi} ( 1 +
  O(1/s))  ~~~~~ s \to \infty, ~~~|\arg(-s)|< \pi
$$
where $ \log(-s) = \log(s) + i \pi $. 

 The  integrand is 
$$ 
  \Gamma^k(-s)e^{i\pi f s} z^{ks} \sim e^{ks(1 + \log z)}e^{-k(s+{1
  \over 2})\log s} e^{-ik \pi (s+{1 \over 2}) + i \pi s f} ~~~s \to
  \infty
$$
Now let $s=-c + i \eta$, $ds = i d \eta$. The integral in $\eta$ is on
the real axis. The dominant part in the
integrand is 
$$ 
   e^{ \eta \left[k \arg(i\eta -c ) + k \pi - k \arg z - \pi f \right]
   } 
$$
The condition of (uniform) convergence of the integral is 
$$ 
      -{\pi \over 2} -f {\pi \over k} < \arg z< {\pi \over 2} - f {\pi
      \over k}
$$ 

\vskip 0.2 cm
{\bf b)} $g$ solves (\ref{30}). In fact (for simplicity $C=1$ here): 
$$ (z \partial_z)^k g(z)= k^k \int_{-c -i\infty}^{-c+i \infty} ds ~
\Gamma^k(-s) e^{i \pi f s } s^k z^{ks}= 
(-1)^k k^k \int_{-c -i\infty}^{-c+i \infty} ds ~\Gamma^k(1-s)e^{i \pi
f s} z^{ks}
$$
where we have used the identity $ s \Gamma(-s)= - \Gamma(1-s)$. Now
let $t= s-1$: 
$$
(z \partial_z)^k g(z)= (-1)^k e^{i \pi f} (zk)^k  \int_{-c-1 -i\infty}^{-c-1+i \infty}
dt ~ \Gamma^{k}(-t) e^{i \pi f t} z^{kt}
$$
$$ = (-1)^k e^{i\pi f} (kz)^k g(z)$$ 
Now we impose $(-1)^k e^{i \pi f}=1$, namely 
$$ 
                    f = k + 2 m~~~~~~~m \in {\bf Z} 
$$ 
Point b) is proved. 

\vskip 0.2 cm
{\bf c)} Asymptotic behaviour of $g(z)$. We use Laplace method for analytic
functions (the so called  ``steepest descent method''). Put $C=1$. 
By Stirling's formula
$$ 
 g(z) = (2 \pi)^{k/2} e^{-i {\pi \over 2} k} ~ \int_{-c-i
 \infty}^{-c+i\infty} ds~ e^{\phi(s)}
$$
$$
~~~~~\phi(s)=s~[ k(1 + \log z) + i \pi f - i \pi k] - k\left( s + {1
\over 2}\right) \log s + O\left( {1 \over s} \right)
$$
The the stationary point of $\phi(s)$ is 
$$ 
   s_0 = z e ^{-i { \pi \over k } (k-f)} -{1 \over 2} + O\left({1
   \over z} \right)  ~~~~~z \to \infty
$$
In the hypothesis $z \to \infty$ a straightforward computation gives 
$$ 
   \phi(s_0)\sim k z e^{-i{\pi \over k} (k-f)}-{k \over 2} \log z + i
   {\pi \over 2} (k-f)~~~~~{d^2 \phi(s_0)\over ds^2}\sim -{k \over z}
   e^{i {\pi \over k}(k-f)}
$$
 We are ready to apply the steepest descent method. We deform the path
 of integration in such a way that it passes through $s_0$. Let us
 call it $\gamma$.
$$ 
  g(z)= (2 \pi)^{k \over 2} e^{-i {\pi \over 2} k} e^{\phi(s_0)} ~ 
\int_{\gamma} ds ~ e^{\phi(s)-\phi(s_0)}
$$
$$
\sim  (2 \pi)^{k \over 2} e^{-i {\pi \over 2} k} e^{\phi(s_0)}~
\int_{\gamma} ds~ e^{{1 \over 2} {d^2 \phi(s_0) \over ds^2} (s-s_0)^2}
$$
Let us divide $\gamma$ in two paths: $\gamma_1$ from $s_0$ to $+ i
\infty$ and $\gamma_2$ from $- i \infty$ to $s_0$. The integrals
becomes the sum of two integrals. In $\int_{\gamma_1}$ we change
variable. Let $\tau>0$ and  
 $$ 
      - \tau^2 =  {1 \over 2} {d^2 \phi(s_0) \over ds^2} (s-s_0)^2= 
                 -{ k \over 2 z} \left[ e^{i {\pi \over k} (k-f)}
                (s-s_0)^2 \right]
$$ 
 Then 
$$ 
\int_{\gamma_1} e^{ {d^2 \phi(s_0) \over ds^2}  (s -s_0)^2}= 
{1 \over 2} \sqrt{ \pi \over k} ( 2 z )^{1\over 2} e^{ - i {\pi \over
2 k} (k-f) -i \pi} + O(e^{- \alpha |z|})~~~~~\alpha>0$$ 
Note that $ \int_{\gamma_1} =\int_{\gamma_2}$. Now, recalling that $f=
k + 2 m$, we conclude that 
$$ 
g(z) \sim ( 2 \pi) ^{k+1 \over 2} e^{-i(m+1) \pi} e^{ -i { \pi \over
2}k} ~ ~ { 1 \over \sqrt{k}}{ e^{i{\pi \over k}m} \over z^{k-1 \over
2}}~ \exp( k e^{i{2\pi \over k}m} z)~~~~~z\to \infty
$$
If we choose 
$$ 
                         m= n-1 +kl
~~~~l \in {\bf Z}
$$ 
and $C= \left[ (2\pi)^{k+1 \over 2} e^{i \pi ( l -{k\over 2} -n-kl)}
\right]^{-1}$ 
we have 
$$ 
   g(z) \sim {1 \over \sqrt{k}} { e^{i { \pi \over k}(n-1)} \over
   z^{k-1\over 2}}
                        \exp( k  e^{i { 2\pi \over k}(n-1)} z)
$$ 
Then $g$ is a solution $\varphi^{(n)}$ on the domain 
$$ 
         -{3 \pi \over 2} -2(n-1) { \pi \over k} -2 \pi l< 
                 \arg z 
                         <    -{ \pi \over 2} -2(n-1) { \pi \over k}
                 -2 \pi l
$$
\vskip 0.2 cm
{\bf d)} We now determine the domain where the analytic extension of $g(z)$
still has the above asymptotic behaviour. From Lemma 4 we derived that
the first entries of the $n(k)^{th}$ columns of $Y_L$ and $Y_R$  are
equal and coincide with $\varphi^{(n(k))}_L(z)\equiv
\varphi^{(n(k))}_R(z)$
 times $z^{k-1
\over 2}$, and $\varphi^{(n(k))}$ has the estabilished asymptotic
behaviour on the enlarged domains 
$$ 
   -{\pi \over k} - \pi < \arg z < \pi + {\pi \over k} 
$$
$$ 
     - \pi < \arg z < \pi + {2\pi \over k} 
$$ 

 If in $g(z)$ we chose $n=n(k)$ and $l=-1$ the integral representation
 holds for $-{\pi \over 2} <\arg z < {\pi \over 2}$ for $k$ even, and
 for   $-{\pi \over 2}+{\pi \over k} < \arg z < {\pi \over 2}+{\pi \over
 k}$ for $k$ odd. Then its analytic continuation is precisely 
$\varphi^{(n(k))}_L(z)\equiv \varphi^{(n(k))}_R(z)$.

\vskip 0.2 cm 
{\bf e)} We prove the identity (\ref{60}). Observe that $(z \partial_z)^k
\varphi = (kz)^k \varphi$ is invariant for $z \mapsto z~e^{2\pi i
\over k}$. A generic solution can be reperesented near $z=0$ as 
$$
       \varphi(z)= \sum_{m=0}^{\infty} {z^{km}\over (m!)^k}~\left[ 
           a^{(1)}_m+a^{(2)}_m\log z + ... + a^{(k)}_m \log^{k-1}z
                      \right]
$$
It follows that the operator $(A \varphi)(z)=\varphi(z e^{2 \pi i
\over k})$ has eigenvalues 1, because on the basis obtained with 
$(a^{(1)}_0=1,~ a^{(2)}_0=...=a^{(k)}_0=0)$, $( a^{(1)}_0=0,
~a^{(2)}_0=1, 
~a^{(3)}_0=...= a^{(k)}_0=0)$, ..., $(a^{(1)}_0=...= a^{(k-1)}_0=0,
~a^{(k)}_0=1)$ it is represented by a lower triangular matrix having
1's on the diagonal. Then $(A-1)^k=0$ and 
$$ 
         (A-1)^k g(z)=0
$$ 
is precisely our identity.   Lemma 4 is
 proved.                  

\rightline{$\Box$}


\vskip 1 cm
\noindent
{\bf APPENDIX 2: }  
\vskip 0.2 cm

First we give $\Phi_R$ and $\Phi_L$. 
{\bf $k$ odd: }
{\tiny
$$
    \Phi_R(z)^T= \left[ \matrix{ 
                                 (-1)^{k-1 \over 2}\left[ 
                          g(ze^{-{2 \pi i \over k} ({k-1 \over 2})})
-\bin{k}{1} g(ze^{-{2 \pi i \over k} ({k-3 \over 2})})+ ...
+\bin{k}{k-1} g(ze^{{2 \pi i \over k} ({k-1 \over 2})})
                                 \right]                    \cr
\cr
\vdots \cr
\cr
 g(ze^{-{4\pi i \over k} })-\bin{k}{1} g(ze^{-{2\pi i \over k} })
 + \bin{k}{2} g(z)-  \bin{k}{3} g(ze^{{2\pi i \over k} })
 +   \bin{k}{4} g(ze^{{4\pi i \over k} })                      \cr
   - g(ze^{-{2\pi i \over k} })  + \bin{k}{1} g(z)- \bin{k}{2} 
g(ze^{{2\pi i \over k} })                                       \cr
g(z)        \cr
 - g(ze^{{2\pi i \over k} })          \cr
  g(ze^{{4\pi i \over k} })                  \cr
\cr
\vdots               \cr
\cr
 (-1)^{k-1 \over 2}  g(ze^{{2\pi i \over k} ({k-1\over 2})}) \cr
                 } \right] 
$$
\vskip 0.2 cm
$$
\Phi_L(z)^T= \left[ \matrix{  
                            (-1)^{k-1 \over 2}  
                              g(ze^{-{2\pi i \over k} ({k-1\over 2})})
                                                           \cr
                         -(-1)^{k-1 \over 2}  
                              g(ze^{-{2\pi i \over k} ({k-3\over 2})})
                                                             \cr
\cr
\vdots  \cr
\cr
 -g(ze^{-{2\pi i \over k} }) \cr
g(z)        \cr
- g(ze^{{2\pi i \over k}} )+\bin{k}{k-1}g(z) \cr
 g(ze^{{4\pi i \over k} })- \bin{k}{k-1} g(ze^{{2\pi i \over k} })+
 \bin{k}{k-2} g(z)- \bin{k}{k-3} g(ze^{-{2\pi i \over k} })        \cr
\cr
\vdots              \cr
\cr
(-1)^{k-1\over 2} \left[    g(ze^{{2\pi i \over k} ({k-1\over 2})})-
                              \bin{k}{k-1}
                               g(ze^{{2\pi i \over k} ({k-3\over 2})})
+  ...-  \bin{k}{2}  g(ze^{-{2\pi i \over k} ({k-3\over 2})})             \right]
                                                                 \cr
}\right]
$$
}

\vskip 0.3 cm
\noindent
{\bf $k$ even:}  
{\tiny
$$
    \Phi_R(z)^T= \left[ \matrix{ 
                                 (-1)^{k \over 2}\left[ 
                          g(ze^{-{i \pi } })
-\bin{k}{1} g(ze^{-i \pi +i{2 \pi  \over k} })+ ...
+\bin{k}{k-1} g(ze^{i(\pi-{2 \pi \over k})})
                                 \right]                    \cr
\cr
\vdots \cr
\cr
 g(ze^{-{4\pi i \over k} })-\bin{k}{1} g(ze^{-{2\pi i \over k} })
 + \bin{k}{2} g(z)-  \bin{k}{3} g(ze^{{2\pi i \over k} })          \cr
   - g(ze^{-{2\pi i \over k} })  + \bin{k}{1} g(z)     \cr
g(z)        \cr
 - g(ze^{{2\pi i \over k} })          \cr
  g(ze^{{4\pi i \over k} })                  \cr
\cr
\vdots               \cr
\cr
 (-1)^{{k \over 2}-1}  g(ze^{{2\pi i \over k} ({k\over 2}-1)}) \cr
                 } \right] 
$$
\vskip 0.2 cm
$$
\Phi_L(z)^T= \left[ \matrix{  
                            (-1)^{k \over 2}  
                              g(ze^{-i\pi  })
                                                           \cr
\vdots  \cr  \cr
                              g(ze^{-i{4\pi  \over k} })
                                                             \cr
 -g(ze^{-{2\pi i \over k} }) \cr
g(z)        \cr
- g(ze^{{2\pi i \over k}}
                              )+\bin{k}{k-1}g(z)-\bin{k}{k-2}
                  g(ze^{-i{2\pi\over k}}) \cr
 g(ze^{{4\pi i \over k} })- \bin{k}{k-1} g(ze^{{2\pi i \over k} })+
 \bin{k}{k-2} g(z)- \bin{k}{k-3} g(ze^{-{2\pi i \over k}
                              })+\bin{k}{k-4} g(ze^{-i{4\pi \over k}})        \cr
\cr
\vdots              \cr
\cr
(-1)^{{k\over 2}-1} \left[    g(ze^{i\pi -i{2\pi  \over k} })-
                              \bin{k}{k-1}
                               g(ze^{i\pi -i{4\pi  \over k}})
+  ...-  \bin{k}{2}  g(ze^{-i\pi +i{2\pi i \over k}} )             \right]
                                                                 \cr
}\right]
$$
}

\vskip 0.3 cm
  We give all the matrices of interest up to $k=10$. 
$S_{upper}$ is
$PSP^{-1}$.  $A$ stends for $A^{\beta}$, $A^{\prime}$ for
  $A^{\beta^{\prime}}$. $S^{\beta}=AS_{upper}A^T$, $S^{
  \beta^{\prime}}= A^{\prime} S_{upper}^{-1}
  [A^{\prime}]^T$. $\sigma_{i,i+1}:=\beta_{i,i+1}^{-1}$

{\small

${\bf CP}^2$

$$
 K_{12} =  \pmatrix{
1 & 0 & 0 \cr
-3 & 1 & 0 \cr
0 & 0 & 1\cr
}
~~~
 K_{13} :=  \pmatrix{
1 & 0 & 0 \cr
0 & 1 & 0 \cr
-3 & 0 & 1\cr
}
~~~
 K_{32} :=  \left[ 
\matrix{
1 & 0 & 0 \cr
0 & 1 & 3 \cr
0 & 0 & 1 \cr
}
 \right] 
$$
$$
T :=  \left[ 
\matrix{
0 & 0 & 1 \cr
-1 & 3 & 3 \cr
0 & -1 & 0 \cr
}
 \right] 
$$
$$
S :=  \left[ 
\matrix{
1 & 0 & 0 \cr
-3 & 1 & 3 \cr
-3 & 0 & 1 \cr
}
 \right]
~~~PSP^{-1}=\left[ 
\matrix{
1 & 3 & -3 \cr
0 & 1 & -3 \cr
0 & 0 & 1 \cr
}
 \right]
$$

\vskip 0.2 cm

  ${\bf CP}^3$

$$
 K_{42} :=  \left[ 
\matrix{
1 & 0 & 0 & 0 \cr
0 & 1 & 0 & 6 \cr
0 & 0 & 1 & 0 \cr
0 & 0 & 0 & 1\cr
}
 \right] 
~~~
 K_{43} :=  \left[ 
\matrix{
1 & 0 & 0 & 0 \cr
-4 & 1 & 0 & 0 \cr
0 & 0 & 1 & 4 \cr
0 & 0 & 0 & 1\cr
}\right]
$$
$$
T :=  \left[ 
\matrix{
0 & 0 & 0 & 1 \cr
-1 & 0 & 6 & 4 \cr
0 & -1 & 4 & 0 \cr
0 & 0 & -1 & 0 \cr
}
 \right]
~~~
 T_1 :=  \left[ 
\matrix{
0 & 0 & 1 & 0 \cr
-1 & 0 & 3 & 0 \cr
0 & -1 & 3 & 0 \cr
0 & 0 & 0 & 1 \cr
}
 \right] 
$$
$$
 S :=  \left[ 
\matrix{
1 & 0 & 0 & 0 \cr
-4 & 1 & 0 & 6 \cr
10 & -4 & 1 & -20 \cr
-4 & 0 & 0 & 1 \cr
}
 \right] 
$$
$$
{\it P    } :=  \left[ 
\matrix{
0 & 0 & 1 & 0 \cr
0 & 1 & 0 & 0 \cr
0 & 0 & 0 & 1 \cr
1 & 0 & 0 & 0 \cr
           }\right]
 ~~~
 S_{upper}= \left[ 
\matrix{
1 & -4 & -20 & 10 \cr
0 & 1 & 6 & -4 \cr
0 & 0 & 1 & -4 \cr
0 & 0 & 0 & 1 \cr
          }
 \right] 
  $$
$$
A :=  \left[ 
\matrix{
0 & 1 & 0 & 0 \cr
1 & 4 & 0 & 0 \cr
0 & 0 & 1 & 0 \cr
0 & 0 & 0 & 1 \cr
}
 \right]
~~~A^{\prime}=\left[\matrix{
1&0&0&0\cr
0&1&0&0\cr
0&0&-4&1\cr
0&0&1&0 \cr
}
\right]
~~~
 S^{\beta}=S^{\beta^{\prime}} :=  \left[ 
\matrix{
1 & 4 & 6 & -4 \cr
0 & 1 & 4 & -6 \cr
0 & 0 & 1 & -4 \cr
0 & 0 & 0 & 1 \cr
           }
 \right] 
$$
where $\beta=\beta_{12}$

 \vskip 0.2 cm                   
                     
${\bf CP}^4$

$$
K_{52} :=  \left[ 
\matrix{
1 & 0 & 0 & 0 & 0 \cr
-5 & 1 & 0 & 0 & 0 \cr
0 & 0 & 1 & 0 & 10 \cr
0 & 0 & 0 & 1 & 0 \cr
0 & 0 & 0 & 0 & 1 \cr
           }
 \right] 
  ~~~
 K_{53} :=  \left[ 
\matrix{                 
1 & 0 & 0 & 0 & 0 \cr
0 & 1 & 0 & 0 & 10 \cr
0 & 0 & 1 & 5 & 0 \cr
0 & 0 & 0 & 1 & 0 \cr
0 & 0 & 0 & 0 & 1 \cr
           }
 \right] 
$$
        $$
T :=  \left[ 
\matrix{                  
0 & 0 & 0 & 0 & 1 \cr
-1 & 0 & 0 & 10 & 5 \cr
0 & -1 & 5 & 10 & 0 \cr
0 & 0 & -1 & 0 & 0 \cr
0 & 0 & 0 & -1 & 0 \cr
           }
 \right] 
~~~
S= \left[ 
\matrix{             
1 & 0 & 0 & 0 & 0 \cr
-5 & 1 & 0 & 0 & 10 \cr
15 & -5 & 1 & 5 & -40 \cr
40 & -10 & 0 & 1 & -95 \cr
-5 & 0 & 0 & 0 & 1 \cr
          }
 \right] 
$$   
 $$
P :=  \left[ 
\matrix{               
0 & 0 & 1 & 0 & 0 \cr
0 & 0 & 0 & 1 & 0 \cr
0 & 1 & 0 & 0 & 0 \cr
0 & 0 & 0 & 0 & 1 \cr
1 & 0 & 0 & 0 & 0 \cr
          }
 \right] 
~~~
S_{upper} :=  \left[ 
\matrix{                 
1 & 5 & -5 & -40 & 15 \cr
0 & 1 & -10 & -95 & 40 \cr
0 & 0 & 1 & 10 & -5 \cr
0 & 0 & 0 & 1 & -5 \cr
0 & 0 & 0 & 0 & 1 \cr
           }
 \right] 
$$
           
          $$
A :=  \left[ 
\matrix{              
0 & 0 & 1 & 0 & 0 \cr
1 & 0 & 5 & 0 & 0 \cr
0 & 1 & 10 & 0 & 0 \cr
0 & 0 & 0 & 1 & 0 \cr
0 & 0 & 0 & 0 & 1 \cr
          }
 \right] 
~~~
A^{\prime}:=\left[ 
\matrix{
0 & 1 & 0 & 0 & 0 \cr
1 & 5 & 0 & 0 & 0 \cr
0 & 0 & 1 & 0 & 0 \cr
0 & 0 & 0 & -5 & 1 \cr
0 & 0 & 0 & 1 & 0 \cr
}
 \right] 
$$
$$
 S^{\beta} =S^{\beta^{\prime}}:=  \left[ 
\matrix{             
1 & 5 & 10 & 10 & -5 \cr
0 & 1 & 5 & 10 & -10 \cr
0 & 0 & 1 & 5 & -10 \cr
0 & 0 & 0 & 1 & -5 \cr
0 & 0 & 0 & 0 & 1 \cr
          }
 \right] 
$$
where $\beta=\beta_{23}\beta_{12}$, $\beta^{\prime}=\beta_{12}\sigma_{45}$.

         \vskip 0.2 cm         
${\bf CP}^5$

$$
 K_{62} :=  \left[ 
\matrix{                 
1 & 0 & 0 & 0 & 0 & 0 \cr
0 & 1 & 0 & 0 & 0 & 15 \cr
0 & 0 & 1 & 0 & 15 & 0 \cr
0 & 0 & 0 & 1 & 0 & 0 \cr
0 & 0 & 0 & 0 & 1 & 0 \cr
0 & 0 & 0 & 0 & 0 & 1 \cr
           }
 \right] 
~~~
 K_{63} :=  \left[ 
\matrix{                
1 & 0 & 0 & 0 & 0 & 0 \cr
-6 & 1 & 0 & 0 & 0 & 0 \cr
0 & 0 & 1 & 0 & 0 & 20 \cr
0 & 0 & 0 & 1 & 6 & 0 \cr
0 & 0 & 0 & 0 & 1 & 0 \cr
0 & 0 & 0 & 0 & 0 & 1 \cr
           }
 \right] 
$$  
  $$
T :=  \left[ 
\matrix{                  
0 & 0 & 0 & 0 & 0 & 1 \cr
-1 & 0 & 0 & 0 & 15 & 6 \cr
0 & -1 & 0 & 15 & 20 & 0 \cr
0 & 0 & -1 & 6 & 0 & 0 \cr
0 & 0 & 0 & -1 & 0 & 0 \cr
0 & 0 & 0 & 0 & -1 & 0 \cr
           }
 \right] 
~~~
S :=  \left[ 
\matrix{                 
1 & 0 & 0 & 0 & 0 & 0 \cr
-6 & 1 & 0 & 0 & 0 & 15 \cr
21 & -6 & 1 & 0 & 15 & -70 \cr
-56 & 21 & -6 & 1 & -84 & 210 \cr
105 & -20 & 0 & 0 & 1 & -294 \cr
-6 & 0 & 0 & 0 & 0 & 1 \cr
           }
 \right] 
$$
            
                  $$
T_1 :=  \left[ 
\matrix{                
0 & 0 & 0 & 0 & 1 & 0 \cr
-1 & 0 & 0 & 0 & 5 & 0 \cr
0 & -1 & 0 & 15 & 25 & 0 \cr
0 & 0 & -1 & 6 & 1 & 0 \cr
0 & 0 & 0 & -1 & -1 & 0 \cr
0 & 0 & 0 & 0 & 0 & 1 \cr
           }
 \right] 
$$
          $$
P :=  \left[ 
\matrix{                 
0 & 0 & 0 & 1 & 0 & 0 \cr
0 & 0 & 1 & 0 & 0 & 0 \cr
0 & 0 & 0 & 0 & 1 & 0 \cr
0 & 1 & 0 & 0 & 0 & 0 \cr
0 & 0 & 0 & 0 & 0 & 1 \cr
1 & 0 & 0 & 0 & 0 & 0 \cr
           }
 \right] 
$$
                  $$
A :=  \left[ 
\matrix{                 
0 & 0 & 0 & 1 & 0 & 0 \cr
0 & 1 & 0 & 6 & 0 & 0 \cr
1 & 6 & 0 & 15 & 0 & 0 \cr
0 & 0 & 1 & 20 & 0 & 0 \cr
0 & 0 & 0 & 0 & 1 & 0 \cr
0 & 0 & 0 & 0 & 0 & 1 \cr
           }
 \right] 
~~~ A^{\prime}:=  \left[ 
\matrix{
1 & 0 & 0 & 0 & 0 & 0 \cr
0 & 1 & 0 & 0 & 0 & 0 \cr
0 & 0 & -20 & 1 & 0 & 0 \cr
0 & 0 & -15 & 0 & -6 & 1 \cr
0 & 0 & 6 & 0 & 1 & 0 \cr
0 & 0 & 1 & 0 & 0 & 0 \cr
}
 \right] 
$$
$$     S^{\beta} :=  \left[ 
\matrix{                 
1 & 6 & 15 & 20 & 15 & -6 \cr
0 & 1 & 6 & 15 & 20 & -15 \cr
0 & 0 & 1 & 6 & 15 & -20 \cr
0 & 0 & 0 & 1 & 6 & -15 \cr
0 & 0 & 0 & 0 & 1 & -6 \cr
0 & 0 & 0 & 0 & 0 & 1 \cr
           }
 \right] ~~~S^{\beta^{\prime}} :=  \left[ 
\matrix{
1 & 6 & 15 & 20 & -15 & -6 \cr
0 & 1 & 6 & 15 & -20 & -15 \cr
0 & 0 & 1 & 6 & -15 & -20 \cr
0 & 0 & 0 & 1 & -6 & -15 \cr
0 & 0 & 0 & 0 & 1 & 6 \cr
0 & 0 & 0 & 0 & 0 & 1 \cr
}
 \right] 
$$
where $\beta= \beta_{12}(\beta_{34}  \beta_{23}    \beta_{12})$,
$\beta^{\prime}= [(\sigma_{34}\sigma_{56})\sigma_{45}]\sigma_{56}$.              
\vskip 0.2 cm

${\bf CP}^6$

$$
K_{72} :=  \left[ 
\matrix{                 
1 & 0 & 0 & 0 & 0 & 0 & 0 \cr
0 & 1 & 0 & 0 & 0 & 0 & 21 \cr
0 & 0 & 1 & 0 & 0 & 35 & 0 \cr
0 & 0 & 0 & 1 & 7 & 0 & 0 \cr
0 & 0 & 0 & 0 & 1 & 0 & 0 \cr
0 & 0 & 0 & 0 & 0 & 1 & 0 \cr
0 & 0 & 0 & 0 & 0 & 0 & 1 \cr
           }
 \right] 
~~~
 K_{73} :=  \left[ 
\matrix{                
1 & 0 & 0 & 0 & 0 & 0 & 0 \cr
-7 & 1 & 0 & 0 & 0 & 0 & 0 \cr
0 & 0 & 1 & 0 & 0 & 0 & 35 \cr
0 & 0 & 0 & 1 & 0 & 21 & 0 \cr
0 & 0 & 0 & 0 & 1 & 0 & 0 \cr
0 & 0 & 0 & 0 & 0 & 1 & 0 \cr
0 & 0 & 0 & 0 & 0 & 0 & 1 \cr
           }
 \right] 
$$

$$
T :=  \left[ 
\matrix{                  
0 & 0 & 0 & 0 & 0 & 0 & 1 \cr
-1 & 0 & 0 & 0 & 0 & 21 & 7 \cr
0 & -1 & 0 & 0 & 35 & 35 & 0 \cr
0 & 0 & -1 & 7 & 21 & 0 & 0 \cr
0 & 0 & 0 & -1 & 0 & 0 & 0 \cr
0 & 0 & 0 & 0 & -1 & 0 & 0 \cr
0 & 0 & 0 & 0 & 0 & -1 & 0 \cr
           }
 \right] 
~~~
 S :=  \left[ 
\matrix{               
1 & 0 & 0 & 0 & 0 & 0 & 0 \cr
-7 & 1 & 0 & 0 & 0 & 0 & 21 \cr
28 & -7 & 1 & 0 & 0 & 35 & -112 \cr
-84 & 28 & -7 & 1 & 7 & -224 & 378 \cr
-378 & 112 & -21 & 0 & 1 & -728 & 1638 \cr
224 & -35 & 0 & 0 & 0 & 1 & -728 \cr
-7 & 0 & 0 & 0 & 0 & 0 & 1 \cr
           }
 \right] 
$$

$$
P :=  \left[ 
\matrix{                  
0 & 0 & 0 & 1 & 0 & 0 & 0 \cr
0 & 0 & 0 & 0 & 1 & 0 & 0 \cr
0 & 0 & 1 & 0 & 0 & 0 & 0 \cr
0 & 0 & 0 & 0 & 0 & 1 & 0 \cr
0 & 1 & 0 & 0 & 0 & 0 & 0 \cr
0 & 0 & 0 & 0 & 0 & 0 & 1 \cr
1 & 0 & 0 & 0 & 0 & 0 & 0 \cr
           }
 \right] 
$$

$$
A :=  \left[ 
\matrix{                  
0 & 0 & 0 & 0 & 1 & 0 & 0 \cr
0 & 0 & 1 & 0 & 7 & 0 & 0 \cr
1 & 0 & 7 & 0 & 21 & 0 & 0 \cr
0 & 1 & 21 & 0 & 35 & 0 & 0 \cr
0 & 0 & 0 & 1 & 35 & 0 & 0 \cr
0 & 0 & 0 & 0 & 0 & 1 & 0 \cr
0 & 0 & 0 & 0 & 0 & 0 & 1 \cr
           }
 \right] 
~~~A^{\prime} :=  \left[ 
\matrix{
0 & 1 & 0 & 0 & 0 & 0 & 0 \cr
1 & 7 & 0 & 0 & 0 & 0 & 0 \cr
0 & 0 & 1 & 0 & 0 & 0 & 0 \cr
0 & 0 & 0 & -35 & 1 & 0 & 0 \cr
0 & 0 & 0 & -21 & 0 & -7 & 1 \cr
0 & 0 & 0 & 7 & 0 & 1 & 0 \cr
0 & 0 & 0 & 1 & 0 & 0 & 0 \cr
           }
 \right] 
$$
$$
     S^{\beta} :=  \left[ 
\matrix{               
1 & 7 & 21 & 35 & 35 & 21 & -7 \cr
0 & 1 & 7 & 21 & 35 & 35 & -21 \cr
0 & 0 & 1 & 7 & 21 & 35 & -35 \cr
0 & 0 & 0 & 1 & 7 & 21 & -35 \cr
0 & 0 & 0 & 0 & 1 & 7 & -21 \cr
0 & 0 & 0 & 0 & 0 & 1 & -7 \cr
0 & 0 & 0 & 0 & 0 & 0 & 1 \cr
           }
 \right] 
~~~S^{\beta^{\prime}}:= \left[ 
\matrix{
1 & 7 & 21 & 35 & 35 & -21 & -7 \cr
0 & 1 & 7 & 21 & 35 & -35 & -21 \cr
0 & 0 & 1 & 7 & 21 & -35 & -35 \cr
0 & 0 & 0 & 1 & 7 & -21 & -35 \cr
0 & 0 & 0 & 0 & 1 & -7 & -21 \cr
0 & 0 & 0 & 0 & 0 & 1 & 7 \cr
0 & 0 & 0 & 0 & 0 & 0 & 1 \cr
}
 \right] 
$$
where $\beta=  (\beta_{23}\beta_{12})(   \beta_{45} \beta_{34}  \beta_{23}
  \beta_{12})$,
  $\beta^{\prime}=\beta_{12}[(\sigma_{45}\sigma_{67})\sigma_{56}]
\sigma_{67}$.            

  \vskip 0.2 cm                 
                      
                                       ${\bf CP}^7$

$$
 K_{82} :=  \left[ 
\matrix{                 
1 & 0 & 0 & 0 & 0 & 0 & 0 & 0 \cr
0 & 1 & 0 & 0 & 0 & 0 & 0 & 28 \cr
0 & 0 & 1 & 0 & 0 & 0 & 70 & 0 \cr
0 & 0 & 0 & 1 & 0 & 28 & 0 & 0 \cr
0 & 0 & 0 & 0 & 1 & 0 & 0 & 0 \cr
0 & 0 & 0 & 0 & 0 & 1 & 0 & 0 \cr
0 & 0 & 0 & 0 & 0 & 0 & 1 & 0 \cr
0 & 0 & 0 & 0 & 0 & 0 & 0 & 1 \cr
           }
 \right] 
~~~
 K_{83} :=  \left[ 
\matrix{                  
1 & 0 & 0 & 0 & 0 & 0 & 0 & 0 \cr
-8 & 1 & 0 & 0 & 0 & 0 & 0 & 0 \cr 
0 & 0 & 1 & 0 & 0 & 0 & 0 & 56 \cr
0 & 0 & 0 & 1 & 0 & 0 & 56 & 0 \cr
0 & 0 & 0 & 0 & 1 & 8 & 0 & 0 \cr
0 & 0 & 0 & 0 & 0 & 1 & 0 & 0 \cr
0 & 0 & 0 & 0 & 0 & 0 & 1 & 0 \cr
0 & 0 & 0 & 0 & 0 & 0 & 0 & 1 \cr
           }
 \right] 
$$

$$
T :=  \left[ 
\matrix{                  
0 & 0 & 0 & 0 & 0 & 0 & 0 & 1 \cr
-1 & 0 & 0 & 0 & 0 & 0 & 28 & 8 \cr
0 & -1 & 0 & 0 & 0 & 70 & 56 & 0 \cr
0 & 0 & -1 & 0 & 28 & 56 & 0 & 0 \cr
0 & 0 & 0 & -1 & 8 & 0 & 0 & 0 \cr
0 & 0 & 0 & 0 & -1 & 0 & 0 & 0 \cr
0 & 0 & 0 & 0 & 0 & -1 & 0 & 0 \cr
0 & 0 & 0 & 0 & 0 & 0 & -1 & 0 \cr
           }
 \right] 
$$
$$
S :=  \left[ 
\matrix{                 
1 & 0 & 0 & 0 & 0 & 0 & 0 & 0 \cr
-8 & 1 & 0 & 0 & 0 & 0 & 0 & 28 \cr
36 & -8 & 1 & 0 & 0 & 0 & 70 & -168 \cr
-120 & 36 & -8 & 1 & 0 & 28 & -504 & 630 \cr
330 & -120 & 36 & -8 & 1 & -216 & 2100 & -1848 \cr
-1512 & 378 & -56 & 0 & 0 & 1 & -3912 & 7476 \cr
420 & -56 & 0 & 0 & 0 & 0 & 1 & -1560 \cr
-8 & 0 & 0 & 0 & 0 & 0 & 0 & 1 \cr
           }
 \right] 
$$

$$
P :=  \left[ 
\matrix{                
0 & 0 & 0 & 0 & 1 & 0 & 0 & 0 \cr
0 & 0 & 0 & 1 & 0 & 0 & 0 & 0 \cr
0 & 0 & 0 & 0 & 0 & 1 & 0 & 0 \cr
0 & 0 & 1 & 0 & 0 & 0 & 0 & 0 \cr
0 & 0 & 0 & 0 & 0 & 0 & 1 & 0 \cr
0 & 1 & 0 & 0 & 0 & 0 & 0 & 0 \cr
0 & 0 & 0 & 0 & 0 & 0 & 0 & 1 \cr
1 & 0 & 0 & 0 & 0 & 0 & 0 & 0 \cr
           }
 \right] 
~~~
 S_{upper} :=  \left[ 
\matrix{                  
1 & -8 & -216 & 36 & 2100 & -120 & -1848 & 330 \cr
0 & 1 & 28 & -8 & -504 & 36 & 630 & -120 \cr
0 & 0 & 1 & -56 & -3912 & 378 & 7476 & -1512 \cr
0 & 0 & 0 & 1 & 70 & -8 & -168 & 36 \cr
0 & 0 & 0 & 0 & 1 & -56 & -1560 & 420 \cr
0 & 0 & 0 & 0 & 0 & 1 & 28 & -8 \cr
0 & 0 & 0 & 0 & 0 & 0 & 1 & -8 \cr
0 & 0 & 0 & 0 & 0 & 0 & 0 & 1 \cr
           }
 \right] 
$$

$$
A :=  \left[ 
\matrix{                
0 & 0 & 0 & 0 & 0 & 1 & 0 & 0 \cr
0 & 0 & 0 & 1 & 0 & 8 & 0 & 0 \cr
0 & 1 & 0 & 8 & 0 & 28 & 0 & 0 \cr
1 & 8 & 0 & 28 & 0 & 56 & 0 & 0 \cr
0 & 0 & 1 & 56 & 0 & 70 & 0 & 0 \cr
0 & 0 & 0 & 0 & 1 & 56 & 0 & 0 \cr
0 & 0 & 0 & 0 & 0 & 0 & 1 & 0 \cr
0 & 0 & 0 & 0 & 0 & 0 & 0 & 1 \cr
           }
 \right] 
~~~
A^{\prime}:= \left[ 
\matrix{
1 & 0 & 0 & 0 & 0 & 0 & 0 & 0 \cr
0 & 1 & 0 & 0 & 0 & 0 & 0 & 0 \cr
0 & 0 & -56 & 1 & 0 & 0 & 0 & 0 \cr
0 & 0 & -70 & 0 & -56 & 1 & 0 & 0 \cr
0 & 0 & -56 & 0 & -28 & 0 & -8 & 1 \cr
0 & 0 & 28 & 0 & 8 & 0 & 1 & 0 \cr
0 & 0 & 8 & 0 & 1 & 0 & 0 & 0 \cr
0 & 0 & 1 & 0 & 0 & 0 & 0 & 0 \cr
         }
 \right] 
$$
$$
 S^{\beta} :=  \left[ 
\matrix{                
1 & 8 & 28 & 56 & 70 & 56 & 28 & -8 \cr
0 & 1 & 8 & 28 & 56 & 70 & 56 & -28 \cr
0 & 0 & 1 & 8 & 28 & 56 & 70 & -56 \cr
0 & 0 & 0 & 1 & 8 & 28 & 56 & -70 \cr
0 & 0 & 0 & 0 & 1 & 8 & 28 & -56 \cr
0 & 0 & 0 & 0 & 0 & 1 & 8 & -28 \cr
0 & 0 & 0 & 0 & 0 & 0 & 1 & -8 \cr
0 & 0 & 0 & 0 & 0 & 0 & 0 & 1 \cr
           }
 \right] ~~~
S^{\beta^{\prime}}:=  \left[ 
\matrix{
1 & 8 & 28 & 56 & 70 & -56 & -28 & -8 \cr
0 & 1 & 8 & 28 & 56 & -70 & -56 & -28 \cr
0 & 0 & 1 & 8 & 28 & -56 & -70 & -56 \cr
0 & 0 & 0 & 1 & 8 & -28 & -56 & -70 \cr
0 & 0 & 0 & 0 & 1 & -8 & -28 & -56 \cr
0 & 0 & 0 & 0 & 0 & 1 & 8 & 28 \cr
0 & 0 & 0 & 0 & 0 & 0 & 1 & 8 \cr
0 & 0 & 0 & 0 & 0 & 0 & 0 & 1 \cr
           }
 \right] 
$$
where $      \beta=  ( \beta_{34}    \beta_{23} \beta_{12}) \beta_{23}
(\beta_{56} 
\beta_{45} \beta_{34}\beta_{23} \beta_{12})$, $\beta^{\prime}=
[(\sigma_{34}\sigma_{56}\sigma_{78})
(\sigma_{45}\sigma_{67})\sigma_{56}]
[(\sigma_{67}\sigma_{78})\sigma_{67}]$.             

\vskip 0.2 cm                  
                   
                                        ${\bf CP}^8$

   $$
K_{92} :=  \left[ 
\matrix{                    
1 & 0 & 0 & 0 & 0 & 0 & 0 & 0 & 0 \cr
0 & 1 & 0 & 0 & 0 & 0 & 0 & 0 & 36 \cr
0 & 0 & 1 & 0 & 0 & 0 & 0 & 126 & 0 \cr
0 & 0 & 0 & 1 & 0 & 0 & 84 & 0 & 0 \cr
0 & 0 & 0 & 0 & 1 & 9 & 0 & 0 & 0 \cr
0 & 0 & 0 & 0 & 0 & 1 & 0 & 0 & 0 \cr
0 & 0 & 0 & 0 & 0 & 0 & 1 & 0 & 0 \cr
0 & 0 & 0 & 0 & 0 & 0 & 0 & 1 & 0 \cr
0 & 0 & 0 & 0 & 0 & 0 & 0 & 0 & 1 \cr
          }
 \right] 
~~~
 K_{93} :=  \left[ 
\matrix{                   
1 & 0 & 0 & 0 & 0 & 0 & 0 & 0 & 0 \cr
-9 & 1 & 0 & 0 & 0 & 0 & 0 & 0 & 0 \cr
0 & 0 & 1 & 0 & 0 & 0 & 0 & 0 & 84 \cr
0 & 0 & 0 & 1 & 0 & 0 & 0 & 126 & 0 \cr
0 & 0 & 0 & 0 & 1 & 0 & 36 & 0 & 0 \cr
0 & 0 & 0 & 0 & 0 & 1 & 0 & 0 & 0 \cr
0 & 0 & 0 & 0 & 0 & 0 & 1 & 0 & 0 \cr
0 & 0 & 0 & 0 & 0 & 0 & 0 & 1 & 0 \cr
0 & 0 & 0 & 0 & 0 & 0 & 0 & 0 & 1 \cr
           }
 \right] 
$$

$$
T :=  \left[ 
\matrix{                 
0 & 0 & 0 & 0 & 0 & 0 & 0 & 0 & 1 \cr
-1 & 0 & 0 & 0 & 0 & 0 & 0 & 36 & 9 \cr
0 & -1 & 0 & 0 & 0 & 0 & 126 & 84 & 0 \cr
0 & 0 & -1 & 0 & 0 & 84 & 126 & 0 & 0 \cr
0 & 0 & 0 & -1 & 9 & 36 & 0 & 0 & 0 \cr
0 & 0 & 0 & 0 & -1 & 0 & 0 & 0 & 0 \cr
0 & 0 & 0 & 0 & 0 & -1 & 0 & 0 & 0 \cr
0 & 0 & 0 & 0 & 0 & 0 & -1 & 0 & 0 \cr
0 & 0 & 0 & 0 & 0 & 0 & 0 & -1 & 0 \cr
           }
 \right] 
$$

$$
S :=  \left[ 
\matrix{                 
1 & 0 & 0 & 0 & 0 & 0 & 0 & 0 & 0 \cr
-9 & 1 & 0 & 0 & 0 & 0 & 0 & 0 & 36 \cr
45 & -9 & 1 & 0 & 0 & 0 & 0 & 126 & -240 \cr
-165 & 45 & -9 & 1 & 0 & 0 & 84 & -1008 & 990 \cr
495 & -165 & 45 & -9 & 1 & 9 & -720 & 4620 & -3168 \cr
3168 & -990 & 240 & -36 & 0 & 1 & -3015 & 25740 & -19932 \cr
-4620 & 1008 & -126 & 0 & 0 & 0 & 1 & -15867 & 25740 \cr
720 & -84 & 0 & 0 & 0 & 0 & 0 & 1 & -3015 \cr
-9 & 0 & 0 & 0 & 0 & 0 & 0 & 0 & 1 \cr
           }
 \right] 
$$

$$
P :=  \left[ 
\matrix{                   
0 & 0 & 0 & 0 & 1 & 0 & 0 & 0 & 0 \cr
0 & 0 & 0 & 0 & 0 & 1 & 0 & 0 & 0 \cr
0 & 0 & 0 & 1 & 0 & 0 & 0 & 0 & 0 \cr
0 & 0 & 0 & 0 & 0 & 0 & 1 & 0 & 0 \cr
0 & 0 & 1 & 0 & 0 & 0 & 0 & 0 & 0 \cr
0 & 0 & 0 & 0 & 0 & 0 & 0 & 1 & 0 \cr
0 & 1 & 0 & 0 & 0 & 0 & 0 & 0 & 0 \cr
0 & 0 & 0 & 0 & 0 & 0 & 0 & 0 & 1 \cr
1 & 0 & 0 & 0 & 0 & 0 & 0 & 0 & 0 \cr
           }
 \right] 
$$

$$
 S_{upper} :=  \left[ 
\matrix{                  
1 & 9 & -9 & -720 & 45 & 4620 & -165 & -3168 & 495 \cr
0 & 1 & -36 & -3015 & 240 & 25740 & -990 & -19932 & 3168 \cr
0 & 0 & 1 & 84 & -9 & -1008 & 45 & 990 & -165 \cr
0 & 0 & 0 & 1 & -126 & -15867 & 1008 & 25740 & -4620 \cr
0 & 0 & 0 & 0 & 1 & 126 & -9 & -240 & 45 \cr
0 & 0 & 0 & 0 & 0 & 1 & -84 & -3015 & 720 \cr
0 & 0 & 0 & 0 & 0 & 0 & 1 & 36 & -9 \cr
0 & 0 & 0 & 0 & 0 & 0 & 0 & 1 & -9 \cr
0 & 0 & 0 & 0 & 0 & 0 & 0 & 0 & 1 \cr
           }
 \right] 
$$

$$
A :=  \left[ 
\matrix{                 
0 & 0 & 0 & 0 & 0 & 0 & 1 & 0 & 0 \cr
0 & 0 & 0 & 0 & 1 & 0 & 9 & 0 & 0 \cr
0 & 0 & 1 & 0 & 9 & 0 & 36 & 0 & 0 \cr
1 & 0 & 9 & 0 & 36 & 0 & 84 & 0 & 0 \cr
0 & 1 & 36 & 0 & 84 & 0 & 126 & 0 & 0 \cr
0 & 0 & 0 & 1 & 126 & 0 & 126 & 0 & 0 \cr
0 & 0 & 0 & 0 & 0 & 1 & 84 & 0 & 0 \cr
0 & 0 & 0 & 0 & 0 & 0 & 0 & 1 & 0 \cr
0 & 0 & 0 & 0 & 0 & 0 & 0 & 0 & 1 \cr
           }
 \right] 
~~~
A^{\prime}:= \left[ \matrix{
0 & 1 & 0 & 0 & 0 & 0 & 0 & 0 & 0 \cr
1 & 9 & 0 & 0 & 0 & 0 & 0 & 0 & 0 \cr
0 & 0 & 1 & 0 & 0 & 0 & 0 & 0 & 0 \cr
0 & 0 & 0 & -126 & 1 & 0 & 0 & 0 & 0 \cr
0 & 0 & 0 & -126 & 0 & -84 & 1 & 0 & 0 \cr
0 & 0 & 0 & -84 & 0 & -36 & 0 & -9 & 1 \cr
0 & 0 & 0 & 36 & 0 & 9 & 0 & 1 & 0 \cr
0 & 0 & 0 & 9 & 0 & 1 & 0 & 0 & 0 \cr
0 & 0 & 0 & 1 & 0 & 0 & 0 & 0 & 0 \cr
}
 \right] 
$$
$$
     S^{\beta}=  \left[ 
\matrix{                 
1 & 9 & 36 & 84 & 126 & 126 & 84 & 36 & -9 \cr
0 & 1 & 9 & 36 & 84 & 126 & 126 & 84 & -36 \cr
0 & 0 & 1 & 9 & 36 & 84 & 126 & 126 & -84 \cr
0 & 0 & 0 & 1 & 9 & 36 & 84 & 126 & -126 \cr
0 & 0 & 0 & 0 & 1 & 9 & 36 & 84 & -126 \cr
0 & 0 & 0 & 0 & 0 & 1 & 9 & 36 & -84 \cr
0 & 0 & 0 & 0 & 0 & 0 & 1 & 9 & -36 \cr
0 & 0 & 0 & 0 & 0 & 0 & 0 & 1 & -9 \cr
0 & 0 & 0 & 0 & 0 & 0 & 0 & 0 & 1 \cr
           }
 \right] ~~~ 
S^{\beta^{\prime}}:= \left[ 
\matrix{
1 & 9 & 36 & 84 & 126 & 126 & -84 & -36 & -9 \cr
0 & 1 & 9 & 36 & 84 & 126 & -126 & -84 & -36 \cr
0 & 0 & 1 & 9 & 36 & 84 & -126 & -126 & -84 \cr
0 & 0 & 0 & 1 & 9 & 36 & -84 & -126 & -126 \cr
0 & 0 & 0 & 0 & 1 & 9 & -36 & -84 & -126 \cr
0 & 0 & 0 & 0 & 0 & 1 & -9 & -36 & -84 \cr
0 & 0 & 0 & 0 & 0 & 0 & 1 & 9 & 36 \cr
0 & 0 & 0 & 0 & 0 & 0 & 0 & 1 & 9 \cr
0 & 0 & 0 & 0 & 0 & 0 & 0 & 0 & 1 \cr
}
 \right] 
$$
where 
$\beta=
(\beta_{45}\beta_{34}\beta_{23}\beta_{12})(\beta_{34}\beta_{23})(\beta_{67}\beta_{56}
\beta_{45} \beta_{34} \beta_{23} \beta_{12})$, and
$\beta^{\prime}=\beta_{12} [(\sigma_{45}
\sigma_{67}\sigma_{89})(\sigma_{56}\sigma_{78}) \sigma_{67} ] [
(\sigma_{78} \sigma_{89})\sigma_{78}]$.

                                        ${\bf CP}^9$

$$
 K_{10,2} :=  \left[ 
\matrix{                  
1 & 0 & 0 & 0 & 0 & 0 & 0 & 0 & 0 & 0 \cr
0 & 1 & 0 & 0 & 0 & 0 & 0 & 0 & 0 & 45 \cr
0 & 0 & 1 & 0 & 0 & 0 & 0 & 0 & 210 & 0 \cr
0 & 0 & 0 & 1 & 0 & 0 & 0 & 210 & 0 & 0 \cr
0 & 0 & 0 & 0 & 1 & 0 & 45 & 0 & 0 & 0 \cr
0 & 0 & 0 & 0 & 0 & 1 & 0 & 0 & 0 & 0 \cr
0 & 0 & 0 & 0 & 0 & 0 & 1 & 0 & 0 & 0 \cr
0 & 0 & 0 & 0 & 0 & 0 & 0 & 1 & 0 & 0 \cr
0 & 0 & 0 & 0 & 0 & 0 & 0 & 0 & 1 & 0 \cr
0 & 0 & 0 & 0 & 0 & 0 & 0 & 0 & 0 & 1 \cr
           }
 \right] 
$$

$$
K_{10,3} :=  \left[ 
\matrix{                     
1 & 0 & 0 & 0 & 0 & 0 & 0 & 0 & 0 & 0 \cr
-10 & 1 & 0 & 0 & 0 & 0 & 0 & 0 & 0 & 0 \cr
0 & 0 & 1 & 0 & 0 & 0 & 0 & 0 & 0 & 120 \cr
0 & 0 & 0 & 1 & 0 & 0 & 0 & 0 & 252 & 0 \cr
0 & 0 & 0 & 0 & 1 & 0 & 0 & 120 & 0 & 0 \cr
0 & 0 & 0 & 0 & 0 & 1 & 10 & 0 & 0 & 0 \cr
0 & 0 & 0 & 0 & 0 & 0 & 1 & 0 & 0 & 0 \cr
0 & 0 & 0 & 0 & 0 & 0 & 0 & 1 & 0 & 0 \cr
0 & 0 & 0 & 0 & 0 & 0 & 0 & 0 & 1 & 0 \cr
0 & 0 & 0 & 0 & 0 & 0 & 0 & 0 & 0 & 1 \cr
           }
 \right] 
$$

$$
T :=  \left[ 
\matrix{                   
0 & 0 & 0 & 0 & 0 & 0 & 0 & 0 & 0 & 1 \cr
-1 & 0 & 0 & 0 & 0 & 0 & 0 & 0 & 45 & 10 \cr
0 & -1 & 0 & 0 & 0 & 0 & 0 & 210 & 120 & 0 \cr
0 & 0 & -1 & 0 & 0 & 0 & 210 & 252 & 0 & 0 \cr
0 & 0 & 0 & -1 & 0 & 45 & 120 & 0 & 0 & 0 \cr
0 & 0 & 0 & 0 & -1 & 10 & 0 & 0 & 0 & 0 \cr
0 & 0 & 0 & 0 & 0 & -1 & 0 & 0 & 0 & 0 \cr
0 & 0 & 0 & 0 & 0 & 0 & -1 & 0 & 0 & 0 \cr
0 & 0 & 0 & 0 & 0 & 0 & 0 & -1 & 0 & 0 \cr
0 & 0 & 0 & 0 & 0 & 0 & 0 & 0 & -1 & 0 \cr
           }
 \right] 
$$

$$
S :=  \left[ 
\matrix{                     
1 & 0 & 0 & 0 & 0 & 0 & 0 & 0 & 0 & 0 \cr
-10 & 1 & 0 & 0 & 0 & 0 & 0 & 0 & 0 & 45 \cr
55 & -10 & 1 & 0 & 0 & 0 & 0 & 0 & 210 & -330 \cr
-220 & 55 & -10 & 1 & 0 & 0 & 0 & 210 & -1848 & 1485 \cr
715 & -220 & 55 & -10 & 1 & 0 & 45 & -1980 & 9240 & -5148 \cr
-2002 & 715 & -220 & 55 & -10 & 1 & -440 & 10395 & -34320 & 15015
 \cr
17160 & -4752 & 990 & -120 & 0 & 0 & 1 & -25190 & 177705 & 
-120120 \cr
-11880 & 2310 & -252 & 0 & 0 & 0 & 0 & 1 & -52910 & 73755 \cr
1155 & -120 & 0 & 0 & 0 & 0 & 0 & 0 & 1 & -5390 \cr
-10 & 0 & 0 & 0 & 0 & 0 & 0 & 0 & 0 & 1 \cr
           }
 \right] 
$$

$$
P :=  \left[ 
\matrix{                    
0 & 0 & 0 & 0 & 0 & 1 & 0 & 0 & 0 & 0 \cr
0 & 0 & 0 & 0 & 1 & 0 & 0 & 0 & 0 & 0 \cr
0 & 0 & 0 & 0 & 0 & 0 & 1 & 0 & 0 & 0 \cr
0 & 0 & 0 & 1 & 0 & 0 & 0 & 0 & 0 & 0 \cr
0 & 0 & 0 & 0 & 0 & 0 & 0 & 1 & 0 & 0 \cr
0 & 0 & 1 & 0 & 0 & 0 & 0 & 0 & 0 & 0 \cr
0 & 0 & 0 & 0 & 0 & 0 & 0 & 0 & 1 & 0 \cr
0 & 1 & 0 & 0 & 0 & 0 & 0 & 0 & 0 & 0 \cr
0 & 0 & 0 & 0 & 0 & 0 & 0 & 0 & 0 & 1 \cr
1 & 0 & 0 & 0 & 0 & 0 & 0 & 0 & 0 & 0 \cr
           }
 \right] 
$$

$$
 S_{upper} :=  \left[ 
\matrix{                    
1 & -10 & -440 & 55 & 10395 & -220 & -34320 & 715 & 15015 & -2002
 \cr
0 & 1 & 45 & -10 & -1980 & 55 & 9240 & -220 & -5148 & 715 \cr
0 & 0 & 1 & -120 & -25190 & 990 & 177705 & -4752 & -120120 & 
17160 \cr
0 & 0 & 0 & 1 & 210 & -10 & -1848 & 55 & 1485 & -220 \cr
0 & 0 & 0 & 0 & 1 & -252 & -52910 & 2310 & 73755 & -11880 \cr
0 & 0 & 0 & 0 & 0 & 1 & 210 & -10 & -330 & 55 \cr
0 & 0 & 0 & 0 & 0 & 0 & 1 & -120 & -5390 & 1155 \cr
0 & 0 & 0 & 0 & 0 & 0 & 0 & 1 & 45 & -10 \cr
0 & 0 & 0 & 0 & 0 & 0 & 0 & 0 & 1 & -10 \cr
0 & 0 & 0 & 0 & 0 & 0 & 0 & 0 & 0 & 1 \cr
           }
 \right] 
$$

  $$
            A=\left[ \matrix{
0 & 0 & 0 & 0 & 0 & 0 & 0 & 1 & 0 & 0 \cr
0 & 0 & 0 & 0 & 0 & 1 & 0 & 10 & 0 & 0 \cr
0 & 0 & 0 & 1 & 0 & 10 & 0 & 45 & 0 & 0 \cr
0 & 1 & 0 & 10 & 0 & 45 & 0 & 120 & 0 & 0 \cr
1 & 10 & 0 & 45 & 0 & 120 & 0 & 210 & 0 & 0 \cr
0 & 0 & 1 & 120 & 0 & 210 & 0 & 252 & 0 & 0 \cr
0 & 0 & 0 & 0 & 1 & 252 & 0 & 210 & 0 & 0 \cr
0 & 0 & 0 & 0 & 0 & 0 & 1 & 120 & 0 & 0 \cr
0 & 0 & 0 & 0 & 0 & 0 & 0 & 0 & 1 & 0 \cr
0 & 0 & 0 & 0 & 0 & 0 & 0 & 0 & 0 & 1 \cr
           }
 \right] 
$$
$$
A^{\prime}  :=  \left[ 
\matrix{
1 & 0 & 0 & 0 & 0 & 0 & 0 & 0 & 0 & 0 \cr
0 & 1 & 0 & 0 & 0 & 0 & 0 & 0 & 0 & 0 \cr
0 & 0 & -120 & 1 & 0 & 0 & 0 & 0 & 0 & 0 \cr
0 & 0 & -210 & 0 & -252 & 1 & 0 & 0 & 0 & 0 \cr
0 & 0 & -252 & 0 & -210 & 0 & -120 & 1 & 0 & 0 \cr
0 & 0 & -210 & 0 & -120 & 0 & -45 & 0 & -10 & 1 \cr
0 & 0 & 120 & 0 & 45 & 0 & 10 & 0 & 1 & 0 \cr
0 & 0 & 45 & 0 & 10 & 0 & 1 & 0 & 0 & 0 \cr
0 & 0 & 10 & 0 & 1 & 0 & 0 & 0 & 0 & 0 \cr
0 & 0 & 1 & 0 & 0 & 0 & 0 & 0 & 0 & 0 \cr
}
 \right]                   
$$

$$
 S^{\beta} :=  \left[ 
\matrix{                  
1 & 10 & 45 & 120 & 210 & 252 & 210 & 120 & 45 & -10 \cr
0 & 1 & 10 & 45 & 120 & 210 & 252 & 210 & 120 & -45 \cr
0 & 0 & 1 & 10 & 45 & 120 & 210 & 252 & 210 & -120 \cr
0 & 0 & 0 & 1 & 10 & 45 & 120 & 210 & 252 & -210 \cr
0 & 0 & 0 & 0 & 1 & 10 & 45 & 120 & 210 & -252 \cr
0 & 0 & 0 & 0 & 0 & 1 & 10 & 45 & 120 & -210 \cr
0 & 0 & 0 & 0 & 0 & 0 & 1 & 10 & 45 & -120 \cr
0 & 0 & 0 & 0 & 0 & 0 & 0 & 1 & 10 & -45 \cr
0 & 0 & 0 & 0 & 0 & 0 & 0 & 0 & 1 & -10 \cr
0 & 0 & 0 & 0 & 0 & 0 & 0 & 0 & 0 & 1 \cr
           }
 \right] 
$$
$$ 
S^{\beta^{\prime}} :=  \left[ 
\matrix{
1 & 10 & 45 & 120 & 210 & 252 & -210 & -120 & -45 & -10 \cr
0 & 1 & 10 & 45 & 120 & 210 & -252 & -210 & -120 & -45 \cr
0 & 0 & 1 & 10 & 45 & 120 & -210 & -252 & -210 & -120 \cr
0 & 0 & 0 & 1 & 10 & 45 & -120 & -210 & -252 & -210 \cr
0 & 0 & 0 & 0 & 1 & 10 & -45 & -120 & -210 & -252 \cr
0 & 0 & 0 & 0 & 0 & 1 & -10 & -45 & -120 & -210 \cr
0 & 0 & 0 & 0 & 0 & 0 & 1 & 10 & 45 & 120 \cr
0 & 0 & 0 & 0 & 0 & 0 & 0 & 1 & 10 & 45 \cr
0 & 0 & 0 & 0 & 0 & 0 & 0 & 0 & 1 & 10 \cr
0 & 0 & 0 & 0 & 0 & 0 & 0 & 0 & 0 & 1 \cr
}
 \right] $$

where $\beta=
(\beta_{56}\beta_{45}\beta_{34}\beta_{23}\beta_{12})(\beta_{45}\beta_{34}
\beta_{23})\beta_{34}
(\beta_{78}\beta_{67}\beta_{56}\beta_{45}\beta_{34}\beta_{23}
\beta_{12})
$,

\noindent
 and
$\beta^{\prime}=[(\sigma_{34}\sigma_{56}\sigma_{78}\sigma_{9,10})  
(\sigma_{45}\sigma_{67}\sigma_{89})(\sigma_{56}\sigma_{78})\sigma_{67}]
[(\sigma_{78}\sigma_{89}
\sigma_{9,10})(\sigma_{78}\sigma_{89})\sigma_{78}]$. 
                  
}


%
%
%

%
%
%

\end{document}